\documentclass[10pt, reqno]{amsart}

\usepackage{amssymb}
\usepackage[pdftex]{graphicx}
\usepackage{epstopdf}
\usepackage{float}
\usepackage{amsmath}
\usepackage[a4paper,margin=1.25in,footskip=0.25in]{geometry}
\usepackage{amsthm}
\usepackage[all]{xy}
\usepackage{amsfonts}
\usepackage{txfonts}
\usepackage{mathrsfs}
\usepackage{mathdots}
\usepackage{mathpazo}

\usepackage{latexsym}
\usepackage[shortlabels]{enumitem}
\usepackage{hyperref}
\usepackage{tikz-cd}
\newtheorem{theorem}{Theorem}[section]
\newtheorem{lemma}[theorem]{Lemma}
\newtheorem{proposition}[theorem]{Proposition}

\theoremstyle{definition}
\newtheorem{definition}[theorem]{Definition}

\theoremstyle{remark}
\newtheorem{remark}[theorem]{Remark}

\numberwithin{equation}{section}

\makeatletter
\g@addto@macro\@floatboxreset\centering
\makeatother

\title{A Gross-Kohnen-Zagier type formula for moduli of shtukas with Iwahori level structures}
\author{Hao Li}
\date{April 22, 2019.}
\begin{document}

\begin{abstract}
In this paper I'm going to study the intersection of two Heegner-Drinfeld cycles coming from two different nonsplit tori on the Yun-Zhang moduli stack of $PGL_2$ Drinfeld stukas with Iwahori level structure. We will see that the intersection number is related to a certain period integral. It is an extension of the result by Howard-Shinidman to the Iwahori case. 
\end{abstract}
\maketitle
\section{Introduction}
In their secomd volumn\cite{YZ2}, Yun-Zhang generalized their previous work\cite{YZ1} to the moduli stack of shtukas with Iwahori level structure and Heegner-Drinfeld cycles coming from a double cover with ramification points away from the level. A natural question is whether one can do a similar thing for Heegner-Drinfeld cycles coming from two different double covers, i.e. relating the intersection number of two Heegner-Drinfeld cycles attached to two different double covers (or rather their certain Hecke eigen-parts) with automorphic L-functions with Iwahori levels, like what Howard and Shnidman did previously in the no level structure case\cite{HS}. This case can be viewed as a function field analogue of the Gross-Kohnen-Zagier formula\cite{GKZ}.

In this paper I will give the answer to this question in the case where the two double covers are still everywhere nonramified. Let $p$ be a prime number greater than $2$, $k$ a finite field of characteristic $p$ with cardinality $q$ (i.e. $\mathbb{F}_q$). Let $X$ be a smooth geometrically connected curve over $k$ of genus $g$ with $F$ its field of rational functions. Take two everywhere nonramified double covers $Y_1$ and $Y_2$ of $X$, and let $Y=Y_1\times_X Y_2$ be the associated fourfold cover $X$. Then there exists a unique third double cover over $X$, $Y_3$, below $Y$ different from $Y_1$ and $Y_2$, and they fit into the following diagram:
\begin{equation}
\xymatrix{ & Y \ar[ld]_{\pi_1} \ar[d]^{\pi_3} \ar[rd]^{\pi_2} & \\
Y_1 \ar[rd]_{\nu_1} & Y_3 \ar[d]^{\nu_3} & Y_2\ar[ld]^{\nu_2}\\
& X &
}
\end{equation}
The generic fiber of this diagram gives the field extension diagram:
\begin{equation}
\xymatrix{ & K \ar[ld]_{\pi_1} \ar[d]^{\pi_3} \ar[rd]^{\pi_2} & \\
K_1 \ar[rd]_{\nu_1} & K_3 \ar[d]^{\nu_3} & K_2\ar[ld]^{\nu_2}\\
& F &
}
\end{equation}
Let $\mathbb{A}$ be the adelic ring of $F$ and $\mathbb{A}_i$ be the adelic ring of $K_i$ for $i=1,2,3$. Also let $\widetilde{T}_i=Res_{Y_i/X}\mathbb{G}_{m,Y_i}$ and $T_i=\widetilde{T}_i/\mathbb{G}_{m,X}$. Then one has $\widetilde{T}_i(F)=K_i^\times$ and $T_i(F)=K_i^\times/F^\times$.

Let $\Sigma\subset |X|$ be a finite set of physical points of $X$ of total degree $N$, this will be the ``level". Like what Yun and Zhang did in their second volumn, one needs to consider the splitting behavior of points in $\Sigma$ in the double covers. In order to make the question under consideration here meaningful, I require the following: there exist a partition of $\Sigma$:
\begin{equation*}
\Sigma=\Sigma_f\sqcup\Sigma_\infty
\end{equation*}
such that: Any $x\in\Sigma_f$ splits in both $Y_1$ and $Y_2$, and any $x\in\Sigma_\infty$ is inert in both the double covers. Let $\mathfrak{S}_\infty=\prod_{x\in\Sigma}\text{Spec}k_x$ and $G=PGL_{2,X}$. Yun-Zhang defined the moduli of shtukas with Iwahori level structures $Sht_G^r(\Sigma;\Sigma_\infty)$. It is a smooth Deligne-Mumford stack of dimension $2r$ equipped with a map:
\begin{equation*}
\Pi^r_G: Sht_G^r(\Sigma;\Sigma_\infty)\longrightarrow X^r\times_k\mathfrak{S}_\infty
\end{equation*}
To define the Heegner-Drinfeld cycle for $Y_1, Y_2$, one needs something similar to $\mathfrak{S}_\infty$: $\mathfrak{S}'_{1,\infty}, \mathfrak{S}'_{2,\infty}$ and $\widetilde{\mathfrak{S}}_\infty$. They are all zero dimensional schemes over $k$ and fit into the following diagram:
\begin{equation*}
\xymatrix{ & \widetilde{\mathfrak{S}}_\infty \ar[ld] \ar[rd] &\\
\mathfrak{S}'_{1,\infty} \ar[rd] & & \mathfrak{S}'_{2,\infty} \ar[ld]\\
& \mathfrak{S}_\infty & 
}
\end{equation*}

I refer you to section 2.6 for the precise meaning of them because $\widetilde{\mathfrak{S}}_\infty$ needs some explaination. Also, one needs a choice of auxilary data: $\underline{\mu}:=(\mu, \mu_{f}, \mu_\infty)\in\{\pm 1\}^r\times\{\pm 1\}^{\Sigma_f}\times\{\pm 1\}^{\Sigma_\infty}$. With its help , one can define a smooth Deligne-Mumford stack of dimesion $r$, $Sht^{\underline{\mu}}_{T_i}(\mu_\infty\cdot\Sigma'_{i,\infty})$ equipped with a map:
\begin{equation*}
Sht^{\underline{\mu}}_{T_i}(\mu_\infty\cdot\Sigma'_{i,\infty})\longrightarrow Y^r_i\times_k\mathfrak{S}'_{i,\infty}
\end{equation*}
For $i=1,2$. For some technical reasons, I will consider the following intersection problem: Let $\xi\in\widetilde{\mathfrak{S}}_\infty(\overline{k})$, consider the following morphism:
\begin{equation*}
Sht^{\underline{\mu}}_{T_i}(\mu_\infty\cdot\Sigma'_{i,\infty})\times_{\mathfrak{S}'_{i,\infty}}\xi\longrightarrow Sht_G^r(\Sigma;\Sigma_\infty)\times_{\mathfrak{S}_\infty}\xi
\end{equation*}
Let $\mathcal{Z}_{i}^\mu(\xi)$ be the push foward of the fundamental cycle of the left hand side along this morphism.

Yun-Zhang proved the following ``coarse spectral decomposition" of the total cohomology of $Sht_G^r(\Sigma;\Sigma_\infty)$ (\cite{YZ2}, Theorem 3.41) (By total cohomology I mean not the cohomology of the generic fiber): Let $V(\xi)=H_c^r(Sht_G^r(\Sigma;\Sigma_\infty)\times_{\mathfrak{S}_\infty}\xi, \mathbb{Q}_l)$, then one has the decompostion:
\begin{equation*}
V(\xi)\otimes\overline{\mathbb{Q}}_l=V(\xi)_{\pi}\oplus (V_{Eis}(\xi)\otimes\overline{\mathbb{Q}}_l)
\end{equation*}
In which $\pi$ runs through all the automorphic representations of $PGL_{2,F}$ with Iwahori level at $\Sigma$. Let $Z_{i,\pi}^\mu$ be the projection to the $\pi$-th part of the cycle class of $\mathcal{Z}_{i}^\mu$.

Following Howard-Shnidman, let $\widetilde{G}_3=GL_{2,\nu_*\mathcal{O}_{Y_3}}$ be the $X$-group scheme whose functor of points on an $X$-scheme $f: U\longrightarrow X$ is given by $H^0(U, End_{\mathcal{O}_U}(f^*(\nu_*\mathcal{O}_{Y_3}))$ and $G_3=\widetilde{G}_3/Z(\widetilde{G}_3)$. Then $T_3$ is canonically sitting inside $G_3$. It's generic fiber is isomorphic to the generic fiber of $PGL_{2,X}$ though they are different $X$-schemes. The choice of auxilary data $(\mu_f, \mu_\infty)$ above determines an isomorphism:
\begin{equation*}
G_3(\mathbb{A})\xrightarrow{\sim}G(\mathbb{A})
\end{equation*}

Let $\pi$ be an automorphic representation of $G(\mathbb{A})$ and $\phi\in\pi$. Using the above isomorphism, one can view $\phi$ as an automorphic form on $G_3(\mathbb{A})$, denoted by $\phi_3$. The double cover $Y\xrightarrow{\pi_3}Y_3$ determines a quadratic character $\eta$ of $T_3(\mathbb{A})=\mathbb{A}^\times_3/\mathbb{A}^\times$. Then it is legitimate to define the following global period integral of $\phi$: 
\begin{equation*}
\mathscr{P}_0(\phi, s)=\int_{[A]}\phi(t_0)|t_0|^sdt_0\quad \mathscr{P}_3(\phi, s)=\int_{[T_3]}\phi_3(t_3)\eta(t_3)dt_3
\end{equation*}
In which $A$ is the standard diagonal torus of $G$, $[A]=A(F)\backslash A(\mathbb{A})$ and $[T_3]=T_3(F)\backslash T_3(\mathbb{A})$. Then the first main result is the following:
\begin{theorem}
Fix the choice of $\underline{\mu}=(\mu,\mu_f,\mu_\infty)$ and $\xi$ and $N=deg\Sigma$ as above. Let $\pi$ be an automorphic represetation of $G(\mathbb{A})$ with Iwahori level at $\Sigma$. Let $\phi\in\pi$ be the vector whose local components are hyperspecial fixed away from $\Sigma$ and Iwahori fixed at $\Sigma$. Then one has:
\begin{equation*}
(Z^\mu_{1,\pi}(\xi), Z^\mu_{2,\pi}(\xi))=(logq)^{-r}(\displaystyle\frac{d}{ds})^r|_{s=0}(q^{Ns}\frac{\mathscr{P}_0(\phi,s)\mathscr{P}_{3,\eta}(\overline{\phi_3})}{\langle \phi, \phi \rangle})
\end{equation*}
\end{theorem}

Let $\chi_1$ and $\chi_2$ be quadratic characters of $\mathbb{A}^\times$ determined by the double covers $Y_1$ and $Y_2$ respectively. Let $L(\pi\otimes\chi_1,s)$ and $L(\pi\otimes\chi_2,s)$ be the twisted L-functions. Also, one defines $\mathscr{L}(\pi, s)=q^{2g-2+N}L(\pi, s)$ to be the normalized L-function. The period integral is actually related to the product of these L-functions, as in the case worked out by Howard-Shnidman. This leads to the following corollary.
\begin{theorem}
In the notations as above, one has:
\begin{equation*}
(Z^\mu_{1,\pi}(\xi), Z^\mu_{1,\pi}(\xi))=0
\end{equation*}
is equivalent to:
\begin{equation*}
\mathscr{L}^{(r)}(\pi, 1/2)\cdot L(\pi\otimes\chi_1, 1/2)\cdot L(\pi\otimes\chi_2, 1/2)=0
\end{equation*}
\end{theorem}

Also, in their paper, Howard and Shnidman computed the plain pairing of the two cycles, i.e. without projecting to the $\pi$-th eigenpart (\cite{HS}, Theorem $\bf{C}$). Here one also has a slightly different result:
\begin{theorem}
Again, with the notations above,
\begin{equation*}
(Z^\mu_1(\xi),Z^\mu_2(\xi))=0
\end{equation*}
\end{theorem}

For the precise meaning of those symbols, please go to the main body of the paper.

The method of proof basically follows the philosophy of Yun-Zhang, with some technical ingredients from Howard-Shnidman. In the second section I'm going to review the moduli stack of shtukas with Iwahori level structure defined by Yun and Zhang, and introduce the Heegner-Drinfeld cycles to be used in this paper. The third section is devoted to the construction of the auxilary moduli spaces, which Yun calls ``the Hitchin type moduli spaces"\cite{Y}. In the fourth section the intersection pairing is computed by applying the Grothendieck-Lefschetz trace formula to the relative cohomology of the ``Hitchin type fiberation". Then I'll move to the analytic side of the picture, relating the geometric side of the relative trace formula with  the relative cohomology of a certain local system for the ``dual Hitchin type fiberation". Then the main result would follow from comparing the cohomological side of the two Grothendieck-Lefschetz formulae. I'm also going to prove the other two theorems in the last section.

The reason I still work with everywhere nonramified covers is that I will use descent along torsors in section 3. In the more general case in which the covers are ramified at some points away from the level, I hope the more general fppf-descent can work. I'll probably return to this topic later. 

\section{The shtuka side set up}
In this section I review Yun-Zhang's definition of moduli of shtukas with Iwahori level structure and define the Heegner-Drinfeld cycles with which I'm going to work later. 

\subsection{Splitting behavior of $\Sigma$ in $Y_3$}
Recall $Y_1, Y_2$ and $Y_3$ are the three everywhere nonramified double covers in the introduction. $Y$ is the associated fourfold cover and $\Sigma$ be the set of closed point in $|X|$ satisfying following condition:
\begin{equation*}
\Sigma=\Sigma_f\sqcup\Sigma_\infty
\end{equation*}
and any $x\in\Sigma_f$ splits in both $Y_1$ and $Y_2$, and any $x\in\Sigma_\infty$ is inert in both the double covers. One also needs to know the splitting behavior of each $x\in\Sigma$ in $Y_3$. Actually , by an easy exercise of class field theory, one has the following:
\begin{proposition}
All the points in $\Sigma$ split in $Y_3$.
\end{proposition}
Then we have the following complete description of the splitting behavior of $x\in\Sigma$ in all the covers:
\begin{itemize}
\item Each $x\in\Sigma_f$ has two physical points (just prime ideals) over it in each $Y_i$ for $i=1,2,3$, and four points over it in the top cover $Y$. Let $y_x, y'_x$ be the two points over $x$ in $Y_1$, $z_x, z_x'$ be the two points over $x$ in $Y_2$. Then one can denote the four preimages of $x$ in $Y$ by $(y_x, z_x), (y'_x, z_x), (y_x, z'_x), (y'_x, z'_x)$. Then let $w_x$ be the point in $Y_3$ below $(y_x, z_x)$ and $(y_x', z'_x)$, $w'_x$ be the point below $(y'_x, z_x), (y_x, z'_x)$.
\item Each point $x\in\Sigma_\infty$ has one physical point with residue extension of degree $2$ in each $Y_i, i=1,2$. Each of these points over $x$ is split in the top cover, so two physical points above. $x$ has two physical points above it in $Y_3$, and these two points are both inert in $Y$.
\end{itemize}

\subsection{Review of the moduli of G-torsor over a curve with Iwahori level structures}
In their second volume, Yun and Zhang defined the moduli stack of $G$-torsors over $X$ with Iwahori level structure at $\Sigma$ (\cite{YZ2}, section 3.1). First, let $Bun_2(\Sigma)$ be the functor from the category of $k$-schemes to the category of groupoids whose $S$ points are the following:
\begin{equation*}
\mathcal{E}^{\dagger}=(\mathcal{E}, \{\mathcal{E}(-\frac{1}{2})\}_{x\in\Sigma})
\end{equation*}
in which:
\begin{itemize}
\item $\mathcal{E}$ is a rank $2$ vector bundle over the product $X\times_kS$.
\item For each $x\in\Sigma$, a rank $2$ vector bundle $\mathcal{E}(-\frac{1}{2}x)$ who fits into the following chain of coherent sheaves:
\begin{equation*}
\mathcal{E}\supset\mathcal{E}(-\frac{1}{2}x)\supset\mathcal{E}(-x):=\mathcal{E}\otimes\mathcal{O}_X(-x)
\end{equation*}
\end{itemize}
such that the quotient sheaf $\mathcal{E}/\mathcal{E}(-\frac{1}{2}x)$ is a skycrayper sheaf supported on the subscheme $\{x\}\times S$.

Let $Pic_X$ be the Picard stack of $X$, i.e. the stack whose $S$ points are line bundles over $X\times_kS$. $Pic_X$ acts $Bun_2(\Sigma)$ by simultenously twisting $\mathcal{E}$ and all $\mathcal{E}(-\frac{1}{2}x)$. Then $Bun_G(\Sigma)$ is defined to be the quotient:
\begin{equation*}
Bun_2(\Sigma)/Pic_X
\end{equation*}
Of course, they defined these stacks for general $GL_n$ and $PGL_n$, and their definition can be easily extended to arbitrary parahoric level structures, but here I only work with $n=2$, so I just review this simpler case.

Actually the dimension of $Bun_G(\Sigma)$ is well known:
\begin{proposition}
$Bun_G(\Sigma)$ is an Artin stack of dimension $3(g-1)+N$
\end{proposition}

One already knows that the plain $Bun_G$ is of dimension is $3(g-1)$, the Iwahori level structure just "adds" $N$ copies of $\mathbb{P}^1_k$ to it, hence the above number.

\subsection{Fractional twists and Atkin-Lehner involution} To define the moduli of shtukas with level structures, Yun-Zhang introduced fractional twists and Atkin-Lehner involution for vector bundles with Iwahori level structures (\cite{YZ2}, Definition 3.2). Again, I only review what I need here, for the more general cases, you can just go to read their paper.

First we define:
\begin{equation*}
\mathfrak{S}_\infty=\prod_{x\in\Sigma}\text{Spec}k_x
\end{equation*}
An $S$ point of $\mathfrak{S}_\infty$ is simply a set of maps: $x^{(1)}\longrightarrow \text{Spec}k_x$ for each $x\in\Sigma_\infty$. Let $x^{(i+1)}=x^{(i)}\circ Fr_S$. Let $d_x=[k_x:k]$, we get $d_x$ maps from $S$ to $\text{Spec}k_x$. 

As in their paper, one has the canonical point:
\begin{equation*}
\textbf{x}^{(1)}:\mathfrak{S}_\infty\longrightarrow \text{Spec}k_x\longrightarrow X
\end{equation*}
for each $x\in\Sigma_\infty$. It is simply projection followed by injection. Let $\textbf{x}^{(2)}=\textbf{x}^{(1)}\circ Fr_{\mathfrak{S}_\infty}$. Let $\Gamma_{\textbf{x}^{(i)}}$ be the graph of $\textbf{x}^{(i)}$. Then we have the following:
\begin{equation*}
\text{Spec}k_x\times_k\mathfrak{S}_\infty=\coprod_{i=1}^{d_x}\Gamma_{\textbf{x}^{(i)}}
\end{equation*}
Take an $S$ point $\mathcal{E}^\dagger$, of $Bun_2(\Sigma)\times\mathfrak{S}_\infty$. Then $S$ acquires a $k_x$ structure for each $x\in\Sigma_\infty$. Therefore we have:
\begin{equation*}
\text{Spec}k_x\times_k S=\coprod_{i=1}^{d_x}S^{(i)}
\end{equation*}
In which $S^{(i)}$ is the image of $x^{(i)}$.

By the definition of $\mathcal{E}^\dagger$, $\mathcal{E}/\mathcal{E}(-\frac{1}{2}x)$ is supported on these $d_x$ copies of $S$s, therefore splits into a direct sum $\oplus_{j=1}^{d_x}\mathscr{L}_i^{(j)}$, where $\mathscr{L}^{(j)}_i$ is supported on the $j$-th copy of $S$, $S^{(j)}$.

Let $D ‎=‎ \sum_{x\in\Sigma_\infty}\sum_{1\leq j\leq d_x}c_x^{(j)}\textbf{x}^{(1)}‎$, In which all $c_x^{(j)}\in\frac{1}{2}\mathbb{Z}$. Yun and Zhang defined the fractional twist by $D$ of $\mathcal{E}^\dagger$ as follows:
First define:
\begin{equation*}
\mathcal{E}(-\frac{1}{2}x^{(j)}):=ker\big(\mathcal{E}\longrightarrow\mathcal{E}/\mathcal{E}(-\frac{1}{2}x)\longrightarrow\mathscr{L}^{(j)}_i\big)
\end{equation*}

Then for any $D$ as above, first rewrite it as $D=D_0-D_1$, in which $D_0$ is an integral divisor on $X\times_k\mathfrak{S}_\infty$ and $D_1=\sum_{x\in\Sigma_\infty}\sum_{1\leq j\leq d_x}\frac{i^{(j)}_x}{2}x^{(j)}\textbf{x}^{(1)}‎$ in which each $i^{(j)}_x$ is either $0$ or $1$. Define $\mathcal{E}(-D_1)$ as the kernel of the following direct sum of projections:
\begin{equation*}
\mathcal{E}\longrightarrow\big(\oplus_{x\in\Sigma_{\infty}, 1\leq j\leq d_x}\mathcal{E}/\mathcal{E}(-\frac{i^{(j)}_x}{2})\big)
\end{equation*}
Then define $\mathcal{E}(D):=\mathcal{E}(-D_1)\otimes_{\mathcal{O}_{X\times_k\mathfrak{S}_\infty}}\mathcal{O}_{\mathcal{O}_{X\times_k\mathfrak{S}_\infty}}$.
For $\mathcal{E}(-\frac{1}{2}x)$, apply the same operation to the chain:
\begin{equation*}
\mathcal{E}(-\frac{1}{2}x)\supset\mathcal{E}(-x)
\end{equation*}
One can get the corresponding $\mathcal{E}(D)(-\frac{1}{2}x)$. This process defines the following map:
\begin{equation}
\widetilde{AL}(D): Bun_2(\Sigma)\times_k\mathfrak{S}_\infty\longrightarrow Bun_2(\Sigma)
\end{equation}
quotient out the simultenous twisting by $Pic_X$, one has the follow:
\begin{equation*}
AL(D): Bun_G(\Sigma)\times_k\mathfrak{S}_\infty\longrightarrow Bun_G(\Sigma)
\end{equation*}
This is the Atkin-Lehner involution of $Bun_G(\Sigma)$.

As in their paper, for the moduli of shtukas to be used later, one only cares about the case $D^{(1)}_\infty=\sum_{x\in\Sigma_\infty}\frac{1}{2}\textbf{x}^{(1)}$, and the Atkin-Lehner involution by $-D^{(1)}_\infty$. Again let $\mathcal{E}^\dagger$ be an $S$ point of $Bun_2(\Sigma)$. Then we simply have:
\begin{equation*}
\mathcal{E}(-D^{(1)}_\infty)=ker\big(\mathcal{E}\longrightarrow\mathcal{E}(-\frac{1}{2}x^{(1)}))
\end{equation*}
i.e.The sections of $\mathcal{E}$ who is belong to $\mathcal{E}(-\frac{1}{2})$ on the first copy of $S$, $S^{(1)}$.

\subsection{The Hecke stack} 
To define moduli of shtukas, one also has to review the Hecke stack for $Bun_G(\Sigma)$ (\cite{YZ2}, Definition 3.3).

Let $r$ be a nonnegative integer, and let $\underline{\mu}=(\mu_1, ..., \mu_r)\in\{\pm 1\}^r$. 
\begin{definition}
Let $Hk^{\underline{\mu}}_2(\Sigma)$ be the stack over $k$ whose $S$-points is the following:
\begin{itemize}
\item A sequence of $r+1$ rank $2$ vector bundles with Iwahori level structure: $\mathcal{E}_i^\dagger=(\mathcal{E}_i, \{\mathcal{E}_i(-\frac{1}{2}x)\}_{x\in\Sigma})\in Bun_G(\Sigma)(S)$ in which $i=0,1,...,r$.
\item $r$ morphisms $x_i: S\longrightarrow X$ for $i=1,...,r$, with the graph $\Gamma_{x_i}\subset X\times_k S$.
\item Isomorphism of vector bundles 
$$f_i: \mathcal{E}_{i-1}|_{X\times_kS-\Gamma_{x_i}}\xrightarrow{\sim}\mathcal{E}_i|_{X\times_kS-\Gamma_{x_i}}$$
for $i=1,...,r$, such that the induced relative position on the formal neighbourhood of $\Gamma_{x_i}$ is $\mu_i$, and $f_i$ respects the filtration given by $\mathcal{E}_{i-1}(-\frac{1}{2}x)$ and $\mathcal{E}_{i}(-\frac{1}{2}x)$ for all $x\in\Sigma$.
\end{itemize}
Let $Hk^{\underline{\mu}}_G(\Sigma):=Hk^{\underline{\mu}}_2(\Sigma)/Pic_X$.
\end{definition}
From the definition, one has the projections from $Hk_2^{\underline{\mu}}$ to $X^r$ recording the points of modifications $x_i$ (people call them pattes, recently the term ``legs" becomes popular). Also, to record the $i$th bundle with Iwahori level structure $\mathcal{E}_i^\dagger$, one has the following projections:
\begin{equation*}
\widetilde{p}_i: Hk_2^{\underline{\mu}}(\Sigma)\longrightarrow Bun_2(\Sigma)
\end{equation*}
for $i=0,...,r$.
Again, quotient out the simultanous twist by $Pic_X$, we get the following:
\begin{equation*}
p_i: Hk_G^{\underline{\mu}}(\Sigma)\longrightarrow Bun_G(\Sigma)
\end{equation*}
Yun and Zhang proved the following geometric properties of the Hecke stacks (\cite{YZ2}, Proposition 3.4).
\begin{proposition}
\begin{enumerate}
\item For $0\leq i \leq r$, the projection map $\widetilde{p}_i: Hk_2^{\underline{\mu}}(\Sigma)\longrightarrow Bun_2(\Sigma)$ is smooth of relative dimension $2r$.
\item For $0\leq i \leq r$, the morphism $(\widetilde{p}_i, \pi^{\underline{\mu}}_Hk): Hk_2^{\underline{\mu}}(\Sigma)\longrightarrow Bun_2(\Sigma)\times_kX^r$ is smooth of relative dimension $r$ when restricted to $Bun_2(\Sigma)\times_k(X-\Sigma)^r$.
\item For $0\leq i \leq r$, the morphism $(\widetilde{p}_i, \pi^{\underline{\mu}}_Hk): Hk_2^{\underline{\mu}}(\Sigma)\longrightarrow Bun_2(\Sigma)\times_kX^r$ is flat of relative dimension $r$.
\item All of the above hold for $Hk_G^{\underline{\mu}}(\Sigma)$ and $Bun_G(\Sigma)$.
\end{enumerate}
\end{proposition}
These properties are needed later.

\subsection{The moduli of shtukas}
Now we are ready to review the definition of the moduli of shtukas (\cite{YZ2}, section 3.2). First define the moduli of shtukas for $GL_2$. Let $\mathscr{D}_\infty=\{\sum_{x\in\Sigma_\infty, 1\leq j\leq d_x}c^{(i)}_x\textbf{x}^{(1)}: c^{(j)}_i\in\frac{1}{2}\mathbb{Z}\}$. Take a $\underline{\mu}\in\{\pm 1\}^r$. To define the moduli of stukas attached to $D_\infty\in\mathscr{D}_\infty$ and $\underline{\mu}$, one furthur requires the following condition:
\begin{equation}
\sum_{i=1}^{r}\mu_i=\sum_{x\in\Sigma_\infty, 1\leq j\leq d_x}c^i_x=2\text{deg}D_\infty
\end{equation}
Once one sees the following definition of the moduli of shtukas, it is clear why there should be such a relation. 
\begin{definition}
Define the stack $Sht^{\underline{\mu}}_2(\Sigma;D_\infty)$ as the following Cartesian diagram:
\begin{equation}
\xymatrix{Sht^{\underline{\mu}}_2(\Sigma;D_\infty) \ar[r] \ar[d]^{\omega_0} & Hk_2^{\underline{\mu}}(\Sigma)\times_k\mathfrak{S}_\infty \ar[d]^{(\widetilde{p}_0, AL(-D_\infty)\circ (\widetilde{p}_r\times id_{\mathfrak{S}_\infty}))} \\
Bun_2(\Sigma) \ar[r]^-{(id, Fr)} & Bun_2(\Sigma)\times Bun_2(\Sigma)}
\end{equation}
\end{definition}
More precisely, an $S$ point of $Sht^{\underline{\mu}}_2(\Sigma;D_\infty)$ consists of the following data:
\begin{itemize}
\item A map: $S\longrightarrow\mathfrak{S}_\infty$, which gives $x^{(1)}: S\longrightarrow Speck_x$ for each $x\in\Sigma_\infty$;
\item $r$ morphisms $x_i: S\longrightarrow X$ for $i=1,...,r$, with the graph $\Gamma_{x_i}\subset X\times_k S$.
\item A sequence of modifications, which starts at $\mathcal{E}_0^\dagger$ and "ends at" $^{\tau}\mathcal{E}_0(D_\infty)$:
\[
\xymatrix{\mathcal{E}_0 \ar@{-->}[r]^{f_1} & \mathcal{E}_1 \ar@{-->}[r]^{f_2} &......\ar@{-->}[r]^{f_r} & \mathcal{E}_r \ar[r]^-{\iota} & ^{\tau}\mathcal{E}_0(D_\infty)}
\]
each $f_i$ is a modification from $\mathcal{E}_{i-1}$ to $\mathcal{E}_{i}$ of relative position $\mu_i$ around $\Gamma_{x_i}$, respecting the filtration as in the definition of the Hecke stack, and $\iota$ at the end is an isomorphism. $^{\tau}\mathcal{E}_0$ is the pull back of $\mathcal{E}_0$ along the map:
\begin{equation*}
X\times_kS\xrightarrow{id_X\times Fr_S} X\times_kS
\end{equation*}
\end{itemize}
Now the restriction above is clear: Modification by $f_i$ changes the degree by $\mu_i$, therefore:
\begin{equation*}
\text{deg}^{\tau}\mathcal{E}_0(D_\infty)-\text{deg}\mathcal{E}_0=\sum_{i=1}^{r}\mu_i
\end{equation*}
This implies that one must have the above restriction, because otherwise the stack would be empty.

Define the moduli of shtukas for $G$ as the quotient by the discrete groupoid $Pic_X(k)$:
\begin{equation*}
Sht^{\underline{\mu}}_G(\Sigma; D_\infty)=Sht^{\underline{\mu}}_2(\Sigma; D_\infty)/Pic_X(k)
\end{equation*}
Then it fits into the following Cartesian diagram:
\begin{equation}
\xymatrix{Sht^{\underline{\mu}}_G(\Sigma;D_\infty) \ar[r] \ar[d]^{\omega_0} & Hk_G^{\underline{\mu}}(\Sigma)\times_k\mathfrak{S}_\infty \ar[d]^{(\widetilde{p}_0, AL(-D_\infty)\circ (\widetilde{p}_r\times id_{\mathfrak{S}_\infty})} \\
Bun_G(\Sigma) \ar[r]^-{(id, Fr)} & Bun_G(\Sigma)\times Bun_G(\Sigma)}
\end{equation}
Let 
\begin{equation}
\Pi^{\underline{\mu}}_{G,D_\infty}=(\pi^{\underline{\mu}}_G, \pi_{G,\infty}): Sht^{\underline{\mu}}_G(\Sigma;D_\infty)\longrightarrow X^r\times_k\mathfrak{S}_\infty
\end{equation}
be the map recording the points of modification and the $\mathfrak{S}_\infty$ structure of $S$.

Yun and Zhang proved the following properties of $Sht^{\underline{\mu}}_G(\Sigma;D_\infty)$ (\cite{YZ2}, Proposition 3.9):
\begin{proposition}
\begin{enumerate}
\item The stack $Sht^{\underline{\mu}}_G(\Sigma; D_\infty)$ is Deligne-Mumford of dimension $2r$.
\item $\Pi^{\underline{\mu}}_{G,D_\infty}=(\pi^{\underline{\mu}}_G, \pi_{G,\infty})$ is separated, and is smooth of relative dimension $r$ when restricted to $(X-\Sigma)^r\times\mathfrak{S}_\infty$
\end{enumerate}
\end{proposition}
Recall that we have fixed a rational divisor on the scheme $X\times_k\mathfrak{S}_\infty$, $D_\infty=\sum_{x\in\Sigma}\frac{1}{2}\textbf{x}^{(1)}$. Define:
\begin{equation}\label{G-shtuka}
Sht^{\underline{\mu}}_G(\Sigma; \Sigma_\infty):=Sht^{\underline{\mu}}_G(\Sigma; D_\infty)
\end{equation}
For the following $r$s, $\mu$s and $D_\infty$s:
\begin{equation}
r=\#\Sigma_\infty \text{    mod}2\quad
D_\infty\in\mathscr{D}_\infty, \text{ and }\sum^r_{i=1}\mu_i=2\text{deg}D_\infty
\end{equation}
They proved that for the above $D_\infty$s, the Atkin-Lehner involutions:
\begin{equation}
AL(-D_\infty): Bun_G(\Sigma)\times_k\mathfrak{S}_\infty\longrightarrow Bun_G(\Sigma)
\end{equation}
all agree with $AL(-D_\infty^{(1)})$ (of course not for $\widetilde{AL}(-D_\infty)$), and the right hand side of \eqref{G-shtuka} is independent of choice of $\underline{\mu}$ and $D_\infty$ as long as the above conditions are satisfied.

\subsection{Hecke correspondence for moduli of shtukas} Now review the Hecke correspondences (\cite{HS}, section 3.3). Let's define the partial spherical Hecke algebra as the following:
\begin{equation}
\mathscr{H}^{\Sigma}_G=\otimes_{x\in|X-\Sigma|}\mathscr{H}_x
\end{equation}
in which $\mathscr{H}_x$ is the local spherical Hecke algebra with $\mathbb{Q}$ coefficients. It has a basis indexed by the group of effective divisors of $X-\Sigma$, i.e. $Div^+(X-\Sigma)$. 

Let $D\in Div^{+}(X-\Sigma)$ be an effective divisor, then define the associated "vertical Hecke modification" stack as follows:
\begin{definition}
$Sht^{\underline{\mu}}_2(\Sigma; D_\infty; h_D)$ is the functor from $k-$schemes to groupoids whose  category of $S$ points is the following data:
\begin{itemize}
\item A map: $S\longrightarrow\mathfrak{S}_\infty$, which gives $x^{(1)}: S\longrightarrow \text{Spec}k_x$ for each $x\in\Sigma_\infty$;
\item $r$ morphisms $x_i: S\longrightarrow X$ for $i=1,...,r$, with the graph $\Gamma_{x_i}\subset X\times_k S$.
\item Two of the points of $Sht^{\underline{\mu}}_2(\Sigma; D_\infty)(S)$, say $(\mathcal{E}^\dagger_i, f_i, \iota, ...)$ and $(\mathcal{E}'^\dagger_i, f'_i, \iota', ...)$ sharing the above data, i.e. map to the same point in $\mathfrak{S}_\infty$, points of modifications $x_i$, together with injections of coherent sheaves $\phi_i$, fitting into the following diagram:
\begin{equation}
\xymatrix{\mathcal{E}_0 \ar@{-->}[r]^{f_1} \ar[d]^{\phi_0} & \mathcal{E}_1 \ar@{-->}[r]^{f_2} \ar[d]^{\phi_1} &......\ar@{-->}[r]^{f_r} & \mathcal{E}_r \ar[r]^-{\iota} \ar[d]^{\phi_r} & ^{\tau}\mathcal{E}_0(D_\infty)\ar[d]^{^{\tau}\phi_0}\\ 
\mathcal{E}'_0 \ar@{-->}[r]^{f'_1} & \mathcal{E}'_1 \ar@{-->}[r]^{f'_2} &......\ar@{-->}[r]^{f'_r} & \mathcal{E}'_r \ar[r]^-{\iota'} & ^{\tau}\mathcal{E}'_0(D_\infty)}
\end{equation}
such that each $\phi_i$ preserves the Iwahori level structure (i.e. sending $\mathcal{E}(-\frac{1}{2}x)$ to $\mathcal{E}'(-\frac{1}{2}x)$ for all $x\in\Sigma$); the induced maps $\phi_i: det(\mathcal{E}_{i-1})\longrightarrow det(\mathcal{E}_i)$ has divisor $D\times_kX\subset X\times_kS$. 
\end{itemize}
\end{definition}
As usual, define:
\begin{equation}
Sht^{\underline{\mu}}_G(\Sigma;\Sigma_\infty; h_D):=Sht^{\underline{\mu}}_2(\Sigma; D_\infty; h_D)/Pic_X(k)
\end{equation}
and the map recording the points of modifications and $\mathfrak{S}_\infty$-structure:
\begin{equation*}
\Pi^r_G(h_D): Sht^{\underline{\mu}}_G(\Sigma; \Sigma_\infty; h_D)\longrightarrow X^r\times_k\mathfrak{S}_\infty
\end{equation*}

By definition, it is a self-correspondence of $Sht^{\underline{\mu}}_G(\Sigma; \Sigma_\infty; h_D)$ over $X^r\times_k\mathfrak{S}_\infty$, i.e. one has the following diagram:
\begin{equation}
\xymatrix{ & Sht^{\underline{\mu}}_G(\Sigma; \Sigma_\infty; h_D) \ar[ld]_{\overleftarrow{p}} \ar[rd]^{\overrightarrow{p}} & \\
Sht^{\underline{\mu}}_G(\Sigma; \Sigma_\infty) \ar[rd] & & Sht^{\underline{\mu}}_G(\Sigma; \Sigma_\infty) \ar[ld]\\
& X^r\times_k\mathfrak{S}_\infty &}
\end{equation}
It has the following properties:
\begin{proposition}
Let $D\in Div^{+}(X-\Sigma)$.
\begin{enumerate}

\item Both $\overleftarrow{p}$ and $\overrightarrow{p}$ are proper and representable;
\item Over $(X-\Sigma)^r$, both of $\overleftarrow{p}$ and $\overrightarrow{p}$ are fintie \'etale;
\item The fibers of $\Pi^r_G(h_D): Sht^{\underline{\mu}}_G(\Sigma; \Sigma_\infty; h_D)\longrightarrow X^r\times_k\mathfrak{S}_\infty$ have dimension $r$.
\end{enumerate}
\end{proposition}
For the detail of the proof, one just reads Yun-Zhang's second volume.

These Hecke correspondences induces the following ring homomorphisms:
\begin{proposition}
One can extend the map:
\begin{equation}
h_D\longrightarrow(\overleftarrow{p}\times\overrightarrow{p})_{*}Sht^{\underline{\mu}}_G(\Sigma; \Sigma_\infty; h_D)
\end{equation}
to a ring homomorphism:
\begin{equation}
\mathscr{H}_G^\Sigma\longrightarrow{_cCh_{2r}(Sht^{\underline{\mu}}_G(\Sigma; \Sigma_\infty)\times Sht^{\underline{\mu}}_G(\Sigma; \Sigma_\infty))_{\mathbb{Q}}}
\end{equation}
which in turn allow $\mathscr{H}_G^\Sigma$ to act on the group of compactly supported Chow cycles.
\end{proposition}

\subsection{The Heegner-Drinfeld cycles}
Now we define the Heegner-Drinfeld cycles for $Y_1$ and $Y_2$, with some auxilary choice of data.

As one has seen in the last section, over a certain point $x\in\Sigma_\infty$, there is one point over it in each $Y_i, i=1,2$ respectively, which we denoted $y_x, z_x$, and there are two points over it in $Y$. We choose one of these two points for each $x\in\Sigma_\infty$ and denoted it by $v_x$, and call the other point by $v'_x$. The corresponding points below them in $Y_3$ are denoted by $w_x$ and $w'_x$ respectively. Then define the following:
\begin{align*}
\mathfrak{S}'_{\infty,1}=\Pi_{x\in\Sigma}\text{Spec}k_{y_x}\\
\mathfrak{S}'_{\infty,2}=\Pi_{x\in\Sigma}\text{Spec}k_{z_x}\\
\widetilde{\mathfrak{S}}_\infty=\Pi_{x\in\Sigma}\text{Spec}k_{v_x}
\end{align*}
From the definition, one can see that there are cononical maps:
\begin{align*}
\widetilde{\mathfrak{S}}_\infty=\Pi_{x\in\Sigma}\text{Spec}k_{v_x}\longrightarrow\mathfrak{S}'_{\infty,1}=\Pi_{x\in\Sigma}\text{Spec}k_{y_x};
\\
\widetilde{\mathfrak{S}}_\infty=\Pi_{x\in\Sigma}\text{Spec}k_{v_x}\longrightarrow\mathfrak{S}'_{\infty,2}=\Pi_{x\in\Sigma}\text{Spec}k_{z_x}
\end{align*}
Similar to the case of $\mathfrak{S}_\infty$, denote the following by $\textbf{y}_x^{(1)}, \textbf{z}_x^{(1)}, \textbf{v}_x^{(1)}$:
\begin{align*}
\mathfrak{S}'_{\infty,1}&\longrightarrow \text{Spec}k_{y_x}\longrightarrow Y_1\\
\mathfrak{S}'_{\infty,2}&\longrightarrow \text{Spec}k_{z_x}\longrightarrow Y_2\\
\widetilde{\mathfrak{S}}_\infty&\longrightarrow \text{Spec}k_{v_x}\longrightarrow Y
\end{align*}
in which the first arrows are the projection and the second arrows are the natural injections. Similar to $\textbf{x}^{(i)}$, define $\textbf{y}^{(i)}_x=\textbf{y}_x^{(i-1)}\circ Fr_{\mathfrak{S}'_{\infty,1}}, \textbf{z}^{(i)}_x=\textbf{z}_x^{(i-1)}\circ Fr_{\mathfrak{S}'_{\infty,2}}, \textbf{v}^{(i)}_x=\textbf{v}_x^{(i-1)}\circ Fr_{\widetilde{\mathfrak{S}}_\infty}$. Then one has the following decomposition:
\begin{equation*}
\xymatrix{ & \text{Spec}k_{v_x}\times_k
\widetilde{\mathfrak{S}}_\infty=\coprod_{i=1}^{2d_x}\Gamma_{\textbf{v}_x^{(i)}}  \ar[ld] \ar[rd] &\\
\text{Spec}k_{y_x}\times_k\mathfrak{S}'_{\infty,1}=\coprod_{i=1}^{2d_x}\Gamma_{\textbf{y}_x^{(i)}} \ar[rd] &  & \text{Spec}k_{z_x}\times_k\mathfrak{S}'_{\infty,2}=\coprod_{i=1}^{2d_x} \ar[ld]\Gamma_{\textbf{z}_x^{(i)}}\\
& \text{Spec}k_x\times_k\mathfrak{S}_\infty=\coprod_{i=1}^{d_x}\Gamma_{\textbf{x}^{(i)}}&
}
\end{equation*}
For $\Gamma_{\textbf{y}_x^{(i)}}$ and $\Gamma_{\textbf{z}_x^{(i)}}$, $\Gamma_{\textbf{v}_x^{(i)}}$ is the unique piece in $\text{Spec}k_{v_x}\times_k
\widetilde{\mathfrak{S}}_\infty$ sitting above them. For this reason, we also use $\Gamma_{(\textbf{y}_x^{(i)}, \textbf{z}_x^{(i)})}$ to denote $\Gamma_{\textbf{v}_x^{(i)}}$.

However, for each $\Gamma_{\textbf{y}_x^{(i)}}$ (or $\Gamma_{\textbf{z}_x^{(i)}}$), $\Gamma_{\textbf{v}_x^{(i)}}$ does not exhaust its preimage in $Y\times_k\widetilde{\mathfrak{S}}_\infty$. Just as  $\text{Spec}k_{v_x}\times_kY$, $\text{Spec}k_{v'_x}\times_kY$ also decomposes into $2d_x$ pieces, and for $\Gamma_{\textbf{y}_x^{(i)}}$ and $\Gamma_{\textbf{z}_x^{(d_x+i)}}$ ($(d_x+i)$ is understood to be mod$2d_x$), there is a unique piece of $\text{Spec}k_{v'_x}\times_kY$ sitting over them, let it be $\Gamma_{(\textbf{y}^{(i)}, \textbf{z}^{(d_x+i)})}$. Then one has:
\begin{equation*}
Speck_{v'_x}\times_kY=\coprod_{i=1}^{2d_x}\Gamma_{(\textbf{y}^{(i)}, \textbf{z}^{(d_x+i)})}
\end{equation*}

For a $k$-scheme $S$ with $S\longrightarrow\widetilde{\mathfrak{S}}_\infty$, one gets $y^{(i)}:S\longrightarrow Speck_{y_x}\longrightarrow Y_1$ and $z^{(i)}:S\longrightarrow Speck_{y_x}\longrightarrow Y_2$. One has a diagrams for $S$ similar to the one above:
\begin{equation*}
\xymatrix{ & \text{Spec}k_{v_x}\times_k
S=\coprod_{i=1}^{2d_x}\Gamma_{(y_x^{(i)}, z_x^{(i)})}  \ar[ld] \ar[d] \ar[rd] &\\
\text{Spec}k_{y_x}\times_kS=\coprod_{i=1}^{2d_x}\Gamma_{y_x^{(i)}} \ar[rd] & \text{Spec}k_{w_x}\times_kS \ar[d] & \text{Spec}k_{z_x}\times_kS=\coprod_{i=1}^{2d_x} \ar[ld]\Gamma_{z_x^{(i)}}\\
& \text{Spec}k_x\times_kS=\coprod_{i=1}^{d_x}\Gamma_{x^{(i)}}&
}
\end{equation*}
and
\begin{equation*}
\xymatrix{ & \text{Spec}k_{v'_x}\times_k
S=\coprod_{i=1}^{2d_x}\Gamma_{(y_x^{(i)}, z_x^{(d_x+i)})}  \ar[ld] \ar[d] \ar[rd] &\\
\text{Spec}k_{y_x}\times_kS=\coprod_{i=1}^{2d_x}\Gamma_{y_x^{(i)}} \ar[rd] &\text{Spec}k_{w'_x}\times_kS \ar[d]   & \text{Spec}k_{z_x}\times_kS=\coprod_{i=1}^{2d_x} \ar[ld]\Gamma_{z_x^{(i)}}\\
& \text{Spec}k_x\times_kS=\coprod_{i=1}^{d_x}\Gamma_{x^{(i)}}&
}
\end{equation*}
For simplicity, I will just use $y_x^{(i)}, z_x^{(i)}$ and $(y_x^{(i)}, z_x^{(i)})$ to denote $\Gamma_{y_x^{(i)}}, \Gamma_{z_x^{(i)}}$ and $\Gamma_{(y_x^{(i)}, z_x^{(i)})}$. For example, I will write $\mathcal{L}(y_x^{(i)})$ for $\mathcal{L}(\Gamma_{y_x^{(i)}})$ to mean twisting the line bundle $\mathcal{L}$ on $X\times_kS$ by the divisor $\Gamma_{y_x^{(i)}}$.

Now let's review the definition of $T_i$-shtukas, I will only do this for $i=1$ and $i=2$ is all the same.

First Let 
\begin{equation}
Bun_{T_i}=Pic_{Y_i}/Pic_X
\end{equation}
and 
\begin{equation}
Hk^{\underline{\mu}}_{T_i}=Hk^{\underline{\mu}}_{1,Y_i}/Pic_X
\end{equation}
for $\underline{\mu}\in\{\pm 1\}^r$ be as usual. Of course, one knows the following properties of them (\cite{YZ1}, Remark 5.2, section 5.4.4):

\begin{proposition}
\begin{enumerate}
\item $Bun_{T_i}$ is a Deligne-Mumford stack of dimension $g-1$.
\item For $0\leq i \leq r$, the projection map $p_i: Hk_{T_1}^{\underline{\mu}}\longrightarrow Bun_{T_1}$ is smooth of relative dimension $r$.
\item For $0\leq i \leq r$, the morphism $(p_i, \pi^{\underline{\mu}}_{Hk}): Hk_{T_1}^{\underline{\mu}}\longrightarrow Bun_{T_1}\times_kY_1^r$ is an isomorphism.
\end{enumerate}
\end{proposition}

Let $\widetilde{D}_\infty$ be a divisor in $Y\times_k\widetilde{\mathfrak{S}}_\infty$ of the following form:
\begin{equation}
\widetilde{D}_\infty=\sum_{x\in\Sigma_\infty, 1\leq i\leq 2d_x}c^{(i)}_x(\textbf{y}_x^{(i)}, \textbf{z}_x^{(i)})
\end{equation}
Project to $Y_1\times_k\mathfrak{S}'_{\infty,1}$, one gets:
\begin{equation}
D_{1,\infty}=\sum_{x\in\Sigma_\infty, 1\leq i\leq 2d_x}c^{(i)}_x\textbf{y}_x^{(i)}
\end{equation}
Then one has the Atkin-Lehner involution:

\begin{align*}
\widetilde{AL}(\widetilde{D}_\infty):=\widetilde{AL}(D_{1,\infty}): Pic_{Y_1}\times_k\widetilde{\mathfrak{S}}_\infty&\longrightarrow Pic_{Y_1}\\
AL(\widetilde{D}_\infty):=AL(D_{1,\infty}):=Bun_{T_1}\times_k\widetilde{\mathfrak{S}}_\infty&\longrightarrow Bun_{T_1}\\
\mathcal{L}\longmapsto \mathcal{L}(\sum_{x\in\Sigma^\infty}c^{(i)}_x\textbf{y}_x^{(i)})
\end{align*}

Let $\underline{\mu}\in\{\pm1\}^r$ and $\widetilde{D}_\infty$ satisfy the equation:
\begin{equation*}
\sum_{i=1}^r\mu_i=\sum_{x\in\Sigma_\infty, 1\leq i\leq 2d_x}c^{(i)}_x
\end{equation*}

Of course, this is due to the same degree reason as in the previous subsection.
Define $Sht^{\underline{\mu}}_{T_i}(D_{1,\infty})'$ to be following fiber product:
\begin{equation}
\xymatrix{Sht^{\underline{\mu}}_{T_i}(D_{1,\infty})' \ar[r] \ar[d]^{\omega_0} & Hk_{T_i}^{\underline{\mu}}(\Sigma)\times_k\widetilde{\mathfrak{S}}_\infty \ar[d]^{(\widetilde{p}_0, AL(-\widetilde{D}_{1,\infty})\circ (\widetilde{p}_r\times id_{\mathfrak{S}_\infty})} \\
Bun_{T_i} \ar[r]^-{(id, Fr)} & Bun_{T_i}\times Bun_{T_i}}
\end{equation}
Given $Sht^{\underline{\mu}}_{T_i}(D_{1,\infty})$, the moduli of shtukas defined using $\mathfrak{S}'_{\infty,1}$, the above one is nothing but 
\begin{equation}
Sht^{\underline{\mu}}_{T_i}(D_{1,\infty})\times_{\mathfrak{S}_{1,\infty}'}\widetilde{\mathfrak{S}}_\infty
\end{equation}

As in Yun-Zhang, one only works with the following divisor in $Y\times_k\widetilde{\mathfrak{S}}_\infty$:
\begin{equation}
\mu_\infty\cdot\widetilde{\Sigma}_\infty=\sum_{x\in\Sigma^\infty}\mu_x(\textbf{y}_x^{(1)}, \textbf{z}_x^{(1)})
\end{equation}
Then projects to $Y_1\times_k\mathfrak{S}'_\infty$, we get
\begin{equation}
\mu_\infty\cdot{\Sigma'_1}_\infty=\sum_{x\in\Sigma^\infty}\mu_x\textbf{y}_x^{(1)}
\end{equation}
Fix $r=\#\Sigma^\infty$ mod $2$, then for this special choice of divisor, we have the following: 

\begin{equation}
Sht^{\underline{\mu}}_{T_i}(\mu_\infty\cdot\Sigma'_{1,\infty})'=Sht^{\underline{\mu}}_{T_i}(\mu_\infty\cdot\Sigma'_{1,\infty})\times_{\mathfrak{S}_{1,\infty}'}\widetilde{\mathfrak{S}}_\infty
\end{equation}
According to Yun-Zhang, one has: 
\begin{equation}
Sht^{\underline{\mu}}_{T_i}(\mu_\infty\cdot\Sigma'_{1,\infty})=
Sht^{\underline{\mu}}_{T_i}(D_{1,\infty})
\end{equation}
for any choice of $D_{1,\infty}=\mu_\infty\cdot{\Sigma'_1}_\infty$ mod $\nu_1^*\mathscr{D}_\infty$. Also, Yun-Zhang showed (\cite{YZ2}, Corollary 4.3):
\begin{proposition}
$Sht^{\underline{\mu}}_{T_i}(\mu_\infty\cdot\Sigma'_{1,\infty})$ is a smooth Deligne-Mumford stack of dimension $r$ 
\end{proposition}

Therefore the stack we are going to consider here $Sht^{\underline{\mu}}_{T_i}(\mu_\infty\cdot\Sigma'_{1,\infty})'$ is also a smooth Deligne-Mumford stack of dimension $r$.

To define Heegner-Drinfeld cycles, one first need to know how to assign an Iwahori level structure to the push foward of a line bundle.
To do this, we now fix a tuple of $1$ and $-1$s:
\begin{equation}
\mu_{\Sigma}=(\mu_f, \mu_\infty)\in\{\pm1\}^{\Sigma_f}\times\{\pm1\}^{\Sigma_\infty}
\end{equation}
define the following morphism:
\begin{definition}\label{phd}
Define:
\begin{equation}
\widetilde{\theta}^{\mu_\Sigma}_{1,Bun}: Pic_{Y_1}\times\widetilde{\mathfrak{S}}_\infty\longrightarrow Bun_2(\Sigma)
\end{equation}
in the following way: For a $k$-scheme $S$ with a map $S\longrightarrow\widetilde{\mathfrak{S}}_\infty$, the image of $\mathcal{L}$, a line bunde over $Y_1\times S$, is $\nu_{1,*}\mathcal{L}$, whose Iwahori level structure, i.e. the "lattice filtration" at $\Sigma$ is given by:
\begin{enumerate}
\item At $x\in\Sigma_f$:
\begin{itemize}
\item If $\mu_x=1$, define $\nu_{1,*}\mathcal{L}(-\frac{1}{2})$ to be $\nu_{1,*}(\mathcal{L}(-y_x))$
\item If $\mu_x=-1$, define $\nu_{1,*}\mathcal{L}(-\frac{1}{2})$ to be $\nu_{1,*}(\mathcal{L}(-y'_x))$
\end{itemize}
\item At $x\in\Sigma_\infty$
\begin{itemize}
\item If $\mu_x=1$, define $\nu_{1,*}\mathcal{L}(-\frac{1}{2})$ to be $\nu_{1,*}(\mathcal{L}(-y^{(1)}_x-y^{(2)}_x-...-y^{(d_x)}))$
\item If $\mu_x=-1$, define $\nu_{1,*}\mathcal{L}(-\frac{1}{2})$ to be $\nu_{1,*}(\mathcal{L}(-y^{(d_x+1)}_x-y^{(d_x+2)}_x-...-y^{(2d_x)}))$
\end{itemize}
\end{enumerate}
\end{definition}
This definition is the same as Yun-Zhang's definition applied to the $\mathfrak{S}_{\infty,1}'$ structure for $S$ induced by its $\widetilde{\mathfrak{S}}_\infty$ with the canonical projection mentioned above.

For $Y_2$ we define the map
\begin{equation}
\widetilde{\theta}^{\mu_\Sigma}_{2,Bun}: Pic_{Y_2}\times\widetilde{\mathfrak{S}}_\infty\longrightarrow Bun_2(\Sigma)
\end{equation}
in completely the same way. 

By quotienting out $Pic_X$, one gets:
\begin{equation}
\theta^{\mu_\Sigma}_{1,Bun}: Bun_{T_1}\times\widetilde{\mathfrak{S}}_\infty\longrightarrow Bun_G(\Sigma)
\end{equation}
and
\begin{equation}
\theta^{\mu_\Sigma}_{2,Bun}: Bun_{T_2}\times\widetilde{\mathfrak{S}}_\infty\longrightarrow Bun_G(\Sigma)
\end{equation}

We can define the morphism on the Hecke stacks:
\begin{equation}\label{3.31}
Hk^{\underline{\mu}}_{T_1}\times_k\widetilde{\mathfrak{S}}_\infty\longrightarrow Hk^{\underline{\mu}}_G(\Sigma)
\end{equation}
By applying $\theta^{\mu_\Sigma}_{Bun,1}$ to each single $\mathcal{L}_i$ in the sequence of line bundles.

As in their paper, define the following $\#$-Atkin-Lehner for $T_1$:
\begin{equation}
AL^{\#}_{T_1, \mu_\infty}:=(AL(-\mu_\infty\cdot\widetilde{\Sigma}_\infty), Fr_{\widetilde{\mathfrak{S}}_\infty}): Bun_{T_1}\times{\widetilde{\mathfrak{S}}_\infty}\longrightarrow Bun_{T_1}\times{\widetilde{\mathfrak{S}}_\infty}
\end{equation}
Then one has the following commutative diagram:
\begin{equation}\label{3.33}
\xymatrix{ Bun_{T_1}\times{\widetilde{\mathfrak{S}}_\infty} \ar[r]^{AL^{\#}_{T_1, \mu_\infty}} \ar[d]_{(\theta^{\mu_\Sigma}_{1, Bun},\nu_\infty)} &  Bun_{T_1}\times{\widetilde{\mathfrak{S}}_\infty} \ar[d]^{\theta^{\mu_\Sigma}_{1, Bun}}\\
Bun_G(\Sigma)\times_k\mathfrak{S}_\infty \ar[r]^-{AL_{G,\infty}} & Bun_G(\Sigma)
}
\end{equation}

By unraveling the meaning of the functor of points of $Sht^{\underline{\mu}}_{T_1}(\mu_\infty\cdot\Sigma_{1,\infty})'$, one finds that it also fits in the following fiber diagram:
\begin{equation}\label{3.34}
\xymatrix{Sht^{\underline{\mu}}_{T_1}(\mu_\infty\cdot\Sigma_{1,\infty})' \ar[r] \ar[d] & Hk_{T_1}^{\underline{\mu}}(\Sigma)\times_k\widetilde{\mathfrak{S}}_\infty \ar[d]^{(p_0\times id_{\widetilde{\mathfrak{S}}_\infty}, AL^{\#}_{T_1,\mu_\infty}\circ (p_r\times id_{\widetilde{\mathfrak{S}}_\infty})} \\
Bun_{T_1} \ar[r]^-{(id, Fr)} & (Bun_{T_1}\times{\widetilde{\mathfrak{S}}_\infty})\times (Bun_{T_1}\times{\widetilde{\mathfrak{S}}_\infty})}
\end{equation}
Then by \eqref{3.33}, \eqref{3.34}, together with the morphism of the Hecke stacks \eqref{3.31}, one gets the Heegner-Drinfeld map for $Y_1$:
\begin{equation}
\theta^{\mu}_1: Sht^{\underline{\mu}}_{T_1}(\mu_\infty\cdot\Sigma_{1,\infty})'\longrightarrow Sht^{\underline{\mu}}_G(\Sigma;\Sigma_{\infty})
\end{equation}
By remembering the $\widetilde{\mathfrak{S}}_\infty$-structure, one actually gets:
\begin{equation}
\theta'^{\mu}_1: Sht^{\underline{\mu}}_{T_1}(\mu_\infty\cdot\Sigma_{1,\infty})'\longrightarrow Sht^{\underline{\mu}}_G(\Sigma;\Sigma_{\infty})\times_{\mathfrak{S}_\infty}\widetilde{\mathfrak{S}}_\infty
\end{equation}

Repeat all the above for $Y_2$, one gets the Heegner-Drinfeld map for $Y_2$:
\begin{equation}
\theta'^{\mu}_2: Sht^{\underline{\mu}}_{T_2}(\mu_\infty\cdot\Sigma_{2,\infty})'\longrightarrow Sht^{\underline{\mu}}_G(\Sigma;\Sigma_{\infty})\times_{\mathfrak{S}_\infty}\widetilde{\mathfrak{S}}_\infty
\end{equation}

Denote the common right hand side by $Sht^{\underline{\mu}}_G(\Sigma;\Sigma_{\infty})'$. Note that unlike Yun-Zhang, there is no base with respect to $Y^r_i\longrightarrow X^r$.

By pushing foward the fundamental cycles of the left hand side of the equation, one gets classes in the middle dimension Chow groups of the right hand side:
\begin{definition} Define the Heegner-Drinfeld cycle for $Y_1$ and $Y_2$ as follows:
\begin{align}
\mathcal{Z}^{\mu}_1:=\theta'^{\mu}_1[Sht^{\underline{\mu}}_{T_1}(\mu_\infty\cdot\Sigma_{1,\infty})']\in Ch_{c,r}(Sht^{\underline{\mu}}_G(\Sigma;\Sigma_{\infty})')\\
\mathcal{Z}^{\mu}_2:=\theta'^{\mu}_2[Sht^{\underline{\mu}}_{T_2}(\mu_\infty\cdot\Sigma_{2,\infty})']\in Ch_{c,r}(Sht^{\underline{\mu}}_G(\Sigma;\Sigma_{\infty})')
\end{align}
\end{definition}
Finishing all the setups, it is the time for us to state the question:

``Compute" the following number:
\begin{equation}
\langle\mathcal{Z}^{\mu}_1, f*\mathcal{Z}^{\mu}_2\rangle
\end{equation}

for an $f\in\mathscr{H}^{\Sigma}_G$. One can first take $f$ to be a basis element indexed by an effective divisor D of $X-\Sigma$, namely $h_D$. $\langle\mathcal{Z}^{\mu}_1, h_D*\mathcal{Z}^{\mu}_2\rangle$ can be computed by the following fiber diagram:
\begin{equation}\label{product}
\xymatrix{\ast \ar[r] \ar[d] & Sht^{\underline{\mu}}_G(\Sigma;\Sigma_\infty; h_D)'\ar[d]\\
Sht^{\underline{\mu}}_{T_1}(\mu_\infty\cdot\Sigma_{1,\infty})'\times Sht^{\underline{\mu}}_{T_2}(\mu_\infty\cdot\Sigma_{2,\infty})' \ar[r]^-{\theta'^{\mu}_1\times\theta'^{\mu}_2} & Sht^{\underline{\mu}}_G(\Sigma;\Sigma^{\infty})'\times Sht^{\underline{\mu}}_G(\Sigma;\Sigma^{\infty})'
}
\end{equation}

In which $Sht^{\underline{\mu}}_G(\Sigma;\Sigma_\infty; h_D)'$ is defined as:
\begin{equation}
Sht^{\underline{\mu}}_G(\Sigma;\Sigma_\infty; h_D)':=Sht^{\underline{\mu}}_G(\Sigma;\Sigma_\infty; h_D)\times_{\mathfrak{S}_\infty}\widetilde{\mathfrak{S}}_\infty
\end{equation}
and we will later give $\ast$ a name. Then $\langle\mathcal{Z}^{\mu}_1, h_D*\mathcal{Z}^{\mu}_2\rangle=deg(\theta'^{\mu}_1\times\theta'^{\mu}_2)^{!}[Sht^{\underline{\mu}}_G(\Sigma;\Sigma_\infty; h_D)']$.

For simplicity, for $f$, let $\mathbb{I}^\mu(f)=(\prod_{x\in\Sigma_\infty}2d_x)^{-1}\langle\mathcal{Z}^{\mu}_1, f*\mathcal{Z}^{\mu}_2\rangle$. The clumsy factor before the intersection pairing is for the nomalization purpose. To remove it, one can fix a $\overline{k}$-point of $\widetilde{\mathfrak{S}}_\infty$, say $\xi$, and consider the base change from $\widetilde{\mathfrak{S}}_\infty$ to $\overline{k}$:
\begin{align*}
Sht^{\overline{\mu}}_G(\Sigma;\xi)&:=Sht^{\underline{\mu}}_G(\Sigma;\Sigma_\infty)\times_{\mathfrak{S}_\infty}\xi\\
Sht^{\underline{\mu}}_{T_i}(\mu_\infty\cdot\xi)&:=Sht^{\underline{\mu}}_T(\mu_\infty\cdot\Sigma_{i,\infty})'\times_{\widetilde{\mathfrak{S}}_\infty}\xi
\end{align*}
Then one has:
\begin{equation*}
Sht^{\underline{\mu}}_{T_i}(\mu_\infty\cdot\xi)\xrightarrow{\theta^{'\mu}_{i,\xi}}Sht^{\underline{\mu}}_G(\Sigma;\Sigma_\infty)'\times_{\widetilde{\mathfrak{S}}_\infty}\xi=Sht^{\underline{\mu}}_G(\Sigma;\Sigma_\infty)\times_{\mathfrak{S}_\infty}\xi=Sht^{\underline{\mu}}_G(\Sigma;\xi)
\end{equation*}

Define $\mathcal{Z}_i^\mu(\xi):=\theta^{'\mu}_{i,\xi}[Sht^{\underline{\mu}}_{T_i}(\mu_\infty\cdot\xi)]\in Ch_{c,r}(Sht^{\mu}_G(\Sigma;\xi))_{\mathbb{Q}}$ for $i=1,2$. Actually the pull back of $\mathcal{Z}^\mu_i$ to $Sht^{\underline{\mu}}_G(\Sigma;\Sigma_\infty)'\times_{\widetilde{\mathfrak{S}}_\infty}\overline{k}$ splits into is just the disjoint union of $\mathcal{Z}^\mu_i(\xi)$ as $\xi$ runs through all the $\prod_{x\in\Sigma_\infty}2d_x$ $\overline{k}$-points of $\widetilde{\mathfrak{S}}_\infty$. Following the method of Yun-Zhang (\cite{YZ2}, Corollary 4.11), one can prove the following:
\begin{lemma}
For any two $\xi, \xi'\in\widetilde{\mathfrak{S}}(\overline{k})$, one has:
\begin{equation}
\langle\mathcal{Z}^{\mu}_1(\xi), f*\mathcal{Z}^{\mu}_2(\xi)\rangle_{Sht^{\underline{\mu}}_G(\Sigma;\xi)}=\langle\mathcal{Z}^{\mu}_1(\xi'), f*\mathcal{Z}^{\mu}_2(\xi')\rangle_{Sht^{\underline{\mu}}_G(\Sigma;\xi')}
\end{equation}
In particular:
\begin{equation}
\mathbb{I}^\mu(f)=\langle\mathcal{Z}^{\mu}_1(\xi), f*\mathcal{Z}^{\mu}_2(\xi)\rangle_{Sht^{\underline{\mu}}_G(\Sigma;\xi)}
\end{equation}
\end{lemma}

By the general theory of Yun-Zhang, one expects that this number can be computed by applying the Grothendieck-Lefschetz trace formula to some cohomological self-correspondence of certain ``Hitchin type" fiberation derived from the moduli of shtukas and the constant sheaf for $D$ of large enough degree. This is the topic of the next section.

\section{The Hitchin type fiberation for the Heegner-Drinfeld cycles}
\subsection{The moduli space $\mathcal{M}^{\mu_\Sigma}_d$} 
Unraveling the definition of diagram used to compute $\langle\mathcal{Z}^{\mu}_1, h_D*\mathcal{Z}^{\mu}_2\rangle$, i.e.\eqref{3.42} and the general philosophy of Yun-Zhang's work, we are led to the following definition. Fixing $\mu_\infty=(\mu_f, \mu_\infty)$ as in the last section.
\begin{definition}
Define $\widetilde{\mathcal{M}}^{\mu_\Sigma}_d$ as the stack whose funtor of point on a $k$-scheme $S$ is the following:
\begin{itemize}
\item A map $S\longrightarrow\widetilde{\mathfrak{S}}_\infty$
\item A line bundle $\mathcal{L}$ on $Y_1\times_kS$ and a line bundle $\mathcal{F}$ on $Y_2\times_kS$.
\item A morphim of coherent sheaves on $X\times_kS$, $\phi$ from $\widetilde{\theta}^{\mu_\Sigma}_{1,Bun}(\mathcal{L})$ to $\widetilde{\theta}^{\mu_\Sigma}_{2,Bun}(\mathcal{F})$, more precisely:
\begin{equation}
\phi: \nu_{1,*}\mathcal{L}\longrightarrow\nu_{2,*}\mathcal{F}
\end{equation}
sending each ``sublattice"  of $\nu_{1,*}\mathcal{L}$ defined in \eqref{phd} to the corresponding one of $\nu_{2,*}\mathcal{F}$. And $\phi$ has the following property: the induced map on the determinant line bundle 
\begin{equation}
\text{det}{\phi}: \text{det}(\nu_{1,*}\mathcal{L})\longrightarrow \text{det}(\nu_{2,*}\mathcal{F})
\end{equation}
has zeroes of total degree $d$ and away from $\Sigma$.
\end{itemize}
Let $\mathcal{M}^{\mu_\Sigma}_d:=\widetilde{\mathcal{M}}^{\mu_\Sigma}_d/Pic_X$.
\end{definition}

Now let $N=deg{\Sigma}$, and $\hat{Y}_{2d-N}$ be the moduli space of the effective divisor of degree $2d-N$ on $Y$, which is simply the $2d-N$-fold symmetric product scheme of $Y$. I want to define a map from $\mathcal{M}^{\mu_\Sigma}_d$ to $\hat{Y}_{2d-N}$. First define a divisor on $Y\times_k\widetilde{\mathfrak{S}}_\infty$, depending on our choice of $(\mu_f, \mu_\infty)$ as follows:
\begin{align}
\widetilde{D}:&=\sum_{x\in\Sigma_f, \mu_x=1}(y'_x, z_x)+\sum_{x\in\Sigma_f, \mu_x=-1}(y_x, z'_x)\\&+\sum_{x\in\Sigma_\infty, \mu_x=1}\sum_{i=1}^{d_x}(y^{(d_x+i)}_x, z^{(i)}_x)+\sum_{x\in\Sigma_\infty, \mu_x=-1}\sum_{i=1}^{d_x}(y^{(i)}_x, z^{(d_x+i)}_x)
\end{align}

And also let $D'_3=\sum_{x\in\Sigma}w'_x$, it is a divisor of $Y_3$. For the definition of $w_x$ and $w'_x$, look at section $(2.1)$ and section $(2.7)$.

Now as Howard-Shnidman did, pull back $\phi: \nu_{1,*}\mathcal{L}\longrightarrow\nu_{2,*}\mathcal{F}$ all the way up to $Y$, this map splits into the following:
\begin{equation}
\widetilde{\mathcal{L}}\oplus\tau_3^*\widetilde{\mathcal{L}}\xrightarrow{\phi}\widetilde{\mathcal{F}}\oplus\tau_3^*\widetilde{\mathcal{F}}
\end{equation}
in which $\widetilde{\phi}$ can be described more precisely:
\begin{align*}
\widetilde{\mathcal{L}}&\xrightarrow{\phi_{11}}\widetilde{\mathcal{F}}\\
\widetilde{\mathcal{L}}\xrightarrow{\sim}\tau_1^*\widetilde{\mathcal{L}}&\xrightarrow{\tau_1(\phi_{11})}\tau_1^*\widetilde{\mathcal{F}}\xrightarrow{\sim}\tau_3^*\widetilde{\mathcal{F}}\\
\tau^*_3\widetilde{\mathcal{L}}\xrightarrow{\sim}\tau^*_2\widetilde{\mathcal{L}}&\xrightarrow{\tau_2(\phi_{11})}\tau_2^*\widetilde{\mathcal{F}}\xrightarrow{\sim}\widetilde{\mathcal{F}}\\
\tau^*_3\widetilde{\mathcal{L}}&\xrightarrow{\tau_3(\phi_{11})}\tau^*_3\widetilde{\mathcal{F}}
\end{align*}
Now $\phi_{11}$ can be viewed as a section of $\widetilde{\mathcal{L}}^{-1}\otimes\widetilde{\mathcal{F}}$. 

Let $\Delta:= det(\nu_{1,*}\mathcal{L})^{-1}\otimes det(\nu_{2,*}\mathcal{F})$. It is a degree $d$ line bundle on $X\times_kS$ since it admits a section $det(\phi)$ of degree $d$ by definition. Then by observing that:
\begin{equation}
\nu_3^*\Delta=Nm(\widetilde{\mathcal{L}}^{-1}\otimes\widetilde{\mathcal{F}})
\end{equation}
one can see that the degree of $\widetilde{\mathcal{L}}^{-1}\otimes\widetilde{\mathcal{F}}$ is $2d$:
\begin{equation}
\nu^*\Delta\xrightarrow{\sim}\pi_3^*\nu_3^*\Delta=\pi_3^*Nm(\widetilde{\mathcal{L}}^{-1}\otimes\widetilde{\mathcal{F}})\xrightarrow{\sim}\widetilde{\mathcal{L}}^{-1}\otimes\widetilde{\mathcal{F}}\otimes\tau^*_3\widetilde{\mathcal{L}}^{-1}\otimes\tau_3^*\widetilde{\mathcal{F}}
\end{equation}
Of course $\nu^*\Delta$ has degree $4d$, so $\widetilde{\mathcal{L}}^{-1}\otimes\widetilde{\mathcal{F}}$ has degree $2d$ since $deg(\widetilde{\mathcal{L}}^{-1}\otimes\widetilde{\mathcal{F}}=deg(\tau^*_3\widetilde{\mathcal{L}}^{-1}\otimes\tau_3^*\widetilde{\mathcal{F}}$.

We also have to take care of the Iwahori level structure. we will see that this results in certain prescribed zeros of $\phi_{11}$:
\begin{enumerate}
\item At $x\in\Sigma_f$.
\begin{itemize}
\item $\mu_x=1$. In this case, the Iwahori level structures are defined by $\nu_{1,*}(\mathcal{L}(-y_x))$ and $\nu_{2,*}(\mathcal{F}(-z_x))$ respectively. Pull them back to the very top curve $Y$, one finds that:
\begin{align*}
\nu^*\nu_{1,*}(\mathcal{L}(-y_x))=\widetilde{\mathcal{L}}(-(y_x, z_x)-(y_x, z'_x))\oplus\tau_3^*\widetilde{\mathcal{L}}(-(y'_x, z_x)-(y'_x, z'_x))\\
\nu^*\nu_{2,*}(\mathcal{F}(-z_x))=\widetilde{\mathcal{F}}(-(y_x, z_x)-(y'_x, z_x))\oplus\tau_3^*\widetilde{\mathcal{F}}(-(y'_x, z'_x)-(y'_x, z'_x))
\end{align*}
Since $\phi$ preserves Iwahori level structure, $\phi_{11}$ should send $\widetilde{\mathcal{L}}(-(y_x, z_x)-(y_x, z'_x))$ to $\widetilde{\mathcal{F}}(-(y_x, z_x)-(y'_x, z_x))$, and by Galois equivariancy this is enough. Therefore one has:
\begin{equation}
\phi_{11}|_{(y'_x,z_x)}=0 \text{ and } \phi_{11}|_{(y_x, z_x)}\neq 0, \phi_{11}|_{(y'_x, z'_x)}\neq 0
\end{equation}
The nonvanishing part comes from the nonvanishing of $det(\phi)$ at $\Sigma$
\item $\mu_x=-1$. In this case, the Iwahori level structures are defined by $\nu_{1,*}(\mathcal{L}(-y'_x))$ and $\nu_{2,*}(\mathcal{F}(-z'_x))$ respectively. The same analysis as above yields:
\begin{equation}
\phi_{11}|_{(y_x,z'_x)}=0 \text{ and } \phi_{11}|_{(y_x, z_x)}\neq 0, \phi_{11}|_{(y'_x, z'_x)}\neq 0
\end{equation}
\end{itemize}
\item At $x\in\Sigma_\infty$.
\begin{itemize}
\item $\mu_x=1$. In this case, the Iwahori level structures are defined by:
\begin{align*}
\nu_{1,*}(\mathcal{L}(-y^{(1)}_x-...-y^{(d_x)}))\\
\nu_{2,*}(\mathcal{F}(-z^{(1)}_x-...-z^{(d_x)}))
\end{align*}
pull back them to the very top curve $Y$, one gets:
\begin{align*}
\widetilde{\mathcal{L}}(-(y^{(1)}_x, z^{(1)}_x)...-(y^{(d_x)}_x, z^{(d_x)}_x)-(y^{(1)}, z^{(d_x+1)}_x)...-(y^{(d_x)}, z^{(2d_x)}_x))\\
\oplus\tau_3^*\widetilde{\mathcal{L}}(-(y^{(d_x+1)}_x, z^{(1)}_x)...-(y^{(2d_x)}_x, z^{(d_x)}_x)...-(y^{(d_x+1)}_x, z^{(d_x+1)}_x)...-(y^{(2d_x)}_x, z^{(2d_x)}_x))\\
\widetilde{\mathcal{F}}(-(y^{(1)}_x, z^{(1)}_x)...-(y^{(d_x)}_x, z^{(d_x)}_x)-(y^{(d_x+1)}, z^{(1)}_x)...-(y^{(2d_x)}, z^{(d_x)}_x))\\
\oplus\tau_3^*\widetilde{\mathcal{F}}(-(y^{(1)}_x, z^{(d_x+1)}_x)...-(y^{(d_x)}_x, z^{(2d_x)}_x)...-(y^{(d_x+1)}_x, z^{(d_x+1)}_x)...-(y^{(2d_x)}_x, z^{(2d_x)}_x))
\end{align*}
Preserving Iwahori level structures translates to:
\begin{equation}
\phi_{11}|_{(y^{(d_x+1)}, z^{(1)}_x)...+(y^{(2d_x)}, z^{(d_x)}_x)}=0
\end{equation}
\item $\mu_x=-1$. In this case, the Iwahori level structure are defined by:
\begin{align*}
\nu_{1,*}(\mathcal{L}(-y^{(d_x+1)}_x-...-y^{(2d_x)}))\\
\nu_{2,*}(\mathcal{F}(-z^{(d_x+1)}_x-...-z^{(2d_x)}))
\end{align*}
The same analysis as above yields the following:
\begin{equation}
\phi_{11}|_{(y^{(1)}, z^{(d_x+1)}_x)...+(y^{(d_x)}, z^{(2d_x)}_x)}=0
\end{equation}
\end{itemize}
\end{enumerate}

Now you must see how the divisor $\widetilde{D}$ is defined. Due to the vanishing conditions, one can view $\phi_{11}$ as a section of the line bundle $(\widetilde{\mathcal{L}}^{-1}\otimes\widetilde{\mathcal{F}})(-\widetilde{D})$, which as degree $2d-N$. 
\begin{definition}
Define the map:
\begin{equation}
\mathcal{M}^{\mu_\Sigma}_d\longrightarrow\hat{Y}_{2d-N}
\end{equation}
as
\begin{equation}
(\mathcal{L}, \mathcal{F}, \phi)\longmapsto((\widetilde{\mathcal{L}}^{-1}\otimes\widetilde{\mathcal{F}})(-\widetilde{D}), \phi_{11})
\end{equation}
\end{definition}

\begin{remark}
Roughly speaking, this map is something to remember the divisor of $\phi_{11}$. Since $\phi_{11}$ must vanish at some prescribed points arising from the Iwahori level structure, one only has to remember the varying parts. This is the meaning of twisting it by $\widetilde{D}$.
\end{remark}

\subsection{The moduli space $\mathcal{A}_d$} The ``Hitchin type" fiberation needs a base space. I'm going to describe it here.
\begin{definition}
Define $\mathcal{A}_d$ to be the stack over $k$ whose functor of point on a $k$-scheme $S$ consisting of:
\begin{itemize}
\item A degree $d$ line bundle $\Delta$ over $X\times_kS$
\item A global section $a\in H^0(Y_3\times_kS, \nu^*_3\Delta)$ who vanishes at the support of $D'_3$, and $Tr(a)$ is a nonzero section of $\Delta$ who doesn't vanish at $\Sigma$.
\end{itemize}
\end{definition}
This is almost the same as the definition of Howard-Shnidman, the new things are the vanishing and nonvanishing conditions, as you know, coming from the Iwahori level structure.

Unlike the case considered by Yun-Zhang, this $\mathcal{A}_d$ is somehow simpler, since it is a scheme on the nose (\cite{HS}, section 3.4).
\begin{proposition}
Consider the trace map: 
\begin{align*}
Tr: \mathcal{A}_d&\longrightarrow\hat{X}_d\\
(\Delta, a)&\longmapsto (\Delta, \sigma_3\text{ invariant of }a+\sigma_3(a))
\end{align*}
This map is quasi-projective, and $\mathcal{A}_d$ is a quasi-projective scheme.
\end{proposition}

As in the case of $\mathcal{M}^{\mu_\Sigma}_d$, one can define the divisor memorizing map for $\mathcal{A}_d$.
\begin{definition}
Define the map:
\begin{equation}
\mathcal{A}_d\longrightarrow\hat{Y}_{3,2d-N}
\end{equation}
as
\begin{equation}
(\Delta, a)\longmapsto(\nu^*_3\Delta(-D'_3), a')
\end{equation}
In which $a'$ means viewing $a$ as a section of $\nu^*_3\Delta(-D'_3)$ rather than $\nu^*_3\Delta$
\end{definition}

Let $D\in Div^{+}(X-\Sigma)$ be an effective divisor on $X-\Sigma$. We can define $\mathcal{A}_D$ as the following fiber product:
\begin{definition}\label{def of AD}
\[\xymatrixcolsep{5pc}
\xymatrix{\mathcal{A}_D \ar[r] \ar[d] & \mathcal{A}_d \ar[d]^{Tr}\\
\text{Spec}k \ar[r]^{(\mathcal{O}_X(D),1)} & \hat{X}_d
}
\]
\end{definition}

From the definitions of $\mathcal{M}^{\mu_\Sigma}_d$ and $\mathcal{A}_d$, it is obvious that there exists a map from $\mathcal{M}^{\mu_\Sigma}_d$ to $\mathcal{A}_d$. 
\begin{definition}
The ``Hitchin type" fiberation for the shtuka side is the following morphism:
\begin{align}
f_d: \mathcal{M}^{\mu_\Sigma}_d&\longrightarrow\mathcal{A}_d\\
(\mathcal{L}, \mathcal{F}, \phi)&\longmapsto (\Delta:= det(\nu_{1,*}\mathcal{L})^{-1}\otimes det(\nu_{2,*}\mathcal{F}), Nm(\phi_{11}))
\end{align}
Preserving the Iwahori level structure makes sure the image lies in the right hand side. Also this map fits into the following diagram:
\begin{equation}\label{Mdiagram}
\xymatrix{\mathcal{M}^{\mu_\Sigma}_d \ar[r] \ar[d] &\hat{Y}_{2d-N} \ar[d]^{Nm}\\
\mathcal{A}_d\ar[r] & \hat{Y}_{3,2d-N}
}
\end{equation}
\end{definition}

This diagram is almost a Cartesian diagram. As we will see.

\subsection{The diagram above is almost Cartesion}
First for convenience, denote by $\mathcal{M}_d$ the Cartesian product, i.e.
\begin{equation}\label{m and a}
\xymatrix{\mathcal{M}_d \ar[r] \ar[d] &\hat{Y}_{2d-N} \ar[d]^{Nm}\\
\mathcal{A}_d\ar[r] & \hat{Y}_{3,2d-N}
}
\end{equation}

Since $\mathcal{A}_d$ is a scheme, $\mathcal{M}_d$ is a scheme. From Howard-Shnidman, one can know the following (\cite{HS}, Lemma 4.4):
\begin{proposition}
When $2d\geq2(2g_3-1)-1+N$, $\mathcal{M}_d$ and $\mathcal{A}_d$ are of dimension $2d-N-g+1$.
\end{proposition}

First Let $U_{2d-N}=\{E\in\hat{Y}_{2d-N}: Tr_{Y_3/X}(Nm_{Y/Y_3}(E))\cap\Sigma=\varnothing\}$, this is an open subscheme of $\hat{Y}_{2d-N}$. The space 
$\mathcal{M}_d$ actually also fits into the following Cartesian diagram:
\begin{equation*}
\xymatrix{\mathcal{M}_d\ar[r] \ar[d] & U_{2d-N} \ar[d]\\
Pic^{d}_X\ar[r] & Pic^{2d-N}_{Y_3}
}
\end{equation*}
The bottom horizontal arrow is simply pull back then twisted by the divisor $D'_3$. 

The right vertical arrow is the restriction to $U_{2d-N}$ of the map: 
\begin{equation*}
\hat{Y}_{2d-N}\xrightarrow{AJ}Pic^{2d-N}_Y\xrightarrow{Nm_{Y/Y_3}}Pic^{2d-N}_{Y_3}
\end{equation*}
in which $AJ$ is the Abel-Jacobi map, which is smooth when $2d-N\geq2(2g_3-1)-1$. And the smoothness of $Nm_{Y/Y_3}$ is shown by Howard-Shnidman. Therefore the right vertical map in the Cartesian diagram is smooth. From this, one can concludes that:
\begin{align*}
\text{dim}\mathcal{M}_{d}&=\text{dim}U_{2d-N}-\text{dim}Pic^{2d-N}_{Y_3}+\text{dim}Pic_X^d\\
&=2d-N-(g_3-1)+(g-1)\\
&=2d-N-(2g-2)+g-1=2d-N-g+1
\end{align*}
We used $g_3=2g-2$ since $Y_3$ is nonramified over $X$. Then from the defining diagram for $\mathcal{M}_d$ \eqref{m and a}, one knows $\mathcal{A}_d$ is of the same dimension. $\square$

By the definition of $\mathcal{M}^{\mu_\Sigma}$ one can see that it is actually a space over $\widetilde{\mathfrak{S}}_\infty$. Actually this is the only difference between $\mathcal{M}^{\mu_\Sigma}_d$ and $\mathcal{M}_d$:
\begin{proposition}
We have an isomorphism:
\begin{equation}
\mathcal{M}^{\mu_\Sigma}_d\xrightarrow{\sim}\mathcal{M}_d\times_k\widetilde{\mathfrak{S}}_\infty
\end{equation}
Therefore $\mathcal{M}^{\mu_\Sigma}_d$ is also a scheme, not just a Deligne-Mumford stack.
\end{proposition}

From the left to the right is simply the maps from $\mathcal{M}^{\mu_\Sigma}_d$ to $\hat{Y}_{2d-N}$, to $\mathcal{A}_d$ defined in the last subsection, together with remembering the map: $S\longrightarrow\widetilde{\mathfrak{S}}_\infty$. We only have to construct its inverse.

What does an $S$ point of $\mathcal{M}_d$ consists of? It is the following
\begin{itemize}
\item a degree $2d-N$ effective divisor $(\mathcal{K}',\psi')$ on $Y\times_kS$;
\item a degree $d$ line bundle on $X\times_kS$, say $\Delta$, and a section $a$ of $\nu_3^*\Delta$, vanishing at $D'_3$, and $Tr(a)$ doesn't vanish at $\Sigma$.
\end{itemize}
such that:
\begin{equation}
Nm((\mathcal{K}', \psi'))=(\nu^*_3\Delta(-D_3'), a')
\end{equation}
In which, again, $a'$ is just $a$ viewed as a section $\nu^*_3\Delta(-D_3')$.

The $\widetilde{\mathfrak{S}}_\infty$ structure allows us to consider the twists $(\mathcal{K}'(\widetilde{D}), \psi)$ in which $\psi$ is $\psi'$ viewed as a section of the more positive line bundle $\mathcal{K}'(\widetilde{D})$. Denote this more positive one by just $\mathcal{K}$. Then from the above equation we find the following equation:
\begin{equation}
\mathcal{K}\otimes\tau_3^*\mathcal{K}=\nu^*\Delta
\end{equation}
This provides $\mathcal{K}\otimes\tau_3^*\mathcal{K}$ with the $\tau_1$ and $\tau_2$ equivariant structure, since it is pulled back from the very bottom:
\begin{equation}
\mathcal{K}\otimes\tau_3^*\mathcal{K}\xrightarrow{\sim}\tau^*_1(\mathcal{K}\otimes\tau_3^*\mathcal{K})\xrightarrow{\sim}\tau^*_1\mathcal{K}\otimes\tau_2^*\mathcal{K}
\end{equation}
This in turn yields:
\begin{equation}
\tau_1^*\mathcal{K}^{-1}\otimes\mathcal{K}\xrightarrow{\sim}\tau^*_3\mathcal{K}^{-1}\otimes\tau_2^*\mathcal{K}\xrightarrow{\sim}\tau_2^*(\tau_1^*\mathcal{K}^{-1}\otimes\mathcal{K})
\end{equation}
Therefore $\tau_1^*\mathcal{K}^{-1}\otimes\mathcal{K}$ descent to $Y_2\times_kS$. Denote the descent by $\mathcal{M}$. 

Since one has the canonical $\tau_i, i=1,2,3$-equivariant isomorphism:
\begin{align*}
\pi_2^*(\mathcal{M}\otimes\sigma_2^*\mathcal{M})\xrightarrow{\sim}\tau_1^*\mathcal{K}^{-1}\otimes\mathcal{K}\otimes\tau^*_3(\tau_1^*\mathcal{K}^{-1}\otimes\mathcal{K})\\\xrightarrow{\sim}\tau^*_1\mathcal{K}^{-1}\otimes\mathcal{K}\otimes\tau^*_2\mathcal{K}^{-1}\otimes\tau_3^*\mathcal{K}\xrightarrow{\sim}\mathcal{O}_{Y\times_kS}
\end{align*}
The last equation comes from the left equation above, moving everything to the left hand side. This equivariant trivialization shows that $Nm_{{Y_2\times_kS}/{X\times_kS}}(\mathcal{M})$ is trivial. Therefore from the exact sequence:
\begin{equation}
1\xrightarrow{}Pic_{Y_2}/Pic_X\xrightarrow{id-\sigma_2} Pic_{Y_2}^0\xrightarrow{Nm}Pic_X^0\xrightarrow{}1
\end{equation}
one knows that there exists a unique line bundle $\mathcal{F}$ over $Y_2\times_kS$ up to twisted by line bundles pulling back from $X\times_kS$, such that 
\begin{equation}
\mathcal{F}\otimes\sigma^*\mathcal{F}^{-1}\xrightarrow{\sim}\mathcal{M}
\end{equation}
Now $\pi_2^*(\mathcal{K}^{-1}\otimes\mathcal{F})$ has a $\tau_1$ equivariant structure:
\begin{align*}
\tau^*_1(\mathcal{K}^{-1}\otimes\pi_2^*\mathcal{F})&\xrightarrow{\sim}\tau^*_1\mathcal{K}^{-1}\otimes\pi^*_2\sigma_2^*\mathcal{F}\\
&\xrightarrow{\sim}\tau^*_1\mathcal{K}^{-1}\otimes\mathcal{M}^{-1}\otimes\pi^*_2\mathcal{F}\\
&\xrightarrow{\sim}\tau^*_1\mathcal{K}^{-1}\otimes\tau^*_1\mathcal{K}\otimes\mathcal{K}^{-1}\otimes\pi^*_2\mathcal{F}\\
&\xrightarrow{\sim}\mathcal{K}^{-1}\otimes\pi_2^*\mathcal{F}
\end{align*}
This produces the desired line bundle $\mathcal{L}$ on $Y_1\otimes_kS$ such that $\pi_1^*\mathcal{L}\xrightarrow{\sim}\mathcal{K}^{-1}\otimes\pi_2^*\mathcal{F}$. 

As in the previous subsections, let $\widetilde{\mathcal{L}}=\pi_1^*\mathcal{L}$ and $\widetilde{\mathcal{F}}=\pi^*_2\mathcal{F}$, one gets
\begin{equation}
\widetilde{\mathcal{L}}^{-1}\otimes\widetilde{\mathcal{F}}\xrightarrow{\sim}\mathcal{K}
\end{equation}
Via this isomorphism, one can view $\psi$ as a global section of $\widetilde{\mathcal{L}}^{-1}\otimes\widetilde{\mathcal{F}}$. Note that $(\mathcal{K}, \psi)$ comes from $(\mathcal{K}'(\widetilde{D}), \psi)$, one can also view $\psi'$ as a global section of $\widetilde{\mathcal{L}}^{-1}\otimes\widetilde{\mathcal{F}}(-\widetilde{D})$. In other words, one can view $a$ as:
\begin{equation}
\psi: \widetilde{\mathcal{L}}\longrightarrow\widetilde{\mathcal{K}}
\end{equation}
vanishing at the divisor $\widetilde{D}$.

By applying $\tau_i$s to it, one can produce:
\begin{align*}
\widetilde{\mathcal{L}}&\xrightarrow{\psi}\widetilde{\mathcal{F}}\\
\widetilde{\mathcal{L}}\xrightarrow{\sim}\tau_1^*\widetilde{\mathcal{L}}&\xrightarrow{\tau_1(\psi)}\tau_1^*\widetilde{\mathcal{F}}\xrightarrow{\sim}\tau_3^*\widetilde{\mathcal{F}}\\
\tau^*_3\widetilde{\mathcal{L}}\xrightarrow{\sim}\tau^*_2\widetilde{\mathcal{L}}&\xrightarrow{\tau_2(\psi)}\tau_2^*\widetilde{\mathcal{F}}\xrightarrow{\sim}\widetilde{\mathcal{F}}\\
\tau^*_3\widetilde{\mathcal{L}}&\xrightarrow{\tau_3(\psi)}\tau^*_3\widetilde{\mathcal{F}}
\end{align*}
This matrix descent to a morphism:
\begin{equation}
\phi: \nu_{1,*}\mathcal{L}\longrightarrow\nu_{2,*}\mathcal{F}
\end{equation}
It is easy to check that $det(\phi)=Tr(a)$. The fact that $\psi$ vanishes at $\widetilde{D}$ makes sure that $\phi$ preserves the Iwahori level structure, and the nonvanishing of $Tr(a)$ at $\Sigma$ makes sure the nonvanishing of $det(\phi)$ at $\Sigma$. This finishes the construction of the inverse map. One can check that this is really the inverse. $\square$

\subsection{The Hecke correspondence for $\mathcal{M}_d^{\mu_{\Sigma}}$ and $\mathcal{M}_d$}
You may have noticed that in the definition of $\mathcal{M}_d^{\mu_\Sigma}$, the $\phi$ comes from the vertical Hecke symmetries of the moduli of shtukas. Now we take care of the horizontal Hecke modifications. As before, fix our choice of $(\underline{\mu},\mu_f,\mu_\infty)$. Define the following stack:
\begin{definition}
Let $\widetilde{Hk}_{\mathcal{M}_d}^{\underline{\mu}}$ be the stack over $k$ whose functor of points on a $k$-scheme $S$ consists of the following:
\begin{itemize}
\item $S\longrightarrow\widetilde{\mathfrak{S}}_\infty$
\item $r$ $S$-points of $Y_1: y_1, y_2, ..., y_r: S\longrightarrow Y_1$ and $S$-points of $Y_2: z_1,z_2, ...,z_r: S\longrightarrow Y_2$\\
such that $\nu_1\circ y_i=\nu_2\circ z_i$ as $S$ points of $X$, for $i=1,...r$. These two points is equivalent to $r$ $S$-points of $Y=Y_1\times_XY_2$.
\item $r+1$ $S$-objects of $\mathcal{M}_d^{\mu_{\Sigma}}$: $(\mathcal{L}_0, \mathcal{F}_0, \phi_0), ...,(\mathcal{L}_r, \mathcal{F}_r, \phi_r)$.
\item Morphisms of coherent sheaves: $f_i: \mathcal{L}_{i-1}\dashrightarrow\mathcal{L}_{i}$ and $f'_i: \mathcal{F}_{i-1}\dashrightarrow\mathcal{F}_{i}$ such that their push foward to $X\times_kS$ fit into the following commutative diagram:
\begin{equation}
\xymatrix{\nu_{1,*}\mathcal{L}_0 \ar@{-->}[r]^{\nu_{1,*}{f_1}} \ar[d]^{\phi_0} & \nu_{1,*}\mathcal{L}_1 \ar@{-->}[r]^{\nu_{1,*}f_2} \ar[d]^{\phi_1} &......\ar@{-->}[r]^{\nu_{1,*}f_r} & \nu_{1,*}\mathcal{L}_r \ar[d]^{\phi_r}\\ 
\nu_{2,*}\mathcal{F}_0 \ar@{-->}[r]^{\nu_{2,*}f'_1} & \nu_{2,*}\mathcal{F}_1 \ar@{-->}[r]^{\nu_{2,*}f'_2} &......\ar@{-->}[r]^{\nu_{2,*}f'_r} & \nu_{2,*}\mathcal{F}_r}
\end{equation}
\end{itemize}
\end{definition}
As usual, let $Hk_{\mathcal{M}_d}^{\underline{\mu}}$ be $\widetilde{Hk}_{\mathcal{M}_d}^{\underline{\mu}}/Pic_X$. 

It is a self-correspondence of $\mathcal{M}_d^{\mu_\Sigma}$ over $\mathcal{A}_d$ (actually even over $\mathcal{A}_d\times_k\widetilde{\mathfrak{S}}_\infty$, though we don't use it)
\begin{equation}
\xymatrix{ & Hk_{\mathcal{M}_d}^{\underline{\mu}} \ar[ld]_{\overleftarrow{p}} \ar[rd]^{\overrightarrow{p}} & \\
\mathcal{M}_d^{\mu_{\Sigma}} \ar[rd] & & \mathcal{M}_d^{\mu_{\Sigma}} \ar[ld]\\
&\mathcal{A}_d\times_k\widetilde{\mathfrak{S}}_\infty &}
\end{equation}

The case $\mu=+1$ is the building block of the mutiple modification points ones, i.e., one has the following isomorphism as a self-correspondence of $\mathcal{M}_d$ over $\mathcal{A}_d$:
\begin{equation}
Hk_{\mathcal{M}_d}^{\underline{\mu}}=Hk_{\mathcal{M}_d}^1\times_{\overrightarrow{p}, \mathcal{M}_d, \overleftarrow{p}}Hk_{\mathcal{M}_d}^1\times_{\overrightarrow{p}, \mathcal{M}_d, \overleftarrow{p}}......\times_{\overrightarrow{p}, \mathcal{M}_d, \overleftarrow{p}}Hk_{\mathcal{M}_d}^1
\end{equation}

The purpose now is to transfer the Hecke correspondence for $\mathcal{M}_d^{\mu_\Sigma}$ to a Hecke correspondence for $\mathcal{M}_d$. But let's first see what that looks like. Recall a definition made by Howard-Shnidman (\cite{HS} section 4.3).
\begin{definition}
Let $H_{2d-N}$ be the $k$-stack whose typical $S$-objects are the following:
\begin{itemize}
\item Two $S$-points of $\hat{Y}_{2d-N}: (\mathcal{K}'_0, \psi'_0), (\mathcal{K}_1', \psi'_1)$.
\item An $S$-point of $Y$, i.e. $v: S\longrightarrow Y$.
\item An isomorphism of sheaves: $s: \mathcal{K}'_0(\tau_2(v)-\tau_1(v))\xrightarrow{\sim}\mathcal{K}'_1$, such that $s(\psi_0')=\psi_1'$, in which vewing $a'_0$ as a rational section of $\mathcal{K}'_0(\tau_2(v)-\tau_1(v))$.
\end{itemize}
\end{definition}
The last condition guarantees that $Nm(\mathcal{K}'_0, \psi'_0)=Nm(\mathcal{K}_1', \psi'_1)$, therefore one gets the following diagram:
\begin{equation}
\xymatrix{ & H_{2d-N} \ar[ld]_{\gamma_0} \ar[rd]^{\gamma_1} & \\
\hat{Y}_{2d-N} \ar[rd] & & \hat{Y}_{2d-N} \ar[ld]\\
&\hat{Y}_{3,2d-N} &}
\end{equation}

Recall that $\mathcal{M}_d$ is defined to be the fiber product: $\mathcal{M}_d:=\hat{Y}_{2d-N}\times_{\hat{Y}_{3, 2d-N}}\mathcal{A}_d$. Pull back the above self-correspondence of $\hat{Y}_{2d-N}$ to a self-correspondence of $\mathcal{M}_d$ over $\mathcal{A}_d$, denoted by $\mathcal{H}_d$:
\begin{equation}
\xymatrix{ & \mathcal{H}_d \ar[ld]_{\gamma_0} \ar[rd]^{\gamma_1} & \\
\mathcal{M}_d \ar[rd] & & \mathcal{M}_d \ar[ld]\\
&\mathcal{A}_d &}
\end{equation}

We know the dimension of $\mathcal{H}_d$:
\begin{proposition}
$\gamma_0$ and $\gamma_1$ are finite and surjective. $\mathcal{H}_d$ is a scheme of dimension $2d-N-g+1$.
\end{proposition}

You may wonder how this $\mathcal{H}_d$ is so defined, the next proposition just tells you this:
\begin{proposition}
$\mathcal{H}_d\times_k\widetilde{\mathfrak{S}}_\infty$ is isomorphic to $Hk^1_{\mathcal{M}_d}$ as a self-correspondence over $\mathcal{A}\times_k\widetilde{\mathfrak{S}}_\infty$ if one identifies $\mathcal{M}_d\times_k\widetilde{\mathfrak{S}}_\infty$ with $\mathcal{M}_d^{\mu_\Sigma}$ via the proposition in the last subsection.
\end{proposition}

Going from $Hk^1_{\mathcal{M}_d}$ to $\mathcal{H}_d\times_k\widetilde{\mathfrak{S}}_\infty$ actually is the definition of $\mathcal{H}_d$. Starts with the point of $Hk^1_{\mathcal{M}_d}$, say: 
\begin{equation*}
S\longrightarrow\widetilde{\mathfrak{S}}_\infty, (\mathcal{L}_0, \mathcal{F}_0, \phi_0), (\mathcal{L}_1, \mathcal{F}_1, \phi_1), f, f', (y,z)
\end{equation*}
one gets two $S$-point of $\mathcal{M}_d$:
\begin{equation*}
(\widetilde{\mathcal{L}}_0^{-1}\otimes\widetilde{\mathcal{F}}_0(-\widetilde{D}), \psi'_0), (\widetilde{\mathcal{L}}_1^{-1}\otimes\widetilde{\mathcal{F}}_1(-\widetilde{D}), \psi_1')
\end{equation*}
Pulling back $f$ and $f'$ to $Y$, one finds:
\begin{align*}
\widetilde{f}:\widetilde{\mathcal{L}}_0((y,z)+(y,z'))\xrightarrow{\sim}\widetilde{\mathcal{L}}_1\\
\widetilde{f'}:\widetilde{\mathcal{F}}_0((y,z)+(y',z))\xrightarrow{\sim}\widetilde{\mathcal{F}}_1
\end{align*}
from this one can get:
\begin{equation*}
\widetilde{f}\otimes\widetilde{f'}^{-1}: \widetilde{\mathcal{L}}_0^{-1}\otimes\widetilde{\mathcal{F}}_0((y',z)-(y,z'))\xrightarrow{\sim}\widetilde{\mathcal{L}}_1^{-1}\otimes\widetilde{\mathcal{F}}_1
\end{equation*}
Therefore the following tuple:
\begin{equation*}
(\widetilde{\mathcal{L}}_0^{-1}\otimes\widetilde{\mathcal{F}}_0(-\widetilde{D}), \psi'_0), (\widetilde{\mathcal{L}}_1^{-1}\otimes\widetilde{\mathcal{F}}_1(-\widetilde{D}), \psi_1'), (\widetilde{f}\otimes\widetilde{f'}^{-1})(-\widetilde{D})
\end{equation*}
 lies in $\mathcal{H}_d$.

Like in the proof of the equivalence between $\mathcal{M}_d^{\mu_\Sigma}$ and $\mathcal{M}_d\times_k\widetilde{\mathfrak{S}}_\infty$, the important part is going backwards. From that equivalence, starting from a point $\hat{Y}_{2d-N}: (\mathcal{K}'_0, \psi'_0), (\mathcal{K}_1', \psi'_1)$ on $(\mathcal{M}_d\times_{\mathcal{A}_d}\mathcal{M}_d)\times_k\widetilde{\mathfrak{S}}_\infty$, one can produce a point on $\mathcal{M}_d^{\mu_\Sigma}\times_{\mathcal{A}_d\times_k\widetilde{\mathfrak{S}}_\infty}\mathcal{M}_d^{\mu_\Sigma}$:
\begin{equation*}
(\mathcal{L}_0, \mathcal{F}_0, \phi_0), (\mathcal{L}_1, \mathcal{F}_1, \phi_1)
\end{equation*}
Now the only thing need to do is to prove that giving an $S$ point $v=(y,z)$ of $Y$ and $s$ identifying $(\mathcal{K}'_0(\tau_2(v)-\tau_1(v)), \psi'_0)$ with $(\mathcal{K}'_1, \psi'_1)$, one can reconstruct the following:
\begin{equation*}
f: \mathcal{L}_0(y)\xrightarrow{\sim}\mathcal{L}_1 \text{ and } f':  \mathcal{F}_0(z)\xrightarrow{\sim}\mathcal{F}_1
\end{equation*}
To do this, as mentioned by Howard-Shnidman, one has to go back to the prove of the equivalence.

As in the prove of the equivalence, $S\longrightarrow\widetilde{\mathfrak{S}}_\infty$ allows one to define: $\mathcal{K}_i=\mathcal{K}'_i(\widetilde{D})$. Apply $\tau_1^*$ to the following the equation and taking the inverse:
\begin{equation*}
s: \mathcal{K}_0((y',z)-(y,z'))\xrightarrow{\sim}\mathcal{K}_1
\end{equation*}
one gets:
\begin{equation*}
\tau_1(s):\tau_1^*\mathcal{K}^{-1}_0((y,z)-(y',z'))\xrightarrow{\sim}\tau^*_1\mathcal{K}^{-1}_1
\end{equation*}
Tensor them up:
\begin{equation*}
\tau_1^*\mathcal{K}^{-1}_0\otimes\mathcal{K}_0((y',z)-(y,z')+ (y,z)-(y',z'))\xrightarrow{\sim}\tau^*_1\mathcal{K}^{-1}_1\otimes\mathcal{K}_1
\end{equation*}
Descent to $Y_2$, one finds:
\begin{equation*}
\mathcal{M}_0(z-z')\xrightarrow{\sim}\mathcal{M}_1
\end{equation*}
By construction, $\mathcal{F}_1$ is the unique line bundle up twists by $Pic_X$ such that $\mathcal{F}_1\otimes\sigma_2^{*}\mathcal{F}_1\xrightarrow{\sim}\mathcal{M}_1$, but the above equation shows that $\mathcal{F}_0(z)$ is another such. This way one gets $f':  \mathcal{F}_0(z)\xrightarrow{\sim}\mathcal{F}_1$. 

Now recall how we reconstructed $\mathcal{L}$. $\mathcal{L}_i$ is defined to be a line bundle on $Y_1\otimes_kS$ pullback to be $\mathcal{K}_i^{-1}\otimes\widetilde{\mathcal{F}}_i$. Using the isomorphism $s$ and $f'$:
\begin{align*}
\widetilde{\mathcal{L}}_1&\xrightarrow{\sim}\mathcal{K}_0((y,z)-(y',z))\otimes\widetilde{\mathcal{F}}_0((y,z')+(y',z))\\
&\xrightarrow{\sim}\mathcal{K}_0\otimes\widetilde{\mathcal{F}}_0((y,z')+(y,z))\\
&\xrightarrow{\sim}\pi^*_1({\mathcal{L}_0(y)})
\end{align*}
One can check that this isomorphism is actually $\tau_1$-equivariant, therefore descent to an isomorphism: $f: \mathcal{L}_0(y)\xrightarrow{\sim}\widetilde{\mathcal{L}}_1$ (I have taken the reverse direction).

The Iwahori level structure part is already built up in the proof of the equivalence of $\mathcal{M}_d^{\mu_\Sigma}$ and $\mathcal{M}_d\times_k\widetilde{\mathfrak{S}}_\infty$. $\square$

Similarly one can form the multiple Hecke correspondence for $\mathcal{M}_d: \mathcal{H}_d^{\underline{\mu}}$. The same proof as above can show that there is an isomorphism of self-correspondences of $\mathcal{M}_d^{\mu_\Sigma}$ over $\mathcal{A}_d$, when $\mathcal{M}_d^{\mu_\Sigma}$ is identified with $\mathcal{M}_d\times_k\widetilde{\mathfrak{S}}_\infty$:
\begin{equation}
Hk^{\underline{\mu}}_{\mathcal{M}_d}\xrightarrow{\sim}\mathcal{H}_d^{\underline{\mu}}\times_k\widetilde{\mathfrak{S}}_\infty
\end{equation}
and $\mathcal{H}_d$ is actually the building block of $\mathcal{H}_d^{\underline{\mu}}$:
\begin{equation}
\mathcal{H}_d^{\underline{\mu}}\xrightarrow{\sim}\mathcal{H}_d\times_{\gamma_1,\mathcal{M}_d,\gamma_0}\mathcal{H}_d\times_{\gamma_1,\mathcal{M}_d,\gamma_0}...\times_{\gamma_1,\mathcal{M}_d,\gamma_0}\mathcal{H}_d
\end{equation}

\section{Going towards the intersection number}

\subsection{The master diagram} The master diagram plays a key role in Yun-Zhang's theory. We first review an auxilary moduli space appearing in this diagram, namely $H_d$.

\begin{definition}
Define $\widetilde{H}_d(\Sigma)$ to be the $k$-stack whose $S$-points consist of the following:
\begin{itemize}
\item two $S$-points of $Bun_G(\Sigma)$, say $\mathcal{E}^{\dagger}=(\mathcal{E}, \{\mathcal{E}(-\frac{1}{2}x)\})$ and $\mathcal{E'}^{\dagger}=(\mathcal{E'}, \{\mathcal{E'}(-\frac{1}{2}x)\})$, over$X\times_kS$; 
\item a morphism of coherent sheaves $\phi: \mathcal{E}\longrightarrow\mathcal{E}'$ preserving the Iwahori level structure, i.e. sending $\{\mathcal{E}(-\frac{1}{2}x)\}$ to $\{\mathcal{E'}(-\frac{1}{2}x)\}$. Also, as before, $det(\phi)$ has degree $d$ and its vanishing locus is away from $\Sigma$.
\end{itemize}
\end{definition}
As you know, define $H_d(\Sigma):=\widetilde{H}_d(\Sigma)/Pic_X$, and let:
\begin{equation*}
\overleftrightarrow{p}_H=(\overleftarrow{p}_H, \overrightarrow{p}_H): H_d(\Sigma)\longrightarrow Bun_G(\Sigma)^2.
\end{equation*}
be the two projections. Yun and Zhang pointed out the following (\cite{HS}, Definition 5.7):
\begin{proposition}
The two projections, $\overleftarrow{p}_H, \overrightarrow{p}_H$ are representable and smooth of relative dimension $2d$. Therefore, $H_d(\Sigma)$ is a smooth Artin stack of pure dimension $2d+3(g-1)+N$.
\end{proposition}

From this definition you can see that the above $\mathcal{M}_d^{\mu_\Sigma}$ is simply the following fiber product:
\begin{equation}
\xymatrix{\mathcal{M}^{\mu_\Sigma}_d \ar[r] \ar[d] &H_d(\Sigma) \ar[d]^{(\overleftarrow{p}, \overrightarrow{p})}\\
Bun_{T_1}\times_kBun_{T_2}\times_k\widetilde{\mathfrak{S}}_\infty\ar[r] & Bun_G(\Sigma)}
\end{equation}

Since the definition of shtukas involves Atkin-Lehner involution, here we also introduce the Atkin-Lehner involution for $H_d(\Sigma)$ and $\mathcal{M}^{\mu_\Sigma}_d$. Recall we have the Atkin-Lehner involution for $Bun_G(\Sigma)$,  namely $AL_{G,\infty}$. For $H_d(\Sigma)$, its Atkin-Lehner involution $AL_{H_\infty}$:
\begin{align}
H_d(\Sigma)\times_k\mathfrak{S}_\infty&\longrightarrow H_d(\Sigma)\\
(\mathcal{E}^{\dagger}\xrightarrow{\phi}\mathcal{E}'^{\dagger}, S\longrightarrow\mathfrak{S}_\infty)&\longmapsto (AL_{G,\infty}(\mathcal{E}^{\dagger})\xrightarrow{\phi}AL_{G,\infty}(\mathcal{E}^{\dagger}))
\end{align}
Recall the relation \eqref{3.33}, define the Atkin-Lehner involution $AL_{\mathcal{M},\infty}$ in the following way:
\begin{align}
\mathcal{M}^{\mu_\Sigma}_d&\longrightarrow \mathcal{M}^{\mu_\Sigma}_d\\
(\mathcal{L}, \mathcal{F}, \phi, S\longrightarrow\widetilde{\mathfrak{S}}_\infty)&\longmapsto (AL_{T_1,\mu_\infty}(\mathcal{L}), AL_{T_2,\mu_\infty}(\mathcal{F}), \phi, S\xrightarrow{Fr_S}S\longrightarrow\widetilde{\mathfrak{S}}_\infty)
\end{align}

Of course as you can see, $H_d(\Sigma)$ is modeled on the vertical Hecke symmetries of moduli shtukas. Considering the horizontal Hecke modifications lead one to the following definition:
\begin{definition}
Let's define $\widetilde{Hk}_{H_d}^{\underline{\mu}}(\Sigma)$ to be the stack whose $S$-points consist of:
\begin{itemize}
\item $r$-$S$ points of $X$, say $x_1,..., x_r$
\item a commutative diagram:
\begin{equation}
\xymatrix{\mathcal{E}_0 \ar@{-->}[r]^{f_1} \ar[d]^{\phi_0} & \mathcal{E}_1 \ar@{-->}[r]^{f_2} \ar[d]^{\phi_1} &......\ar@{-->}[r]^{f_r} & \mathcal{E}_r \ar[d]^{\phi_r}\\ 
\mathcal{E}'_0 \ar@{-->}[r]^{f'_1} & \mathcal{E}'_1 \ar@{-->}[r]^{f'_2} &......\ar@{-->}[r]^{f'_r} & \mathcal{E}'_r}
\end{equation}
in which each column is a point of $\widetilde{H}_d(\Sigma)(S)$.  $f_i$ and $f'_i$ are isomorphisms away from $\Gamma_{x_i}$ for $ i=1,...,r$, and have relative position $\mu_i$.
\end{itemize}
\end{definition}
Let $Hk^{\underline{\mu}}_{H_d}(\Sigma):=\widetilde{Hk}_{H_d}^{\underline{\mu}}(\Sigma)/Pic_X$. One has:
\begin{proposition}
$Hk^{\underline{\mu}}_{H_d}(\Sigma)$ is an Artin stack of dimension $2d+2r+3(g-1)+N$
\end{proposition}

Now we are ready for Yun-Zhang's master diagram, adapted to our current case:
\begin{equation}
\xymatrixcolsep{5pc}\xymatrix{Hk^{\underline{\mu}}_{T_1}\times Hk^{\underline{\mu}}_{T_2}\times\widetilde{\mathfrak{S}}_\infty \ar[r]^{\theta_{1,Hk}\times\theta_{2,Hk}\times Id_{\widetilde{\mathfrak{S}}_\infty}} \ar[d]^{p^{\underline{\mu}}_{T_1,0}\times p^{\underline{\mu}}_{T_2,0}\times Id_{\widetilde{\mathfrak{S}}_\infty}, \alpha_T} & Hk^r_G(\Sigma)^2\times\widetilde{\mathfrak{S}}_\infty \ar[d]^{p_{G,0}^2, \alpha_G} & \ar[l]_{\overleftrightarrow{q}\times Id_{\widetilde{\mathfrak{S}}_\infty}} Hk^r_{H_d}(\Sigma)\times\widetilde{\mathfrak{S}}_\infty \ar[d]^{p_{H,0}, \alpha_H}\\
(Bun_{T_1}\times Bun_{T_2}\times\widetilde{\mathfrak{S}}_\infty)^2 \ar[r]^{(\theta^{\mu_\Sigma}_{1,Bun}\times\theta^{\mu_\Sigma}_{2,Bun})^2}& Bun_G(\Sigma)^2\times Bun_G(\Sigma)^2 &\ar[l]_{\overleftrightarrow{p_H}\times\overleftrightarrow{p_H}} H_d(\Sigma)\times H_d(\Sigma)\\
Bun_{T_1}\times Bun_{T_2}\times\widetilde{\mathfrak{S}}_\infty \ar[r]^{\theta^{\mu_\Sigma}_{1,Bun}\times\theta^{\mu_\Sigma}_{2,Bun}} \ar[u]_{id, Fr} & Bun_G(\Sigma)^2 \ar[u]_{id, Fr} & \ar[l]_{\overleftrightarrow{p_H}} \ar[u]_{id,Fr} H_d(\Sigma)
}
\end{equation}
Let's review the meaning of those $\alpha$s appeared above.
\begin{enumerate}
\item
\begin{align}
\alpha_T: Hk^{\underline{\mu}}_{T_1}\times Hk^{\underline{\mu}}_{T_2}\times\widetilde{\mathfrak{S}}_\infty&\xrightarrow{p^{\underline{\mu}}_{T_1,r}\times p^{\underline{\mu}}_{T_1,r}\times id_{\widetilde{\mathfrak{S}}_\infty}}Bun_{T_1}\times Bun_{T_2}\times{\widetilde{\mathfrak{S}}_\infty}\\
&\xrightarrow{AL_{T_1, T_2, \mu_\infty}}Bun_{T_1}\times Bun_{T_2}\times{\widetilde{\mathfrak{S}}_\infty}
\end{align}
In which $AL_{T_1, T_2, \mu_\infty}$ sends $(\mathcal{L}, \mathcal{F}, v^{(1)}:S\longrightarrow\widetilde{\mathfrak{S}}_\infty)$ to:
\begin{equation}
(\mathcal{L}(-\sum_{x\in\Sigma_\infty}\mu_xy^{(1)}_x), \mathcal{F}(-\sum_{x\in\Sigma_\infty}\mu_xz^{(1)}_x), v^{(2)}: S\xrightarrow{Fr}S\xrightarrow{v^{(1)}}\widetilde{\mathfrak{S}}_\infty)
\end{equation}
Notice that there is a Frobenius shift on the last factor. From the master diagram, this is obvious.
\item
$\alpha_G: Hk_G^r(\Sigma)^2\times\widetilde{\mathfrak{S}}_\infty\xrightarrow{p_{G,r}^2\times id_{\widetilde{\mathfrak{S}}_\infty}}Bun_G(\Sigma)^2\times\widetilde{\mathfrak{S}}_\infty\xrightarrow{AL^2_{G,\infty}}Bun_G(\Sigma)^2$
\item
$\alpha_{H}: Hk_{H_d}^r(\Sigma)\times\widetilde{\mathfrak{S}}_\infty\xrightarrow{p_{H,r}^2\times id_{\widetilde{\mathfrak{S}}_\infty}}H_d(\Sigma)\times\widetilde{\mathfrak{S}}_\infty\xrightarrow{AL_{H,\infty}}H_d(\Sigma)$
\end{enumerate}

How this diagram plays its role in Yun-Zhang's theory? Because taking the fiber products of the columns is how we reached the ``thing we want to compute", taking the fiber product of the rows is how we are going to compute it.

The fiber product of the columns:
\begin{equation}\label{fibcol}
Sht^{\underline{\mu}}_{T_1}(\mu_\infty\cdot\Sigma_{1,\infty})'\times_{\widetilde{\mathfrak{S}}_\infty}Sht^{\underline{\mu}}_{T_2}(\mu_\infty\cdot\Sigma_{2,\infty})'\xrightarrow{\theta^\mu_1\times\theta^\mu_2}Sht^{\underline{\mu}}_G(\Sigma;\Sigma_{\infty})'\times_{\widetilde{\mathfrak{S}}_\infty}Sht^{\underline{\mu}}_G(\Sigma;\Sigma_{\infty})'\xleftarrow{}Sht^{\underline{\mu}}_{H_d}(\Sigma;\Sigma_{\infty})'
\end{equation}
In which:
\begin{equation}
\xymatrix{Sht^{\underline{\mu}}_{H_d}(\Sigma;\Sigma_{\infty})\ar[r]\ar[d] &  Hk_{H_d}^r(\Sigma)\times\mathfrak{S}_\infty \ar[d]^{p_{H,0}, \alpha_H}\\
H_d(\Sigma) \ar[r]^-{id, Fr} & H_d(\Sigma)\times H_d(\Sigma)
}
\end{equation}
and $Sht^{\underline{\mu}}_{H_d}(\Sigma;\Sigma_{\infty})'=Sht^{\underline{\mu}}_{H_d}(\Sigma;\Sigma_{\infty})\times_{\mathfrak{S}_\infty}\widetilde{\mathfrak{S}}_\infty$.

This looks familiar right? Actually by taking the fiber product of the two arrows one obtains all $\ast$s in \eqref{product} for all $D\in Div^{+}(X-\Sigma)$ of degree $d$ collectively. To see it, just observe that one has the following decomposition:
\begin{equation}\label{decom}
Sht^{\underline{\mu}}_{H_d}(\Sigma;\Sigma_{\infty})=\coprod_{D\in Div^d(X-\Sigma)(k), degD=d}Sht^{\underline{\mu}}_G(\Sigma;\Sigma_{\infty}, h_D)
\end{equation}
Here for use, base change this equality to $\widetilde{\mathfrak{S}}_\infty$ with respect to $\mathfrak{S}_\infty$:
\begin{equation*}
Sht^{\underline{\mu}}_{H_d}(\Sigma;\Sigma_{\infty})'=\coprod_{D\in Div^d(X-\Sigma)(k), degD=d}Sht^{\underline{\mu}}_G(\Sigma;\Sigma_{\infty}, h_D)'
\end{equation*}

Now taking the fiber product of each row, what you see? The following:
\begin{equation}
\xymatrix{Hk^{\underline{\mu}}_{\mathcal{M}_d}\ar[d]^{p_{\mathcal{M},0}, \alpha_M}\\
\mathcal{M}_d^{\mu_\Sigma}\times\mathcal{M}_d^{\mu_\Sigma}\\
\ar[u]^{id, Fr} \mathcal{M}_d^{\mu_\Sigma}
}
\end{equation}
in which: $\alpha_\mathcal{M}:=AL_{\mathcal{M},\infty}\circ p_{\mathcal{M},r}$.

Let the fiber product of these two arrows be:
\begin{definition}\label{shtm}
\begin{equation}
\xymatrix{Sht^{\underline{\mu}}_{\mathcal{M}_d} \ar[r] \ar[d] & Hk^{\underline{\mu}}_{\mathcal{M}_d} \ar[d]^{p_{\mathcal{M},0}, \alpha_\mathcal{M}}\\
\mathcal{M}_d^{\mu_\Sigma} \ar[r]^-{id, Fr} & \mathcal{M}_d^{\mu_\Sigma}\times\mathcal{M}_d^{\mu_\Sigma}
}
\end{equation}
\end{definition}

Of course, $Sht^{\underline{\mu}}$ should coincide with \eqref{fibcol}. Then \eqref{decom} can induce the following decomposition of $Sht^{\underline{\mu}}_{\mathcal{M}_d}$:
\begin{equation}
Sht^{\underline{\mu}}_{\mathcal{M}_d}=\coprod_{D\in Div^+(X-\Sigma)(k), degD=d}Sht^{\underline{\mu}}_{\mathcal{M},D}
\end{equation} 
in which $Sht^{\underline{\mu}}_{\mathcal{M},D}$ means the preimage of $Sht^{\underline{\mu}}_{H_d}(\Sigma;\Sigma_{\infty}, h_D)'$ in $Sht^{\underline{\mu}}_{\mathcal{M},D}$.

The expectation is the following:
\begin{equation}\label{eq1}
(id, Fr)^![Hk^r_{H_d}(\Sigma)\times\widetilde{\mathfrak{S}}_\infty]=[Sht^{\underline{\mu}}_{H_d}(\Sigma;\Sigma_{\infty})']
\end{equation}
 then one has the following equation:
\begin{equation}
(\theta'^\mu_1\times\theta'^\mu_2)^![Sht^{\underline{\mu}}_{H_d}(\Sigma;\Sigma_{\infty})']=(\theta'^\mu_1\times\theta'^\mu_2)^!(id, Fr)^![Hk^r_{H_d}(\Sigma)\times\widetilde{\mathfrak{S}}_\infty]
\end{equation}

If one can check those conditions listed by Yun-Zhang in their first volumn concerning the master diagram, one should get:
\begin{equation}
(\theta'^\mu_1\times\theta'^\mu_2)^!(id, Fr)^![Hk^r_{H_d}(\Sigma)\times\widetilde{\mathfrak{S}}_\infty]=(id, Fr)^!(\theta_{1,Hk}\times\theta_{2,Hk}\times id_{\widetilde{\mathfrak{S}}_\infty})^![Hk^r_{H_d}(\Sigma)\times\widetilde{\mathfrak{S}}_\infty]
\end{equation}

If one can check the following:
\begin{equation}
(\theta_{1,Hk}\times\theta_{2,Hk}\times id_{\widetilde{\mathfrak{S}}_\infty})^![Hk^r_{H_d}(\Sigma)\times\widetilde{\mathfrak{S}}_\infty]=[Hk^{\underline{\mu}}_{\mathcal{M}_d}]
\end{equation}
the right hand side of the equation would become:
\begin{equation}
(id, Fr)^![Hk^{\underline{\mu}}_{\mathcal{M}_d}]
\end{equation}
which can be shown is accessible to the Grothendieck-Lefschetz trace formula.

\subsection{Checking the equations}
This section is devoted to check those equations mentioned above. They basically follow from Yun-Zhang's fundamental work.

First, checking \eqref{eq1}. Let $((X-\Sigma)^d\times X^r)^\circ\subset ((X-\Sigma)^d\times X^r)$ be the locus $\{(D; x_1, ...x_r): D \text{ is disjoint from } x_1, ..., x_r\}$. Then restrict 
$Sht^{\underline{\mu}}_{H_d}(\Sigma;\Sigma_{\infty})$ to $((X-\Sigma)^d\times X^r)^\circ\subset ((X-\Sigma)^d\times X^r)$:
\begin{equation}
\xymatrix{Sht^{\underline{\mu}, \circ}_{H_d}(\Sigma;\Sigma_{\infty}) \ar[r] \ar[d] & Sht^{\underline{\mu}}_{H_d}(\Sigma;\Sigma_{\infty}) \ar[d]\\
((X-\Sigma)^d\times X^r)^\circ\times\mathfrak{S}_\infty \ar[r] & ((X-\Sigma)^d\times X^r)\times\mathfrak{S}_\infty
}
\end{equation}
Then the preimage of each  piece on the right hand of the decomposition \eqref{decom} in $Sht^{\underline{\mu}}_{H_d}(\Sigma;\Sigma_{\infty})^\circ$ is simply $Sht^{\underline{\mu}}_G(\Sigma;\Sigma_{\infty}, h_D)|_{(X-D)^r\times\mathfrak{S}_\infty}$, which is \'etale over $Sht^{\underline{\mu}}_G(\Sigma;\Sigma_{\infty})$, therefore is smooth of dimension $2r$, since $Sht^{\underline{\mu}}_G(\Sigma;\Sigma_{\infty})$ is. Then $Sht^{\underline{\mu}}_G(\Sigma;\Sigma_{\infty}, h_D)'|_{(X-D)^r\times\widetilde{\mathfrak{S}}_\infty}$, base change from $\mathfrak{S}_\infty$ to $\widetilde{\mathfrak{S}}_\infty$, is also smooth of dimension $2r$. Then the following diagram is a complete intersection diagram:
\begin{equation}
\xymatrix{Sht^{\underline{\mu},\circ}_{H_d}(\Sigma;\Sigma_{\infty})'\ar[r]\ar[d] &  Hk_{H_d}^{r, \circ}(\Sigma)\times\widetilde{\mathfrak{S}}_\infty \ar[d]^{p_{H,0}, \alpha_H}\\
H_d(\Sigma) \ar[r]^-{id, Fr} & H_d(\Sigma)\times H_d(\Sigma)
}
\end{equation}
since $Sht^{\underline{\mu},\circ}_{H_d}(\Sigma;\Sigma_{\infty})'$ is smooth of the expected dimension: 
\begin{align}
\text{dim}{Hk_{H_d}^{r, \circ}(\Sigma)\times\widetilde{\mathfrak{S}}_\infty}+\text{dim}(Im(id, Fr))-2\cdot \text{dim} H_d(\Sigma)\times H_d(\Sigma)\\
=2d+2r+3(g-1)+N-(2d+(3g-1)+N)=2r
\end{align}
Therefore $(id, Fr)^![Hk^{r,\circ}_{H_d}(\Sigma)\times\widetilde{\mathfrak{S}}_\infty]=[Sht^{\underline{\mu}, \circ}_{H_d}(\Sigma;\Sigma_{\infty})']$. But one already knows that the dimension of $Sht^{\underline{\mu}}_G(\Sigma;\Sigma_{\infty}, h_D)'$ is $2r$ (not necessarily smooth though), the equation holds not only after restricted to $((X-\Sigma)^d\times X^r)^\circ\times\widetilde{\mathfrak{S}}_\infty$, but entirely, that is to say, one has the expected equation:
\begin{equation}
(id, Fr)^![Hk^r_{H_d}(\Sigma)\times\widetilde{\mathfrak{S}}_\infty]=[Sht^{\underline{\mu}}_{H_d}(\Sigma;\Sigma_{\infty})']  
\end{equation}
$\square$

Next, check the conditions concerning the master diagram.
\begin{enumerate}
\item Condition one: In the master diagram, all stacks except the right top corner, namely $Hk^r_{H_d}(\Sigma)\times\widetilde{\mathfrak{S}}_\infty$ are smooth. Basically Yun and Zhang proved all of these in their second volumn. Here we just take the product of some of them with $\widetilde{\mathfrak{S}}_\infty$, this doesn't affect the smoothness since $\widetilde{\mathfrak{S}}_\infty$ is just a product of points.
\item  Condition two: The following fiber products are of the expected dimension:
\begin{itemize}
\item The first column: it is $Sht^{\underline{\mu}}_{T_1}(\mu_\infty\cdot\Sigma_{1,\infty})'\times_{\widetilde{\mathfrak{S}}_\infty}Sht^{\underline{\mu}}_{T_2}(\mu_\infty\cdot\Sigma_{2,\infty})'$. It is of dimension $2r$. The expected dimension is: 
\begin{align}
&\text{dim}Hk^{\underline{\mu}}_{T_1}+\text{dim}Hk^{\underline{\mu}}_{T_2}-\text{dim}Bun_{T_1}-\text{dim}Bun_{T_2}\\
&=r+\text{dim}Bun_{T_1}+r+\text{dim}Bun_{T_2}-\text{dim}Bun_{T_1}-\text{dim}Bun_{T_2}=2r
\end{align}
\item The last row: The fiber product is $\mathcal{M}_d^{\mu_\Sigma}$. For $d$ large enough, as we have seen dim$\mathcal{M}_d^{\mu_\Sigma}=2d-N-g+1$. The expected dimension is: 
\begin{align}
&\text{dim}H_d(\Sigma)-\text{dim}Bun_{T_1}-\text{dim}Bun_{T_2}\\
&=2d+3(g-1)+N-4(g-1)-2N\\
&=2d-g+1-N
\end{align}
\item The middle column: the fiber product is $Sht^{\underline{\mu}}_G(\Sigma;\Sigma_{\infty})'\times_{\widetilde{\mathfrak{S}}_\infty}Sht^{\underline{\mu}}_G(\Sigma;\Sigma_{\infty})'$, so its dimension is simply $4r$, and the expected dimension is:
\begin{align}
&2\text{dim}Hk^{\underline{\mu}}_G(\Sigma)-2\text{dim}Bun_G(\Sigma)\\
&=6(g-1)+2N+4r-6(g-1)-2N\\
&=4r
\end{align}
\item The middle row. It is simply double of the last row.
\end{itemize}
\item Condition three: 
\begin{itemize}
\item $Hk^{\underline{\mu}}_{T_1}\times Hk^{\underline{\mu}}_{T_2}\times\widetilde{\mathfrak{S}}_\infty \xrightarrow{\theta_{1,Hk}\times\theta_{2,Hk}\times id_{\widetilde{\mathfrak{S}}_\infty}} Hk^r_G(\Sigma)^2\times\widetilde{\mathfrak{S}}_\infty$ can be factored into a regular local immersion followed by a smooth relative Deligne-Mumford type morphism. Yun and Zhang showed this is true for $Hk^{\underline{\mu}}_{T_i}\longrightarrow Hk^r_G(\Sigma)$. Therefore the product also has the desired property.
\item $H_d(\Sigma)\xrightarrow{id, Fr}H_d(\Sigma)\times H_d(\Sigma)$ satisfy the conditions in (A.2.10) in their first volumn\cite{YZ1}. They proved it in their second volumn, and nothing changed here.
\end{itemize}
\item Condition four:
\begin{itemize}
\item $Sht^{\underline{\mu}}_{\mathcal{M}_d}$ admits a finite flat presentation. $Sht^{\underline{\mu}}_{\mathcal{M}_d}\longrightarrow\mathcal{M}_d^{\mu_\Sigma}$ is representable, and in our case $\mathcal{M}_d^{\mu_\Sigma}\xrightarrow{\sim}\mathcal{M}_d\times\widetilde{S}_\infty$ is actually a smooth scheme, as have seen from last section. Therefore $Sht^{\underline{\mu}}_{\mathcal{M}_d}$ admits a finite flat presentation. For the same reason $\mathcal{M}_d^{\mu_\Sigma}\xrightarrow{id, Fr}\mathcal{M}_d^{\mu_\Sigma}\times\mathcal{M}_d^{\mu_\Sigma}$ is a regular local immersion.
\item $Sht^{\underline{\mu}}_{T_i}(\mu_\infty\cdot\Sigma_{1,\infty})'$ for $i=1,2$ and $Sht^{\underline{\mu}}_G(\Sigma;\Sigma_{\infty})'$ are all smooth Deligne-Mumford stacks. So
\begin{align*}
 Sht^{\underline{\mu}}_{T_1}(\mu_\infty\cdot\Sigma_{1,\infty})'\times_{\widetilde{\mathfrak{S}}_\infty}Sht^{\underline{\mu}}_{T_2}(\mu_\infty\cdot\Sigma_{2,\infty})'\\
\xrightarrow{\theta^\mu_1\times\theta^\mu_2} Sht^{\underline{\mu}}_G(\Sigma;\Sigma_{\infty})'\times_{\widetilde{\mathfrak{S}}_\infty}Sht^{\underline{\mu}}_G(\Sigma;\Sigma_{\infty})'
\end{align*}
can automatically be factored as a regular local immersion followed by a smooth relative Deligne-Mumford type morphism, according to the observation by Yun and Zhang in their first volumn (\cite{YZ1}, Remark A4). $\square$
\end{itemize}
\end{enumerate}
 Finally checking $(\theta_{1,Hk}\times\theta_{2,Hk}\times id_{\widetilde{\mathfrak{S}}_\infty})^![Hk^r_{H_d}(\Sigma)\times\widetilde{\mathfrak{S}}_\infty]=[Hk^{\underline{\mu}}_{\mathcal{M}_d}]$. First consider the restriction of $Hk^r_{H_d}(\Sigma)$ to $(X-\Sigma)_r\times X^r$, which means $div(\phi)$ and the horizontal modification points are disjoint, and denoted by $Hk^{r, \circ}_{H_d}(\Sigma)$. As Yun-Zhang pointed out (\cite{YZ2}, Lemma 5.19), the map: $Hk^r_{H_d}(\Sigma)\longrightarrow Hk_G^r(\Sigma)\times_{Bun_G(\Sigma)}H_d(\Sigma)$ is an isomorphism when restricted to $Hk^r_{H_d}(\Sigma)^\circ$, so $Hk^{r, \circ}_{H_d}(\Sigma)$ is a smooth stack of dimension $3(g-1)+N+2r+2d$. Therefore $Hk^{\underline{\mu}, \circ}_{\mathcal{M}_d}$, which means the preimage of $Hk^{r, \circ}_{H_d}(\Sigma)\times\widetilde{\mathfrak{S}}_\infty$ in $Hk^{\underline{\mu}}_{\mathcal{M}_d}$, should have dimension:
\begin{align*}
&\text{dim}Hk^{r, \circ}_{H_d}(\Sigma)+\text{dim}Hk^{\underline{\mu}}_{T_1}+\text{dim}Hk^{\underline{\mu}}_{T_1}-2\text{dim}Hk_G^r(\Sigma)\\
&=3(g-1)+N+2r+2d+2(g-1+r)-6(g-1)-2N-4r=2d-g+1-N
\end{align*}
Therefore $(\theta_{1,Hk}\times\theta_{2,Hk}\times id_{\widetilde{\mathfrak{S}}_\infty})^![Hk^{r, \circ}_{H_d}(\Sigma)\times\widetilde{\mathfrak{S}}_\infty]=[Hk^{\underline{\mu}, \circ}_{\mathcal{M}_d}]$.

The remaining thing to check is that: $\text{dim}(Hk^{\underline{\mu}}_{\mathcal{M}_d}-Hk^{\underline{\mu}, \circ}_{\mathcal{M}_d})<\text{dim}Hk^{\underline{\mu}}_{\mathcal{M}_d}$. Since there is the isomorphism: $Hk^{\underline{\mu}}_{\mathcal{M}_d}\xrightarrow{\sim}\mathcal{H}_d^{\underline{\mu}}\times\widetilde{\mathfrak{S}}_\infty$, one only has to check a similar inequality for $\mathcal{H}_d^{\underline{\mu}}\times\widetilde{\mathfrak{S}}_\infty$. Recall the $\mathcal{H}^{\underline{\mu}}_d$ from the last section, its points can be considered as: a tuple of divisors $(E_0, E_1,... E_r)$, each of degree $2d-N$, and a sequence of points $(y_1,..., y_r)$ such that $E_i$ is obtained from $E_{i-1}$ by changing a point $y_i\in E_i$ to $\tau_3(y_i)$; together with $\Delta\in Pic^d_X(S)$ and $a\in H^0(Y_3\times S, \nu_3^*\Delta)$ such that $D_3'\in div(a)$ and $div(a)-D_3'=\pi_3(E_i):=E$ for all $i=1,...r$.
Let $\mathcal{H}^{\underline{\mu},\circ}_d$ be the subscheme of $\mathcal{H}^{\underline{\mu}}_d$ whose points have the property that $div(Tr_{Y_3/X})\cap(\sum_{i=1}^r\nu_3\pi_3(y_i))$. Then a point in $\mathcal{H}^{\underline{\mu}}_d-\mathcal{H}^{\underline{\mu},\circ}_d$ is characterized by the following:
\begin{equation*}
div(Tr_{Y_3/X}(a))\cap \nu_3(E)\neq\varnothing
\end{equation*}
This implies that there exist an $x\in |X|$ such that $\nu_3^{-1}(x)\subset E$. 
Now consider the pair $(\Delta(-x), div(a)-\nu_3^{-1}(x))$, it is a point of $\mathcal{A}_{d-1}$ by the construction. There is a map:
\begin{align*}
X\times\mathcal{A}_{d-1}&\longrightarrow\mathcal{A}_d\\
(x, (\Delta, div(a)))&\longmapsto (\Delta(x), div(a)+\nu_3(x))
\end{align*}
Then obviously $(\Delta, a)$ is the image of $(\Delta(-x), div(a)-\nu_3^{-1}(x))$. Therefore the projection of $\mathcal{H}^{\underline{\mu}}_d-\mathcal{H}^{\underline{\mu},\circ}_d$ to $\mathcal{A}_d$ lies in the image of $X\times\mathcal{A}_{d-1}$ under the above map. However, one has $\text{dim}(X\times\mathcal{A}_d)=1+2(d-1)-N-g+1=2d-N-g$. By the finiteness of the projection, one has dim$(\mathcal{H}^{\underline{\mu}}_d-\mathcal{H}^{\underline{\mu},\circ}_d)\leq\text{dim}(X\times\mathcal{A}_d)=2d-N-g<2d-N-g+1$. $\square$

 Recall that we have the following diagram:
\begin{equation}
\xymatrix{Sht^{\underline{\mu}}_{\mathcal{M}_d} \ar[r] \ar[d] & Hk^{\underline{\mu}}_{\mathcal{M}_d} \ar[d]^{p_{\mathcal{M},0}, \alpha_\mathcal{M}}\\
\mathcal{M}_d^{\mu_\Sigma} \ar[r]^-{id, Fr} & \mathcal{M}_d^{\mu_\Sigma}\times\mathcal{M}_d^{\mu_\Sigma}
}
\end{equation}
Let $\zeta=[Hk^{\underline{\mu}}_{\mathcal{M}_d}]=(\theta_{1,Hk}\times\theta_{2,Hk}\times id_{\widetilde{\mathfrak{S}}_\infty})^![Hk^r_{H_d}(\Sigma)\times\widetilde{\mathfrak{S}}_\infty]$. Then from the smoothness of $\mathcal{M}_d^{\mu_\Sigma}$, one knows that:
\begin{equation}
(id, Fr)^!\zeta\in Ch_0(Sht^{\underline{\mu}}_{\mathcal{M},d})
\end{equation}
is well defined. 

Using the decomposition \eqref{decom}, define $((id, Fr)^!\zeta)_D\in  Ch_0(Sht^{\underline{\mu}}_{\mathcal{M},D})$ be the $D$-component. Since $Sht^{\underline{\mu}}_{\mathcal{M}_d}$ is proper, one can take the degree, and define:
\begin{equation}
\langle \zeta, \Gamma(Fr_{\mathcal{M}_d^{\mu_\Sigma}}) \rangle_D:=\text{deg}((id, Fr_{\mathcal{M}_d^{\mu_\Sigma}})^!\zeta)_D
\end{equation}
 
From all have been checked above, one has:
\begin{proposition}
Suppose $D$ is an effective divisor on $X-\Sigma$ of degree $d$ such that $2d\geq 2(2g_3-1)-1+N$ (as we have seen, this is for the dimension purpose of ${\mathcal{M}_d^{\mu_\Sigma}}$), then one has: 
\begin{equation}
\langle\mathcal{Z}_1^{\mu}, h_D*\mathcal{Z}_2^{\mu}\rangle=\langle \zeta, \Gamma(Fr_{\mathcal{M}_d^{\mu_\Sigma}}) \rangle_D
\end{equation}
\end{proposition}
 
From all the facts in the last subsection, one has:
\begin{align}
(\theta'^\mu_1\times\theta'^\mu_2)^![Sht^{\underline{\mu}}_{H_d}(\Sigma;\Sigma_{\infty})']&=(\theta'^\mu_1\times\theta'^\mu_2)^!(id, Fr)^![Hk^r_{H_d}(\Sigma)\times\widetilde{\mathfrak{S}}_\infty]\\
&=(id, Fr)^!(\theta_{1,Hk}\times\theta_{2,Hk}\times id_{\widetilde{\mathfrak{S}}_\infty})^![Hk^r_{H_d}(\Sigma)\times\widetilde{\mathfrak{S}}_\infty]\\
&=(id, Fr)^![Hk^{\underline{\mu}}_{\mathcal{M}_d}]=(id, Fr)^!\zeta
\end{align}
Extracting the $D$-component and taking the degree, one gets:
\begin{equation}
\text{deg}((\theta'^\mu_1\times\theta'^\mu_2)^![Sht^{\underline{\mu}}_{H_d}(\Sigma;\Sigma_{\infty}; h_D)'])=\text{deg}((id, Fr)^!\zeta)_D
\end{equation}
By definition, left hand side is $\langle\mathcal{Z}_1^{\mu}, h_D*\mathcal{Z}_2^{\mu}\rangle$ while the right hand side is $\langle \zeta, \Gamma(Fr_{\mathcal{M}_d^{\mu_\Sigma}}) \rangle_D$.

\subsection{The intersection number and the Grothendieck-Lefschetz trace formula}
Define the ``moduli of shtukas" $\mathcal{S}^{\mu}_d$, for $\mathcal{M}_d$ as follow first:
\begin{equation}\label{shth}
\xymatrix{\mathcal{S}^{\mu}_d \ar[r] \ar[d] & \mathcal{H}^{\mu}_d \ar[d]^{(p_{\mathcal{H},0}, p_{\mathcal{H},r})}\\
\mathcal{M}_d \ar[r]^-{(id,Fr)} & \mathcal{M}_d\times\mathcal{M}_d
}
\end{equation}
Base change all the spaces above to $\widetilde{\mathfrak{S}}_\infty$:
\begin{equation}\label{shth'}
\xymatrixcolsep{5pc}\xymatrix{\mathcal{S}^{\mu}_d\times\widetilde{\mathfrak{S}}_\infty \ar[r] \ar[d] & \mathcal{H}^{\mu}_d\times\widetilde{\mathfrak{S}}_\infty \ar[d]^{(p_{\mathcal{H},0}\times id_{\widetilde{\mathfrak{S}}_\infty}, p_{\mathcal{H},r}\times id_{\widetilde{\mathfrak{S}}_\infty})}\\
\mathcal{M}_d\times\widetilde{\mathfrak{S}}_\infty \ar[r]^-{(id,Fr\times id_{\widetilde{\mathfrak{S}}_\infty})} & (\mathcal{M}_d\times\widetilde{\mathfrak{S}}_\infty)\times(\mathcal{M}_d\times\widetilde{\mathfrak{S}}_\infty)
}
\end{equation}

We want to compare this diagram with \eqref{shtm}. From meaning of its functor of points, 
\begin{equation}\label{shtm'}
\xymatrixcolsep{5pc}\xymatrix{Sht^{\underline{\mu}}_{\mathcal{M}_d} \ar[r] \ar[d] & Hk^{\underline{\mu}}_{\mathcal{M}_d} \ar[d]^{(p_{\mathcal{M},0}, p_{\mathcal{M},r})}\\
\mathcal{M}_d^{\mu_\Sigma} \ar[r]^-{(id, AL^{-1}_{\mathcal{M},\infty}\circ Fr)} & \mathcal{M}_d^{\mu_\Sigma}\times\mathcal{M}_d^{\mu_\Sigma}
}
\end{equation}

Actually \eqref{shth'} and \eqref{shtm'} are isomorphic, with the assistance of the following lemma:
\begin{lemma}
The following diagram is commutative:
\begin{equation}
\xymatrixcolsep{5pc}\xymatrix{\mathcal{M}_d\times\widetilde{\mathfrak{S}}_\infty\ar[r]^{Fr\times id_{\widetilde{\mathfrak{S}}_\infty}} \ar[d] &\mathcal{M}_d\times\widetilde{\mathfrak{S}}_\infty \ar[d]\\
\mathcal{M}_d^{\mu_\Sigma} \ar[r]^-{(id, AL^{-1}_{\mathcal{M},\infty}\circ Fr)} & \mathcal{M}_d^{\mu_\Sigma}}
\end{equation}
\end{lemma}
I'm going to check the diagram with the vertical arrows reversed. Since one can checks it locally, here I only check one $x\in\Sigma_\infty$. In other words I assume $\Sigma_\infty=\{x\}$ and $\Sigma_f=\varnothing$. Take an $S$-object of the left bottom corner:
\begin{equation*}
(\mathcal{L},\mathcal{F},x:S\longrightarrow\text{Spec}k_{v_x}, \phi:\nu_{1,*}\mathcal{L}\longrightarrow\nu_{2,*}\mathcal{F}). 
\end{equation*}
Going up along the reversed left arrow, one gets: 
\begin{equation*}
\widetilde{\mathcal{L}}\otimes\widetilde{\mathcal{F}}(-(y^{(d+1)},z^{(1)})-...-(y^{(2d)},z^{(d)})), x:S\longrightarrow\text{Spec}k_{v_x}. 
\end{equation*}
Then going along the top horizontal arrow, one gets:
\begin{equation*}
{^{\tau}\widetilde{\mathcal{L}}}\otimes{^{\tau}\widetilde{\mathcal{F}}}(-(y^{(d+2)},z^{(2)})-...-(y^{(2d)}, z^{(d)})-(y^{(1)},z^{(d+1)}))
\end{equation*}
On the other hand, if the original objects goes along the bottom horizontal arrow, it becomes:
\begin{equation*}
^\tau\mathcal{L}(y^{(1)}), {^\tau\mathcal{F}(z^{(1)})}, x:S\longrightarrow\widetilde{\mathfrak{S}}_\infty
\end{equation*}
Then going up long the right vertical arrow, one gets:
\begin{align*}
^\tau\widetilde{\mathcal{L}}\otimes{^\tau\widetilde{\mathcal{F}}}(-(y^{(1)},z^{(1)})-(y^{(1)},z^{(d+1)})+(y^{(1)},z^{(1)})+(y^{(d+1)},z^{(1)})\\
-(y^{(d+1)},z^{(1)})-(y^{(d+2)},z^{(2)})...-(y^{(2d)},z^{(d)}))
\end{align*}
One can see that $(y^{(1)}, z^{(1)})$ and $(y^{(d+1)}, z^{(1)})$ both cancel away, so the remaining is:
\begin{equation*}
{^{\tau}\widetilde{\mathcal{L}}}\otimes{^{\tau}\widetilde{\mathcal{F}}}(-(y^{(1)},z^{(d+1)})-(y^{(d+2)},z^{(2)})-...-(y^{(2d)}, z^{(d)}))
\end{equation*}
It is the same as you going the other way around as above. $\square$

Now via the identification $Hk^{\underline{\mu}}_{\mathcal{M}_d}\xrightarrow{\sim}\mathcal{H}^{\mu}_d\times\widetilde{\mathfrak{S}}_\infty$, we still use $\zeta$ for $[\mathcal{H}^{\mu}_d\times\widetilde{\mathfrak{S}}_\infty]$. Since we just saw that one can identify \eqref{shtm'} with \eqref{shth'}, it is ok to instead consider 
\begin{align}
(id, Fr\times id_{\widetilde{\mathfrak{S}}_\infty})^!\zeta&=(id, Fr\times id_{\widetilde{\mathfrak{S}}_\infty})^![\mathcal{H}^{\mu}_d\times\widetilde{\mathfrak{S}}_\infty]\\
&=((id, Fr)^![\mathcal{H}^\mu_d]\times[\widetilde{\mathfrak{S}}_\infty])\in Ch_0(\mathcal{S}^\mu_d\times\widetilde{\mathfrak{S}}_\infty)
\end{align}
Extracting the $D$ component one gets:
\begin{equation}
\langle(\zeta, \Gamma(Fr))\rangle_D=\text{deg}(\widetilde{\mathfrak{S}}_\infty)\cdot\langle [\mathcal{H}^\mu_d], \Gamma(Fr)\rangle_D
\end{equation}

As you know, the image of the map $\mathcal{S}_d^\mu\longrightarrow\mathcal{M}_d\xrightarrow{f_d}\mathcal{A}_d$ is in the rational points of $\mathcal{A}_d$, i.e. $\mathcal{A}(k)$. This is because any point in its image is ``the same" as its Frobenius twist, therefore must descent to $k$. $\mathcal{S}_d^\mu$ can be decomposed in the following way:
\begin{equation}
\mathcal{S}^\mu_d=\coprod_{a\in\mathcal{A}(k)}\mathcal{S}^\mu_d(a)
\end{equation}
Under the identification of \eqref{shtm'} with \eqref{shth'}, one has:
\begin{equation}
Sht^{\underline{\mu}}_{\mathcal{M}_D}\xrightarrow{\sim}\coprod_{a\in\mathcal{A}_D(k)}\mathcal{S}^\mu_d(a)
\end{equation}

$[\mathcal{H}_d^\mu]$ defines a cohomological self-correspondence of $(\mathcal{M}_d, \mathbb{Q}_l)$, therefore induces an endomorphism of the cohomology sheaf (complex) $\mathrm{R}f_{d,!}\mathbb{Q}_l$:
\begin{equation}
f_{d,!}[\mathcal{H}_d^\mu]: \mathrm{R}f_{d,!}\mathbb{Q}_l\longrightarrow\mathrm{R}f_{d,!}\mathbb{Q}_l
\end{equation}

According to Yun-Zhang, one can apply the Grothendieck-Lefschetz trace formula to the diagram \eqref{shth}, getting the following formula:
\begin{equation}
\langle [\mathcal{H}^\mu_d], \Gamma(Fr)\rangle_D=\sum_{a\in\mathcal{A}_D(k)}\text{Tr}(f_{d,!}[\mathcal{H}_d^\mu]\circ Fr_a, (\mathrm{R}f_{d,!}\mathbb{Q}_l)_{\overline{a}})
\end{equation}
In which $\overline{a}$ is the geometric point sitting over $a$. Since $\mathcal{H}^\mu_d$ is built up from composing the building block $\mathcal{H}_d$ with itself, one can furthur writes $f_{d,!}[\mathcal{H}_d^\mu]=(f_{d,!}[\mathcal{H}_d])^r$.

In summary, we reach the following:
\begin{proposition}
Suppose $D$ has degree greater than $2g_3-1+N$, then one has:
\begin{equation}
\langle\mathcal{Z}_1^{\mu}, h_D*\mathcal{Z}_2^{\mu}\rangle=(\prod_{x\in\Sigma_\infty}2d_x)\cdot\sum_{a\in\mathcal{A}_D(k)}\text{Tr}(f_{d,!}[\mathcal{H}_d^\mu]\circ Fr_a, (\mathrm{R}f_{d,!}\mathbb{Q}_l)_{\overline{a}})
\end{equation}
in other words:
\begin{equation}
\mathbb{I}^{\mu}(h_D)=\sum_{a\in\mathcal{A}_D(k)}\text{Tr}(f_{d,!}[\mathcal{H}_d^\mu]\circ Fr_a, (\mathrm{R}f_{d,!}\mathbb{Q}_l)_{\overline{a}})
\end{equation}
\end{proposition}

\section{The analytic side}
\subsection{The moduli space on the analytic side} Recall that we fixed a choice of $\widetilde{\mathfrak{S}}_\infty$, for each $x\in\Sigma_f$ an over point in $Y$, and $(\mu_f,\mu_\infty)\in\{\pm1\}^{\Sigma_f}\times\{\pm1\}^{\Sigma_\infty}$. This choice defined a divisor $D_3'=\sum_{x\in\Sigma_f}w'_x$ on $Y_3$. 

Now fix a nonnegative integer $d$. For any pair of nonnegative integers $(d_1,d_2)$ such that $d_1+d_2=2d$, define the following moduli space:
\begin{definition}
Let $\widetilde{\mathcal{N}}_{(d_1,d_2)}$ be the stack over $k$ whose functor of points on a $k$-scheme $S$ consisting of the following:
\begin{itemize}
\item A line bundle $\mathcal{L}$ on $Y_3\times_kS$, a splitted rank $2$ vector bundle $\mathcal{L}_1\oplus\mathcal{L}_2$ on $X\times_kS$. 
\item A map of coherent sheaves over $X\times_kS$, $\phi:\nu_{3,*}\mathcal{L}\longrightarrow\mathcal{L}_1\oplus\mathcal{L}_2$, such that deg$\phi=d$ and $\text{div}(\text{det}(\phi))$ is away from $\Sigma$. Also we require $\phi(\nu_{3,*}(\mathcal{L}(-w'_x)))\subset \mathcal{L}_1\oplus\mathcal{L}_2(-x)$ for each $x\in\Sigma$.
\item $2\text{deg}\mathcal{L}_1-\text{deg}\mathcal{L}=d_1$ and $2\text{deg}\mathcal{L}_2-\text{deg}\mathcal{L}=d_2$
\end{itemize}
\end{definition}
Then as usual, $\mathcal{N}_{(d_1,d_2)}$ is defined to be $\widetilde{\mathcal{N}}_{(d_1,d_2)}$ modulo the action of $Pic_X$. Collecting all $(d_1, d_2)$ such that $d_1+d_2=2d$, one can denote: 
\begin{equation}
\mathcal{N}_d=\coprod_{d_1+ d_2=2d}\mathcal{N}_{(d_1,d_2)}.
\end{equation}

This $\mathcal{N}_{(d_1,d_2)}$ fits into a commutative diagram similar to \eqref{Mdiagram}, in particular, admits a ``Hitchin type fiberation" over the moduli space $\mathcal{A}_d$. To get it, first base change $\phi$ to $Y_3$, one gets:
\begin{equation}
\mathcal{L}\oplus\sigma_3^*\mathcal{L}\longrightarrow\nu^*_3\mathcal{L}_1\oplus\nu^*_3\mathcal{L}_2
\end{equation}
Let $\phi_{11}:\mathcal{L}\longrightarrow\nu^*_3\mathcal{L}_1$ and $\phi_{21}:\mathcal{L}\longrightarrow\nu^*_3\mathcal{L}_2$. Then the condition $\phi(\nu_{3,*}(\mathcal{L}(-w'_x)))\subset \mathcal{L}_1\oplus\mathcal{L}_2(-x)$ for each $x\in\Sigma$ becomes:
\begin{equation}
\mathcal{L}(-w'_x)\oplus\sigma_3^*\mathcal{L}(-w_x)\longrightarrow\nu^*_3{\mathcal{L}}_1\oplus\nu^*_3{\mathcal{L}}_2(-w_x-w'_x)
\end{equation}
This puts restrictions on $\phi_{21}$:
\begin{equation}
\sigma_3(\phi_{21})|_{w'_x}=0 \text{ and } \phi_{21}|_{w_x}=0
\end{equation}
Therefore one can view $\sigma_3(\phi_{21})$ as a section of the line bundle $\sigma_3^*\mathcal{L}^{-1}\otimes\nu^*_3{\mathcal{L}}_2(-D_3')$, which is of degree $d_2-N$. By recording $\phi_{11}$ and $\sigma_3(\phi_{21})$, one gets the following morphism:
\begin{align}
\mathcal{N}_{(d_1,d_2)}&\longrightarrow\hat{Y}_{3,d_1}\times_k\hat{Y}_{3,d_2-N}\\
(\mathcal{L},\mathcal{L}_1,\mathcal{L}_2)&\longmapsto ((\mathcal{L}^{-1}\otimes\nu^*_3{\mathcal{L}}_1, \phi_{11}), (\sigma_3*\mathcal{L}^{-1}\otimes\nu^*_3{\mathcal{L}}_2(-D'_3), \sigma_3(\phi_{21}))
\end{align}
Taking the tensor product of $\phi_{11}$ and $\sigma(\phi_{21})$, one gets $\phi_{11}\otimes\sigma_3(\phi_{21})$ as a section of $\mathcal{L}^{-1}\otimes\sigma_3^*\mathcal{L}^{-1}\otimes\nu^*_3{\mathcal{L}}_1\otimes\nu^*_3{\mathcal{L}}_2$ vanishing at $D_3'$. But $\mathcal{L}^{-1}\otimes\sigma_3^*\mathcal{L}^{-1}\otimes\nu^*_3{\mathcal{L}}_1\otimes\nu^*_3{\mathcal{L}}_2$ is nothing but $\nu^*_3\Delta$, in which $\Delta$ is $\text{det}(\nu_{3,*}\mathcal{L})^{-1}\otimes\text{det}(\mathcal{L}_1\oplus\mathcal{L}_2)$. 

The fact that $\mathcal{L}^{-1}\otimes\nu^*_3{\mathcal{L}}$ admits a global section $\phi_{11}$ implies that $d_1\geq 0$ and $\sigma_3(\phi_{21})$ must vanish at $D'_3$ forces $d_2\geq N$, one has:
\begin{equation}
\mathcal{N}_d=\coprod_{d_1+d_2=2d, d_1\geq 0, d_2\geq N}\mathcal{N}_{(d_1,d_2)}
\end{equation}

From the definition, one obviously has the commutative diagram:
\begin{equation}
\xymatrix{\mathcal{N}_{(d_1,d_2)} \ar[r] \ar[d] & \hat{Y}_{3,d_1}\times_k\hat{Y}_{3,d_2-N} \ar[d]\\
\mathcal{A}_d \ar[r] & \hat{Y}_{3,2d-N}}
\end{equation}

Actually one has the following:
\begin{proposition}
This diagram is Cartesian.
\end{proposition}
It is amount to showing that given: $(\Delta, a)\in\mathcal{A}_d(S)$ and $(\mathcal{K}_1, \phi_{11})\in\hat{Y}_{3,d_1}$, $(\mathcal{K}_2, \phi_{22})\in\hat{Y}_{3,d_2-N}$ together with an isomorphism $\nu^*_{3}\Delta(-D_3')\xrightarrow{\sim}\mathcal{K}_1\otimes\mathcal{K}_2$ carrying $a'$ to $\phi_{11}\otimes \phi_{22}$ (recall that I use $a'$ to mean $a$ viewed as a section of $\nu^*_3\Delta(-D'_3)$.), one can reconstruct a unique object in $\mathcal{N}_{(d_1,d_2)}$. It is more or less the same as what Howard-Shnidman did, but paying attention to the sublattice or vanishing of the sections at $D_3'$.

Let $\mathcal{K}'_2=\mathcal{K}_2(D'_3)$. Twisting the isomorphism $\nu^*_{3}\Delta(-D_3')\xrightarrow{\sim}\mathcal{K}_1\otimes\mathcal{K}_2$ by $D_3'$, one gets: 
\begin{equation*}
\nu^*_{3}\Delta\xrightarrow{\sim}\mathcal{K}_1\otimes\mathcal{K}_2(D'_3)=\mathcal{K}_1\otimes\mathcal{K}_2'
\end{equation*}
Take $\mathcal{L}_1$ to be any line bundle over $X\times_kS$. Then you can guess that $\mathcal{L}$ is simply defined to be $\mathcal{K}_1^{-1}\otimes\widetilde{\mathcal{L}}_1$. Then from the definition of the map $\mathcal{N}_{(d_1,d_2)}\longrightarrow\hat{Y}_{3,2d-N}$, one can guess that $\nu_3^*\mathcal{L}_2$ should be defined as $\sigma^*_3\mathcal{L}\otimes\mathcal{K}'_2$. So one has to show that $\sigma_3^*\mathcal{L}\otimes\mathcal{K}'_2$ can descent to $X$ to get $\mathcal{L}_2$, but:
\begin{align*}
\sigma_3^*(\sigma^*_3\mathcal{L}\otimes\mathcal{K}_2')&=\sigma^*_3(\sigma^*_3\mathcal{K}^{-1}_1\otimes\sigma_3^*\nu_3^*\mathcal{L}_1)\otimes\sigma^*_3\mathcal{K}'_2\\
&\xrightarrow{\sim}\mathcal{K}_1^{-1}\otimes\nu^*_3\mathcal{L}_1\otimes\sigma^*_3\mathcal{K}'_2\xrightarrow{\sim}\mathcal{K}_1^{-1}\otimes\sigma^*_3\mathcal{K}_1^{-1}\otimes\nu^*_3\Delta\otimes\nu_3^*\mathcal{L}_1\\
&\xrightarrow{\sim}(\mathcal{K}_1^{-1}\otimes\nu^*_3\Delta)\otimes\sigma^*_3\mathcal{K}^{-1}_1\otimes\nu^*_3\mathcal{L}_1\xrightarrow{\sim}\mathcal{K}'_2\otimes\sigma^*_3(\mathcal{K}_1^{-1}\otimes\nu^*_3\mathcal{L}_1)\\
&\xrightarrow{\sim}\sigma^*_3\mathcal{L}\otimes\mathcal{K}'_2
\end{align*}
All the isomorphisms are either canonical or from the definition. Therefore, $\sigma^*_3\mathcal{L}\otimes\mathcal{K}'_2$ is an $\sigma_3$-equivariant line bundle, hence descent to $X$. One can view $\phi_{11}$ as a section of $Hom_{\mathcal{O}_{Y_3}}(\mathcal{L},\nu^*_3\mathcal{L}_1)$ and $\phi_{22}$ as a section of $Hom_{\mathcal{O}_{Y_3}}(\mathcal{L},\nu^*_3\mathcal{L}_2)$ vanishing at $D'_3$. Now the homomorphisms of line bundles over $Y_3$:
\begin{align*}
\phi_{11}: \mathcal{L}&\longrightarrow\nu^*_3\mathcal{L}_1\\
\phi_{12}=\sigma_3(\phi_{11}): \sigma^*_3\mathcal{L}&\longrightarrow\sigma^*_3\nu^*_3\mathcal{L}_1\xrightarrow{\sim}\nu^*_3\mathcal{L}_1\\
\phi_{21}=\sigma_3(\phi_{22}): \sigma^*_3\mathcal{L}&\longrightarrow\sigma^*_3\nu^*_3\mathcal{L}_2\xrightarrow{\sim}\nu^*_3\mathcal{L}_2\\
\phi_{22}: \mathcal{L}&\longrightarrow\nu^*_3\mathcal{L}_2\\
\end{align*}
descent to a homomorphism: $\phi: \nu_{3,*}\mathcal{L}\longrightarrow\mathcal{L}_1\oplus\mathcal{L}_2$. Since $\phi_{22}$ vanishes at $D_3'$ and $\phi_{21}$ vanishes at $\sigma_3(D_3')$, $\phi$ actually sends $\nu_{3,*}(\mathcal{L}(-w_x'))$ to $\mathcal{L}_1\oplus\mathcal{L}_2(-x)$ for each $x\in\Sigma$. $\square$

\subsection{The orbital integral}
Recall $\widetilde{G}_3=GL_{2,\nu_*\mathcal{O}_{Y_3}}$ be the group scheme over $X$ whose functor of points on an $X$-scheme $U$ is $H^0(U, End_{\mathcal{O}_U}(f^*(\nu_*\mathcal{O}_{Y_3})))$, $G_3=\widetilde{G}_3/Z(\widetilde{G}_3)$. Also, $G=PGL_{2,X}$.

The pushfoward $\nu^*_3\mathcal{O}_{Y_3}$ is not isomorphic to $\mathcal{O}_X\oplus\mathcal{O}_X$, but their generic fibers are simply $K_3$ and $F\oplus F$, both $2$-dimensional vector space over $F$. Fix an isomorphism $\rho: F\oplus F\longrightarrow K_3$. This isomorphism induces isomorphisms: $ F_x\oplus F_x\longrightarrow(K_3)_x$ for all $x\in |X|$. Almost all of them are integral, except at finitely many $x$s. Therefore it induces an isomorphism: $\rho: \mathbb{A}\oplus\mathbb{A}\longrightarrow\mathbb{A}_3$, in which $\mathbb{A}$ is the adelic ring of $X$, and $\mathbb{A}_3=K_3\otimes_F\mathbb{A}$. It carries $\mathbb{O}\oplus\mathbb{O}$ not exactly to $\mathbb{O}_3$. One can take an $h\in G_3(\mathbb{A})$ such that $\rho(\mathbb{O}\oplus\mathbb{O})=h\cdot\mathbb{O}_3$. Any $x\in\Sigma$ splits in $K_3$ into $w_x$ and $w'_x$ as we have seen above. Therefore one has $(K_3)_x=K_{w_x}\oplus K_{w'_x}$ and each one of them is isomorphic to $F_x$ as local fields. For later use, I furthur require that $h^{-1}$ carries $\rho(O_x\oplus\varpi_xO_x)$ into $O_{w_x}\oplus \varpi_{w_x'}O_{w'_x}$. Use $\rho$ and $h$, one can define the following maps:
\begin{align}
&\rho: G_{F}\longrightarrow G_{3,F}, g\longmapsto\rho\circ g\circ\rho^{-1}\\
&\rho\cdot h: G(\mathbb{A})\longrightarrow G_3(\mathbb{A}), g\longmapsto\rho\circ g\circ\rho^{-1}\cdot h\label{match adelic groups}
\end{align}
The first map is an isomorphism of group schemes over $F$, not over $X$, and the second map is an isomorphism of topological spaces, $\textbf{not}$ a group homomorphism. Via the first isomorphism of group schemes, one gets a one to one correspondence between the Borel subgroups of $G_0(\mathbb{A})$ and $G_3(\mathbb{A})$. Via the second isomorphism of topological spaces, one identify the space of automorphic forms on them:
\begin{align}
L^2(G(F)\backslash G(\mathbb{A}))&\longrightarrow L^2(G_3(F)\backslash G_3(\mathbb{A}))\\
\phi&\longmapsto \phi_3: y\longmapsto\phi(\rho^{-1}\circ(y\cdot h^{-1})\circ\rho)
\end{align}
Even though there is such a weird $h$ here, one still has the following:
\begin{proposition}\label{isom of l2}
The above map of function spaces enjoys the properties:
\begin{enumerate}
\item It identifies $L^2_{cusp}(G(F)\backslash G(\mathbb{A}))$ with $L^2_{cusp}(G_3(F)\backslash G_3(\mathbb{A}))$.
\item Let $\widetilde{Iw}_x$ be the Iwahori subgroup of $G(O_x)$ stablizing the lattice chain: $O_x\oplus O_x\supset O_x\oplus\varpi_xO_x$ and $\widetilde{Iw}_{w'_x}$ be the Iwahori subgroup of  $G_3(O_x)$ stablizing the lattice chain: $O_x\oplus O_{w_x}\supset O_{w'_x}\oplus\varpi_{w'_x}O_{w'_x}$. Let $Iw_x$ and $Iw_{w'_x}$ be their image in $G(O_x)$ and $G_3(O_x)$. Then the above maps induces:
\begin{equation}
L^2(G(F)\backslash G(\mathbb{A}))^{G(\mathbb{O}^\Sigma)\times\prod_{x\in\Sigma}Iw_x}\longrightarrow L^2(G_3(F)\backslash G_3(\mathbb{A}))^{G_3(\mathbb{O}^\Sigma)\times\prod_{x\in\Sigma}Iw_{w'_x}} 
\end{equation}
\end{enumerate}
\end{proposition}

For the first part, recall that $\phi\in L^2_{cusp}(G(F)\backslash G_0(\mathbb{A}))$ if it satisfy the following equation:
\begin{equation*}
\int_{F\backslash\mathbb{A}}\phi(\begin{pmatrix} 
1 & x \\
0 & 1 
\end{pmatrix}g)dx=0
\end{equation*}
for any Borel subgroup (the above equation is for the standard upper triangular one) and any $g$. Since $\rho: G_{F}\longrightarrow G_{3,F}$ is an isomorphism of group schemes, it induces a one to one correspondence between the Borels of $G_{F}$ and those of $G_{3,F}$. Now let $\phi_3$ be the function on $G_3(\mathbb{A})$ correspond to $\phi$ as above, consider the following integral:

\begin{equation*}
\int_{F\backslash\mathbb{A}}\phi_3(\rho\circ\begin{pmatrix} 
1 & x \\
0 & 1 
\end{pmatrix}\circ\rho^{-1}\cdot y)dx=0
\end{equation*}
Using the relation between $\phi$ and $\phi_3$:
\begin{align*}
&\phi_3(\rho\circ\begin{pmatrix} 
1 & x \\
0 & 1 
\end{pmatrix}\circ\rho^{-1}\cdot y)\\
&=\phi(\rho^{-1}\circ\rho\circ\begin{pmatrix} 
1 & x \\
0 & 1 
\end{pmatrix}\circ\rho^{-1}\cdot y\cdot h^{-1}\circ\rho)\\
&=\phi(\begin{pmatrix} 
1 & x \\
0 & 1 
\end{pmatrix}\cdot\rho^{-1}\circ y\cdot h^{-1}\circ\rho)
\end{align*}
Therefore one has the equation:
\begin{equation*}
\int_{F\backslash\mathbb{A}}\phi_3(\rho\circ\begin{pmatrix} 
1 & x \\
0 & 1 
\end{pmatrix}\circ\rho^{-1}\cdot y)dx=\int_{F\backslash\mathbb{A}}\phi(\begin{pmatrix} 
1 & x \\
0 & 1 
\end{pmatrix}\cdot\rho^{-1}\circ y\cdot h^{-1}\circ\rho)dx
\end{equation*}
Now $\rho^{-1}\circ y\cdot h^{-1}\circ\rho$ is an element of $G(\mathbb{A})$, so by the cupidality of $\phi$, the right hand side integral vanishes, therefore also the left hand side one.

For the second part, one needs to check that $\phi_3(y\cdot g)=\phi_3(y)$ for any $g\in G_3(\mathbb{O}^\Sigma)\times\prod_{x\in\Sigma}Iw_{w'_x}$, if $\phi$ is fixed by $G(\mathbb{O}^\Sigma)\times\prod_{x\in\Sigma}Iw_x$. First let's check for $g\in id^{\Sigma}\times\prod_{x\in\Sigma}Iw_{w'_x}$. To see this:
\begin{align*}
\phi_3(y\cdot g)&=\phi(\rho^{-1}\circ y\cdot g\cdot h^{-1}\circ\rho)\\
&=\phi((\rho^{-1}\circ y\cdot h^{-1}\circ\rho)\cdot(\rho^{-1}\circ h\cdot g\cdot h^{-1}\circ\rho))
\end{align*}
As I required, $h^{-1}\cdot\rho(O_x\oplus \varpi_xO_x)=O_{w_x}\oplus\varpi_{w'_x}O_{w'_x}$ for every $x\in\Sigma$. Since $g\in\prod_{x\in\Sigma}Iw_{w_x}$, one gets: $g\cdot h^{-1}\cdot\rho(O_x\oplus \varpi_xO_x)=O_{w_x}\oplus\varpi_{w'_x}O_{w'_x}$. Therefore $\rho^{-1}\circ h\cdot g\cdot h^{-1}\circ\rho$ stablizes the lattice chain: $O_x\oplus O_x\supset O_x\oplus\varpi_xO_x$, therefore is in $\prod_{x\in\Sigma}Iw_x$. Therefore one finds:
\begin{align*}
&\phi((\rho^{-1}\circ y\cdot h^{-1}\circ\rho)\cdot(\rho^{-1}\circ h\cdot g\cdot h^{-1}\circ\rho))\\
&=\phi(\rho^{-1}\circ y\cdot h^{-1}\circ\rho)=\phi_3(y)
\end{align*}
The spherical part is the same, only easier because you only have to care about the maximal lattices not lattice chains. $\square$

To relate the orbital side of relative trace formula to the moduli space on the analytic side, one needs the following $J$ space, introduced by Howard and Shnidman.
\begin{definition}
Let $\widetilde{J}$ be the following functor over $X$:
\begin{align*}
(Sch/X)&\longrightarrow (Ens)\\
U&\longmapsto Isom_{\mathcal{O}_U}(\nu_{3,*}\mathcal{O}_{Y_3},\mathcal{O}_X\oplus\mathcal{O}_X)
\end{align*}
\end{definition}
It is representable by a scheme over $X$. Let $J$ be $\widetilde{J}/\mathbb{G}_{m,X}$. It is a $G\times G_3$ bitorsor. 

One can take its adelic point:
\begin{equation*}
J(\mathbb{A})=Iso_{\mathbb{A}}(\mathbb{A}_3,\mathbb{A}\oplus\mathbb{A})/\mathbb{A}^*
\end{equation*}
Using $\rho$ and $h$, it can be identified with $G(\mathbb{A})$ as topological spaces:
\begin{align}
J(\mathbb{A})&\longrightarrow G(\mathbb{A})\\
\psi&\longmapsto \psi\circ h^{-1}\circ\rho
\end{align}

This map has the property:
\begin{proposition}\label{match of hecke}
Let $\widetilde{Iw}_{J,x}$ be the subspace of $\widetilde{J}(O_x)$ sending the lattice chain $O_x\oplus O_{w_x}\supset O_{w'_x}\oplus\varpi_{w'_x}O_{w'_x}$ to the lattice chain $O_x\oplus O_x\supset O_x\oplus\varpi_xO_x$ for each $x\in\Sigma$. Let $Iw_{J,x}$ be its image in $J(O_x)$. Then the above isomorphism of topological spaces identifies $Iw_{J,x}$ with $Iw_x$. Away from $\Sigma$, it identifies $J(O_x)$ with $G(O_x)$.
\end{proposition}

Let $C_c^\infty(G(\mathbb{A}))$ be the compactly supported group algebra (test functions) of $G(\mathbb{A})$. Via the above isomorphism of topological spaces, one can get the an isomorphism of function spaces:
\begin{align}\label{match of test functions}
C_c^\infty(G(\mathbb{A}))&\longrightarrow C_c^\infty(J(\mathbb{A}))\\
f&\longmapsto f_J: \psi\longmapsto f(\psi\circ h^{-1}\circ\rho)
\end{align}
We also identify $J(F)$ and $G(F)$ as follows:
\begin{align}
G(F)&\longrightarrow J(F)\\
\gamma&\longmapsto\gamma\circ\rho
\end{align}
Note that it is different from the map on the adelic points. 

Recall the integration kernel for $f\in C_c^\infty(G(\mathbb{A})$ in the trace formulae:
\begin{align}
\mathbb{K}_f: G(\mathbb{A})\times G(\mathbb{A})&\longrightarrow\mathbb{Q}_l\\
(g_0, g_0')&\longmapsto \sum_{\gamma\in G(F)}f(g_0^{-1}\gamma g'_0)
\end{align}
To relate it with the moduli space on the analytic side, we want to ``make" it a kernel function on $G(\mathbb{A})\times J(\mathbb{A})$. 
\begin{proposition}
Let $f$ and $f_J$ be the functions on $G(\mathbb{A})$ and $J(\mathbb{A})$ respectively, matching each other as \eqref{match of test functions}. Let $\gamma'\in G(F)$ and $\gamma=\gamma'\circ\rho^{-1}\in J(\mathbb{A})$. Let $g_3\in G_3(\mathbb{A})$ and $g_0'\in G(\mathbb{A})$ match each other as in \eqref{match adelic groups}. Then one has:
\begin{equation}
f_J(g_0^{-1}\gamma g_3)=f(g_0^{-1}\gamma' g'_0)
\end{equation}
\end{proposition}
To see it:
\begin{align*}
f_J(g_0^{-1}\gamma g_3)&=f(g^{-1}_0\gamma g_3\circ h^{-1}\circ\rho)\\
&=f(g_0^{-1}\gamma\circ\rho\circ\rho^{-1}\circ g_3\circ h^{-1}\circ\rho)
\end{align*}
By \eqref{match adelic groups}, one has $\rho^{-1}\circ g_3\circ h^{-1}\circ\rho=g'_0$. Also since $\gamma'=\gamma\circ\rho$, one sees the right hand side is exactly $f(g_0^{-1}\gamma' g_0')$. $\square$

For $f_J$, define: 
\begin{align}
\mathbb{K}_{f_J}: G(\mathbb{A})\times G_3(\mathbb{A})&\longrightarrow\mathbb{Q}_l\\
(g_0, g_3)&\longmapsto \sum_{\gamma\in J(F)}f_J(g_0^{-1}\gamma g_3)
\end{align}

Let $A\subset G$ be the diagonal torus and $T_3=Res_{Y_3/X}\mathbb{G}_m/\mathbb{G}_m$ be the nonsplit torus determined by $Y_3$ naturally sitting inside $G_3$ as in the introduction. By Howard-Shnidman, the fourfold cover $Y$ defines a character $\eta$ of $T_3(\mathbb{A})=\mathbb{A}_3^*/\mathbb{A}^*$. Let:
\begin{equation}
A(\mathbb{A})_n=\{t_0\in A(\mathbb{A}): \text{deg}(t_0)=n\}
\end{equation}
and $[A]_n=A(F)\backslash A(\mathbb{A})_n$. For any $f_J\in C_c^\infty(J(\mathbb{A})$, define the integral:
\begin{align}
\mathbb{J}_n(f_J,s)=\int_{[A]_n\times [T_3]}\mathbb{K}_{f_J}(t_0, t_3)|t_0|^{2s}\eta(t_3)dt_0dt_3\\
=q^{-2ns}\int_{[A]_n\times [T_3]}\mathbb{K}_{f_J}(t_0, t_3)\eta(t_3)dt_0dt_3
\end{align}
and one has the fact:
\begin{proposition}
The integral $\mathbb{J}_n(f_J,s)$ vanishes for $|n|$ sufficiently large.
\end{proposition}
Therefore it is legitimate to define:
\begin{equation}
\mathbb{J}(f_J,s)=\sum_{n\in\mathbb{Z}}\mathbb{J}_n(f_J,s)
\end{equation}

Also, let:
\begin{equation}
\mathbb{K}_{f_J,\gamma}(g_0,g_3)=\sum_{\delta\in A(F)\gamma T_3(F)}f_J(g^{-1}_0\delta g_3)
\end{equation}
For any $\gamma\in J(F)$ and $g_0\in G(F)\backslash G(\mathbb{A}), g_3\in G_3(F)\backslash G_3(\mathbb{A})$. One also define:
\begin{equation}
\mathbb{J}_n(\gamma, f_J, s)=\int_{[A]_n\times [T_3]}\mathbb{K}_{f_J, \gamma}(t_0, t_3)|t_0|^{2s}\eta(t_3)dt_0dt_3
\end{equation}
and
\begin{equation}
\mathbb{J}(f_J, s)=\sum_{n\in\Sigma}\mathbb{J}(\gamma, f_J, s)
\end{equation}
so that:
\begin{equation}
\mathbb{J}_n(f_J,s)=\sum_{\gamma\in A(F)\backslash J(F)/T_3(F)}\mathbb{J}_n(\gamma, f_J, s)
\end{equation}
and
\begin{equation}
\mathbb{J}(f_J,s)=\sum_{\gamma\in A(F)\backslash J(F)/T_3(F)}\mathbb{J}(\gamma, f_J, s)
\end{equation}

Recall that Howard and Shnidman defined the invariants of $A(F)\backslash J(F)/T_3(F)$:
\begin{proposition}
There is a bijection:
\begin{equation}
\text{inv}: A(F)\backslash J(F)/T_3(F)\longrightarrow\{a\in K_3: \text{Tr}_{K_3/F}(a)=1\}
\end{equation}
\end{proposition}
Let me review their proof. Actually this is the generic fiber version of the Hitchin type fiebration introduced above. Take a morphism of $F$-vector spaces:
\begin{equation*}
\phi: K_3\longrightarrow F\oplus F
\end{equation*}
base change it to Spec$K_3$ (the double cover) to split $K_3$ (the rank two vector bundle):
\begin{equation*}
K_3\oplus K_3\xrightarrow{\sim}K_3\otimes K_3\xrightarrow{\phi\otimes id}K_0\otimes_FK_3\xrightarrow{\sim}K_3\oplus K_3
\end{equation*}
As in its geometrization, one gets the maps of $K_3$-vector spaces:
\begin{align*}
&\phi_{11}: K_3\longrightarrow K_3, \sigma_3(\phi_{11}): K_3\longrightarrow K_3\\
&\phi_{21}: K_3\longrightarrow K_3, \sigma_3(\phi_{21}): K_3\longrightarrow K_3
\end{align*}
In which $\sigma_3(a)$ and $\sigma_3(c)$ are $\sigma_3$-linear. Let:
\begin{equation*}
\Delta=Hom_F(det_F(K_3), det_F(F\oplus F))
\end{equation*}
Then $det_F(\phi)$ is an element of this vector space (a section of this line bundle over the generic point of $X$).

Now the invariant map is defined to be $\text{inv}(\phi)=(\phi_{11}\cdot\sigma_3(\phi_{21}))/\text{det}(\phi)$. The map $\phi$ is called regular if $a=\text{inv}(\phi)$ is not $0$ in $K_3$. $\square$

Therefore one can also index the double cosets in $A(F)\backslash J(F)/T_3(F)$ by their invariants, i.e.:
\begin{equation}
\mathbb{J}(f_J,s)=\sum_{\gamma\in A(F)\backslash J(F)/T_3(F)}\mathbb{J}(\gamma, f_J, s)=\sum_{a\in K_3, \text{Tr}_{K_3/F}(a)=1}\mathbb{J}(a, f_J, s)
\end{equation}

\subsection{The orbital integral and the period integral} The relative trace formula relates the orbital integral with certain period integrals summed over automorphic forms. 
Let $\pi$ be a cuspidal automorphic representation of $G(\mathbb{A})$. For any $\phi\in\pi$, define the period integral along $A(\mathbb{A})$:
\begin{equation}
\mathscr{P}_0(\phi,s)=\int_{[A]}\phi(t_0)|t_0|^{2s}dt_0
\end{equation}
Using \eqref{isom of l2}, one can move acusp form $\phi$ on $G(\mathbb{A})$ to a cusp form $\phi_3$ on $G_3(\mathbb{A})$. Define the period of $\phi_3$ along $T_3(\mathbb{A})$:
\begin{equation}
\mathscr{P}_{3,\eta}(\phi_3)=\int_{[T_3]}\phi_3(t_3)\eta(t_3)dt_3
\end{equation}
Both integrals are absolutely convergent. For $T_3$ there is no variable since it is anisotropic.

Define the global spherical character relative to $(A\times T_3, 1\times\eta)$ as:
\begin{equation}
\mathbb{J}_\pi(f_J, s)=\sum_\phi\displaystyle\frac{\mathscr{P}_0(\pi(f)\phi,s)\mathscr{P}_{3,\eta}(\overline{\phi_3})}{\langle \phi, \phi \rangle}
\end{equation}
in which $f\in C_c^\infty(G(\mathbb{A}))$ is the function correspond to $f_J$ and $\phi_3$ is the automorphic form of $G_3(\mathbb{A})$ corresponding to $\phi$. The sum is over an orthogonal basis of $\phi$, and $\langle \phi, \phi \rangle$ is the Petersson inner product of $\phi$ with itself.

Now let's choose the test function to be used. In their second volumn, Yun and Zhang defined the Eisenstein ideal $\mathcal{I}^\Sigma_{Eis}\in\mathscr{H}_G^\Sigma$. And they proved the following theorem:
\begin{theorem}
Let $f\in\mathcal{I}^\Sigma_{Eis}$. Take any $f_\Sigma\in C_c^\infty(G(\mathbb{A_\Sigma}))$ such that it is left invariant under the Iwahori subgroup $\prod_{x\in\Sigma}Iw_x$. Then the kernel function of $f^\Sigma=f\otimes f_\Sigma$ has vanishing Eisenstein part, i.e.:
\begin{equation*}
\mathbb{K}_{f^\Sigma}=\mathbb{K}_{f^\Sigma,cusp}+\mathbb{K}_{f^\Sigma,sp}
\end{equation*}
In other words, $\mathbb{K}_f$ only has cuspidal and residual parts.
\end{theorem}

Please notice that $f^\Sigma$ $\bf{doesn't}$ mean the product of the components away from $\Sigma$ as usual, it is a product over all primes. Consider the test function $f^{\Sigma}=f\otimes(\bigotimes_{x\in\Sigma}1_{Iw_x})$ for $f\in\mathcal{I}_{Eis}^\Sigma$. Then $f_J$ has the form $f_J^\Sigma=f_J\otimes(\bigotimes_{x\in\Sigma}1_{Iw_{J,x}})$. Define the 
\begin{proposition}
For each $f^\Sigma$ of the form above, one has:
\begin{equation}
\mathbb{J}(f^\Sigma_J, s)=\sum_\pi\quad\mathbb{J}_\pi(f^\Sigma_J, s) 
\end{equation}
where the sum is over all $\pi$ with $G(\mathbb{A}^\Sigma)\times\prod_{x\in\Sigma}Iw_{x}$ fixed vector.
\end{proposition}

First on $G(\mathbb{A})$, one has:
\begin{align}
\mathbb{K}_{f^\Sigma}(g_0,g'_0)=&\sum_\pi\quad\lambda_\pi(f)\cdot\displaystyle\frac{\phi(g_0)\overline{\phi(g'_0)}}{\langle \phi, \phi \rangle}\\
&+\sum_\chi\quad\lambda_\chi(f)\cdot\cdot\chi(\text{det}(g_0))\cdot\chi(\text{det}(g'_0))
\end{align}
where the sums are over $\pi$s and $\chi$s with Iwahori fixed vectors, and $\phi$ is the $G(\mathbb{A}^\Sigma)\times\prod_{x\in\Sigma}Iw_{x}$-fixed vector. Now transfer to $J(\mathbb{A})$:
\begin{align*}
\mathbb{K}_{f_J^\Sigma}(g_0,g_3)&=\sum_\pi\quad\lambda_\pi(f)\cdot\displaystyle\frac{\phi(g_0)\overline{\phi_3(g_3)}}{\langle \phi, \phi \rangle}\\
&+\sum_\chi\quad\lambda_\chi(f)\cdot\cdot\chi(\text{det}(g_0))\cdot\chi(\text{det}(\rho\circ g_3\cdot h^{-1}\circ\rho))\\
&=\sum_\pi\quad\lambda_\pi(f)\cdot\displaystyle\frac{\phi(g_0)\overline{\phi_3(g_3)}}{\langle \phi, \phi \rangle}\\
&+\sum_\chi\quad\lambda_\chi(f)\cdot\cdot\chi(\text{det}(g_0))\cdot\chi(\text{det}(\rho^{-1}\circ g_3\circ\rho\circ\rho^{-1}\cdot h^{-1}\circ\rho))
\end{align*}
Let $\text{det}_3: G_{3,F}\longrightarrow\mathbb{G}_{m,F}$ be the character defined by $\text{det}_3(g_3):=\text{det}(\rho^{-1}\circ g_3\circ\rho)$. It is a generator of $X^*(G_{3,F})$. The last term above becomes:
\begin{equation}
\chi(\text{det}(\rho^{-1}\circ g_3\circ\rho\circ\rho^{-1}\cdot h^{-1}\circ\rho))=\chi(\text{det}_3(g_3))\cdot\chi(\text{det}_3(h^{-1}))
\end{equation}
in which $\chi(\text{det}_3(h^{-1}))$ is a nonzero number.

Now one only has to show that either one of the following is $0$:
\begin{equation*}
\int_{[A]_n}\chi(\text{det}(t_0))|t_0|^{2s}dt_0, \quad \int_{[T_3]}\chi(\text{det}_3(t_3))\eta(t_3)dt_3
\end{equation*}
If $\chi$ is non trivial on $\mathbb{A}^1=\mathbb{A}_0^*$, take a $t_0'$ such that $\chi(\text{det}(t'_0))\neq 1$. Then one has:
\begin{align*}
\int_{[A]_n}\chi(\text{det}(t_0))|t_0|^{2s}dt_0&=\int_{[A]_n}\chi(\text{det}(t'_0t_0))|t'_0t_0|^{2s}d(t'_0t_0)\\
&=\int_{[A]_n}\chi(\text{det}(t'_0))\chi(\text{det}(t_0))|t_0|^{2s}dt_0
\end{align*}
From this one gets:
\begin{equation*}
(1-\chi(\text{det}(t'_0)))\cdot\int_{[A]_n}\chi(\text{det}(t_0))|t_0|^{2s}dt_0=0
\end{equation*}
This implies the second factor is $0$.

If $\chi$ is trivial on $\mathbb{A}^1$, then the integral on $[T_3]$ becomes:
\begin{equation}
(\text{a non zero constant})\cdot\int_{[T_3]}\eta(t_3)dt_3
\end{equation}
The constant is something in $q$. Now since $\eta$ is the quadratic character defining $K$ over $K_3$, it is non trivial on $\mathbb{A}_3^1$. Applying the same trick above to $\eta$ and $\mathbb{A}_3^1$ proves the vanishing of this integral. $\square$

We want to study the $r$-th derivative of $q^{Ns}\mathbb{J}(f_J^\Sigma,s)$, hence define.:
\begin{equation}
\mathbb{J}_r(f_J^\Sigma, s)=\displaystyle\frac{d^r}{ds^r}|_{s=0}(q^{Ns}\mathbb{J}(f_J^\Sigma,s))
\end{equation}
Similarly define:
\begin{equation}
\mathbb{J}_r(a, f_J^\Sigma, s)=\displaystyle\frac{d^r}{ds^r}|_{s=0}(q^{Ns}\mathbb{J}(a, f_J^\Sigma, s))
\end{equation}
The $q^{Ns}$ is to make it match the geometric side later.

Then one has:
\begin{equation}
\mathbb{J}_r(f_J^\Sigma, s)=\sum_{a\in K_3, \text{Tr}_{K_3/F}(a)=1}\mathbb{J}_r(a, f_J^\Sigma, s)
\end{equation}
Actually one needs $f_J$'s spherical part to be in the Eisenstein ideal, but that is not directly related to the moduli spaces. Rather, one first consider those $h_D\in\mathscr{H}^\Sigma$, then use the linear combination of these elements indexed by the effective divisors on $X-\Sigma$ to approximate the elements in the Eisenstein ideal (in some representation of the Hecke algebra). With the assit of \eqref{match of hecke}, we have:
\begin{proposition}
Let $h_D\in\mathscr{H}^\Sigma$ be the element indexed by $D\in Div^+(X-\Sigma)$, then $h_D\otimes(\bigotimes_{x\in\Sigma}1_{Iw_x})$ matches $h_{D,J}\otimes(\bigotimes_{x\in\Sigma}1_{Iw_{J,x}})$ where $h_{D,J}$ is the element in $C^\infty_c(G(\mathbb{O}^\Sigma)\backslash G(\mathbb{A}^\Sigma)/G_3(\mathbb{O}^\Sigma))$ indexed by $D$.
\end{proposition}

From now on, I use $h^\Sigma_D$ to denote $h_{D,J}\otimes(\bigotimes_{x\in\Sigma}1_{Iw_{J,x}})$.

\subsection{The orbit integral and the moduli space}
The general spirit of Yun-Zhang's theory is to ``upgrade" the equality of two numbers to an ``equality" between two complexes of sheaves on the Hitchin base $\mathcal{A}_d$ first. After proving the equality of complexes of sheaves, one decategorify it via Grothendieck's faisceaux-fonctions correspondence to get the expected equality between two numbers. Now we have to express the numbers $\mathbb{J}_r(a, h^\Sigma_D, s)$ as the Frobenius traces of some complexes of sheaves on $\mathcal{A}_d$.

What that complex of sheaf is? First $Y$ is a double cover of $Y_3$, the pushforward of the constant sheaf $\pi_{3,*}\mathbb{Q}_l$ decomposes as $\mathbb{Q}_l\oplus L_{Y/Y_3}$, in which $L_{Y/Y_3}$ is the local system on $Y_3$ defined by the representation:
\begin{equation}
\pi_1^{et}(Y_3)\longrightarrow\langle \tau_3 \rangle\longrightarrow\{\pm1\}
\end{equation}
in which the second arrow is non trivial.

From the geometric class field theory for $Y_3$, there is a corresponding local system on $Pic_{Y_3}$, therefore a local system on $Pic^d_{Y_3}$ for each $d$. Pull it back along the cover:
\begin{equation}
\hat{Y}_{3,d}\longrightarrow Pic_{Y_3}
\end{equation}
one gets a local system on $\hat{Y}_{3,d}$ for every $d$. Denote it by $L_d$.

Recall that we have the morphism:
\begin{equation}
j_{(d_1,d_2)}: \mathcal{N}_{(d_1,d_2)}\longrightarrow\hat{Y}_{3,d_1}\times_k\hat{Y}_{d_2-N}
\end{equation}
Let $L_{(d_1,d_2)}=j^*_{(d_1,d_2)}(L_{d_1}\boxtimes\mathbb{Q}_l)$

Recall that we defined $\mathcal{A}_D$ before \eqref{def of AD}:
\begin{definition}
\[\xymatrixcolsep{5pc}
\xymatrix{\mathcal{A}_D \ar[r] \ar[d] & \mathcal{A}_d \ar[d]^{Tr}\\
Speck \ar[r]^{(\mathcal{O}_X(D),1)} & \hat{X}_d
}
\]
\end{definition}

Then we have the proposition parallel to the one in Howard-Shnidman:
\begin{proposition}
There is a canonical bijection:
\begin{equation*}
\text{inv}_D:\mathcal{A}_D(k)\xrightarrow{\sim}
\left\{\begin{array}{l|l}
& Tr_{K_3/F}(a)=1\\ 
a\in K_3 &\\
& div(a)+\nu_3^*D-D'_3\geq 0
\end{array}\right\}
\end{equation*}
\end{proposition}

From the definition of $\mathcal{A}_D$, a point of it consists of a line bundle $\Delta\in Pic^d(X)$, a global section $a$ of $\nu_3^*\Delta$, vanishing at $D'_3$ and an isomorphism: $(\Delta, Tr_{Y_3/X}(a))\xrightarrow{\sim}(\mathcal{O}_X(D), 1)$. The identification here $\Delta\xrightarrow{\sim}\mathcal{O}_X(D)$ induces an identification $\nu_3^*\Delta\xrightarrow{\sim}\mathcal{O}_{Y_3}(\nu_3^*D)$. Via this latter isomorphism, $a$ is identified with a global section of $\mathcal{O}_{Y_3}(\nu_3^*D)$ vanishing at $D'_3$. This is nothing but a $a\in K_3$ such that $\text{div}(a)+\nu_3^*D-D'_3\geq 0$. $\square$

One call the set on the right hand side the invariants for $\mathcal{A}_D(k)$.

Now it is the time to state what we want:
\begin{proposition}\label{analytic trace}
Let $a\in K_3$ such that $\text{Tr}_{K_3/F}(a)=1$. Identify $\mathcal{A}_d(k)$ as a subset of $K_3$ via $\text{inv}_D$ as above. Then one has the following:
\begin{enumerate}
\item If $a\notin\mathcal{A}_D(k)$, then $\mathbb{J}(a,h_D^\Sigma, s)=0$; 
\item If $a\in\mathcal{A}_D(k)$, then one has:
\begin{equation}
\mathbb{J}(a, h_D^\Sigma, s)=\sum_{(d_1+d_2=2d, d_1\geq 0, d_2\geq N)}q^{(d_1-d_2)\cdot s}\text{Tr}(Fr_a, (\mathrm{R}g_!L_{(d_1,d_2)})_{\overline{a}})
\end{equation}
where $\overline{a}$ is the geometric point sitting over the $k$-point $a:\text{Spec}k\longrightarrow\mathcal{A}_D$.
\end{enumerate}
\end{proposition}

First let $\widetilde{h}_D\in\mathscr{H}(GL_2(\mathbb{O}^\Sigma)\backslash GL_2(\mathbb{A}^\Sigma)/GL_2(\mathbb{O}^\Sigma))$ be the characteristic function of those double cosets whose determinants have divisor $D$ on $X-\Sigma$. Let:
\begin{equation}
\widetilde{h}_D^{\Sigma}:=\widetilde{h}_D\otimes\bigotimes_{x\in\Sigma}1_{\widetilde{Iw}_{J,x}}
\end{equation}
Let $\gamma$ be an element with invariant $\xi$, and $\widetilde{\gamma}$ be a preimage of $\gamma$ in GL$_2(F)$, then one has:
\begin{equation}
\mathbb{J}(\gamma, h_D^\Sigma, s)=\int_{Z(\mathbb{A})\backslash \widetilde{A}(\mathbb{A})\times\widetilde{T}_3(\mathbb{A})}\widetilde{h}_D^\Sigma(t_0^{-1}\widetilde{\gamma} t_3)|t_0|^{2s}\eta(t_3)dt_0dt_3
\end{equation}
One rewrite it as:
\begin{align*}
\int_{Z(\mathbb{A})\backslash \widetilde{A}(\mathbb{A})\times\widetilde{T}_3(\mathbb{A})/\widetilde{A}(\mathbb{O})\times\widetilde{T}_3(\mathbb{O})}(\int_{\widetilde{A}(\mathbb{O})\times\widetilde{T}_3(\mathbb{O})}\widetilde{h}_D^\Sigma(t_0^{-1}x^{-1}\gamma y t_3)|xt_0|^{2s}\eta(yt_3)dxdy)dt_0dt_3
\end{align*}
Since $|x|=1$ and $\eta$ is trivial on $\mathbb{O}_3$, one can simply it:
\begin{align*}
\int_{Z(\mathbb{A})\backslash \widetilde{A}(\mathbb{A})\times\widetilde{T}_3(\mathbb{A})/\widetilde{A}(\mathbb{O})\times\widetilde{T}_3(\mathbb{O})}\widetilde{h}_D^\Sigma(t_0^{-1}\gamma t_3)|t_0|^{2s}\eta(t_3)\cdot(\int_{\widetilde{A}(\mathbb{O})\times\widetilde{T}_3(\mathbb{O})}dxdy)dt_0dt_3\\
=\int_{Z(\mathbb{A})\backslash \widetilde{A}(\mathbb{A})\times\widetilde{T}_3(\mathbb{A})/\widetilde{A}(\mathbb{O})\times\widetilde{T}_3(\mathbb{O})}\widetilde{h}_D^\Sigma(t_0^{-1}\gamma t_3)|t_0|^{2s}\eta(t_3)dt_0dt_3
\end{align*}

To furthur simply this integral and relate it to geometry, we introduce the set $\widetilde{\mathfrak{X}}_{D,\widetilde{\gamma}}$ whose points consist of the tuples $(E, E_1, E_2)\in Div(Y_3)\times Div(X)\times Div(X)$ who satisfy the following:
\begin{itemize}
\item The map induced by $\widetilde{\gamma}: \nu_{3,*}\mathcal{O}_{Y_3}\dashrightarrow\mathcal{O}_X\oplus\mathcal{O}_X$:
\begin{equation*}
\phi_{\widetilde{\gamma}}: \nu_{3,*}(\mathcal{O}_{Y_3}(-E))\longrightarrow\mathcal{O}_X(-E_1)\oplus\mathcal{O}_X(-E_2)
\end{equation*}
is everywhere defined.
\item $\text{div}(\text{det}(\phi))$ is $D$.\\
\item it sends $\nu_{3,*}(\mathcal{O}_{Y_3}(-E)(-w'_x))$ into $\mathcal{O}_X(-E_1)\oplus\mathcal{O}_X(-E_2)(-x)$ for every $x\in\Sigma$.
\end{itemize}

Let $\mathfrak{X}_{D,\widetilde{\gamma}}=\widetilde{\mathfrak{X}}_{D,\widetilde{\gamma}}/Div(X)$. The integral furthur simplies to:
\begin{align}\label{counting point}
\int_{Z(\mathbb{A})\backslash \widetilde{A}(\mathbb{A})\times\widetilde{T}_3(\mathbb{A})/\widetilde{A}(\mathbb{O})\times\widetilde{T}_3(\mathbb{O})}\widetilde{h}_D^\Sigma(t_0^{-1}\gamma t_3)|t_0|^{2s}\eta(t_3)dt_0dt_3\\
=\sum_{(E,E_1,E_2)\in\mathfrak{X}_{D,\widetilde{\gamma}}}q^{-deg(E_1-E_2)\cdot 2s}\eta(E)
\end{align}

We have the invariant $\chi$ attached to $\widetilde{\gamma}$, let's denote $\mathfrak{X}_{D,\widetilde{\gamma}(\xi)}$.

The set $\mathfrak{X}_{D,\widetilde{\gamma}(a)}$ is closedly related to the moduli space $\mathcal{N}_{d}$ for $d=\text{deg}D$. For a $k$-rational point of $\mathcal{A}_D$, $a$, furthur define $\mathcal{N}_{d,a}$ to be the fiber product:
\begin{equation}\label{xi fiber}
\xymatrix{\mathcal{N}_{d,a} \ar[r] \ar[d] & \mathcal{N}_d \ar[d]^g\\
\text{Spec}k\ar[r]^{a} &\mathcal{A}_d
}
\end{equation}

\begin{proposition} Let $a$ be the invariant of $\widetilde{\gamma}$. If $a$ is not in the image of the invariants of $\mathcal{A}_D(k)$, then $\mathfrak{X}_{D,\widetilde{\gamma}}=\varnothing$; If $a$ is an invariant of $\mathcal{A}_D(k)$, then there exists a bijection:
\begin{equation}
\lambda: \mathfrak{X}_{D,\widetilde{\gamma}(a)}\xrightarrow{\sim}\mathcal{N}_{d,a}(k)
\end{equation}
\end{proposition}

Let's first define a map:
\begin{equation}
\lambda:  \mathfrak{X}_{D,\widetilde{\gamma}(a)}\xrightarrow{\sim}\mathcal{N}_d(k)
\end{equation}
by sending $(E,E_1,E_2)$ to $(\mathcal{O}_{Y_3}(-E), \mathcal{O}_X(-E_1), \mathcal{O}_X(-E_2), \phi_{\widetilde{\gamma}})$. Since the induced map $\phi_{\widetilde{\gamma}}$ has devisor $D$ degree $d$, the image really lies in $\mathcal{N}_d(k)$, and its image in $\mathcal{A}_d(k)$ lies in $\mathcal{A}_D(k)$. From the definition of the invariant of $\widetilde{\gamma}$ and of $\phi_{\widetilde{\gamma}}$, one can see they are the same. So $a$ is in the image of $\mathcal{A}_D(k)$. Therefore if $a$ is not an invariant of $\mathcal{A}_D(k)$, $\mathfrak{X}_{D,\widetilde{\gamma}(a)}$ must be empty.

Now one has to show that any element of $\mathcal{N}_{d,a}$ arise as a image of some point of $\mathfrak{X}_{D,\widetilde{\gamma}(a)}$. This basically follows from Howard-Shnidman. Given a tuple $(\mathcal{L}, \mathcal{L}_1, \mathcal{L}_2, \phi)\in\mathcal{N}_{d,a}(k)$, fix isomorphisms of the generic fiber of $\nu_3^*\mathcal{L}$ with $K_3$, $\mathcal{L}_1, \mathcal{L}_2$ with $F$ respectively:
\begin{equation}
g_3: \nu_{3,*}\mathcal{L}|_{\eta_X}\xrightarrow{\sim} K_3,\quad g_1:\mathcal{L}_1|_{\eta_X}\xrightarrow{\sim}F,\quad g_2:\mathcal{L}_2|_{\eta_X}:\xrightarrow{\sim}F
\end{equation}
Then $\phi$ transfer to a map of $F$-vector spaces $\gamma: K_3\longrightarrow F\oplus F$. The inverse images of $1$ under $g_3$ and $ g_1,g_2$ defines rational sections of $\nu_{3,*}\mathcal{L}$ and $\mathcal{L}_1, \mathcal{L}_2$, say $s_1, s_2, s_3$. Then $(E, E_1, E_2)$ can be selected as $(-\text{div}(s_3), -\text{div}(s_1), -\text{div}(s_2))$. The requirement at $\Sigma$ is automatic from the same requirement for $(\mathcal{L}, \mathcal{L}_1, \mathcal{L}_2, \phi)$. $\square$

One can furthur decompose $\mathcal{N}_{d,a}$ into $\mathcal{N}_{(d_1,d_2),a}$, getting:
\begin{equation}
\mathcal{N}_{d,a}=\coprod_{d_1+d_2=2d, d_1\geq 0, d_2\geq N}\mathcal{N}_{(d_1,d_2,a)}
\end{equation}
 The point counding formula \eqref{counting point} on $\mathfrak{X}_{\widetilde{\gamma}, D}(k)$ can be reinterpreted as counting points on $\mathcal{N}_{d,a}(k)$:
\begin{align}
&\sum_{(\mathcal{L}, \mathcal{L}_1, \mathcal{L}_2, \phi)\in\mathcal{N}_{d,a}(k)}q^{-deg(\mathcal{L}_1-\mathcal{L}_2)\cdot 2s}\cdot\eta(\mathcal{L})\\
&=\sum_{(d_1+d_2=2d, d_1\geq 0, d_2\geq N)}q^{(d_1-d_2)\cdot s}\sum_{(\mathcal{L}, \mathcal{L}_1, \mathcal{L}_2, \phi)\in\mathcal{N}_{(d_1,d_2),a}(k)}\eta(\mathcal{L})
\end{align}
Apply the Grothendieck-Lefschetz trace formula to the diagram \eqref{xi fiber} and the complex of vector spaces $(\mathrm{R}g_!L_{(d_1,d_2)})_a$, one get the desired equation:
\begin{equation*}
\mathbb{J}(a, h_D^\Sigma, s)=\sum_{(d_1+d_2=2d, d_1\geq 0, d_2\geq N)}q^{(d_1-d_2)\cdot s}\text{Tr}(Fr_a, (\mathrm{R}g_!L_{(d_1,d_2)})_{\overline{a}})
\end{equation*}
This finish the proof of \eqref{analytic trace}.

\section{Conclusions}
\subsection{Comparision of sheaves}
Now both the intersection number and the orbital integral are translated into Frobenius traces on some complexes of sheaves, one has to compare these sheaves.

First let me review some combinatorial facts proved by Howard-Shnidman. Let $d_2'=d_2-N$ and $d'=2d-N$. Consider the group:
\begin{equation}
\Gamma_{d'}=\{\pm1\}^{d'}\rtimes S_{d'}
\end{equation}
Let $\textbf{1}$ be the trivial representation of the group $S_{d'}$, and:
\begin{equation}
\text{Ind}^{\Gamma_{d'}}_{S_{d'}}\textbf{1}=\{\Phi: S_{d'}\backslash\Gamma_{d'}\longrightarrow\mathbb{Q}_l\}
\end{equation}
$S_{d'}$ is embedded into $\Gamma_{d'}$ as the subgroup $\{1\}^{d'}\rtimes S_{d'}$. For each $x\in\{\pm 1\}^{d'}$, let $\Phi_x$ be characteristic function of the coset $x\cdot S_{d'}$. When $x$ runs through all of these tuples formed by $1$ and $-1$, one get a basis of the induced representation above. Let:
\begin{equation*}
e_i=(1,...,1,-1,1,...,1)\in\{\pm 1\}^{d'}
\end{equation*}
in which $-1$ appears at the $i$-th spot. Let $H$ be the operator on $\text{Ind}^{\Gamma_{d'}}_{S_{d'}}\textbf{1}$ defined by:
\begin{equation*}
H\cdot\Phi_x=\sum_{i=1}^{d'}\Phi_{e_ix}
\end{equation*}
i.e. $H$ is the element $\sum_{i=1}^{d'}e_i\in\mathbb{Z}[\Gamma_{d'}]$ acting on the induced representation. 

Similarly, for a pair of numbers $d_1, d_2'$ such that $d_1+d'_2=d'$, define:
\begin{equation*}
\Gamma_{d_1}=\{\pm 1\}^{d_1}\rtimes S_{d_1}\quad \Gamma_{d_2'}=\{\pm 1\}^{d_2'}\rtimes S_{d_2'}
\end{equation*}

and characters:
\begin{equation*}
\eta_{d_1}:\{\pm 1\}^{d_1}\longrightarrow\{\pm1\}\quad\eta_{d_2}:\{\pm 1\}^{d'_2}\longrightarrow\{\pm 1\}
\end{equation*}

Then one has the following proposition of $\Gamma_{d'}$-module:
\begin{proposition}\label{repdecomp}
There is a decomposition:
\begin{equation}
\text{Ind}^{\Gamma_{d'}}_{S_{d'}}\textbf{1}=\bigoplus_{d_1,d'_2\geq 0, d_1+d'_2=d'}V_{(d_1,d'_2)}
\end{equation}
\end{proposition}
such that:
\begin{itemize}
\item Each $V_{(d_1,d'_2)}$ is an irreducible representation of $\Gamma_{d'}$;
\item $H$ acts on $V_{(d_1,d'_2)}$ by the scalar $(d_1-d'_2)$;
\item $V_{(d_1,d'_2)}$ is isomorphic $V_{(d'_2,d_1)}$;
\item $V_{(d_1,d'_2)}$ is isomorphic to $\text{Ind}^{\Gamma_{d'}}_{\Gamma_{d_1}\times\Gamma_{d'_2}}(\eta_{d_1}\boxtimes\textbf{1})$ and $\text{Ind}^{\Gamma_{d'}}_{\Gamma_{d_1}\times\Gamma_{d'_2}}(\textbf{1}\boxtimes\eta_{d_1})$.
\end{itemize}

Now we can state the basic fact about sheaves comparison. Recall one has the following diagram:
\begin{equation}
\xymatrix{ & \mathcal{H}_d \ar[ld]_{\gamma_0} \ar[rd]^{\gamma_1} & \\
\mathcal{M}_d \ar[rd]_{f_d} & & \mathcal{M}_d \ar[ld]^{f_d}\\
&\mathcal{A}_d &}
\end{equation}
\begin{proposition}
Assume $d\geq 2g_3-1+N$. One has the following isomorphism:
\begin{equation}
\mathrm{R}f_{d,!}\mathbb{Q}_l\xrightarrow{\sim}\bigoplus_{d_1\geq 0, d_2\geq N, d_1+d_2=2d}\mathrm{R}g_{(d_1,d_2),!}L_{(d_1, d_2)}
\end{equation}
\end{proposition}
in the category $D^b_c(\mathcal{A}_d,\mathbb{Q}_l)$. Each summand on the right is stable under the Hecke correspondence $\mathcal{H}_d$, which acts on $\mathrm{R}g_{(d_1,d_2),!}L_{(d_1,d_2)}$ by $d_1-d_2+N$.

The prove is the same as Yun-Zhang and Howard-Shnidman. The idea is first compare them on an open dense part of $\mathcal{A}_d$. Let $\hat{Y}^\circ_{3,d'}\subset\hat{Y}_{3,d'}$ be the open dense subscheme classifying divisors on $Y_3$ of degree $d'$ without multiplicity. Similarly define $\hat{Y}^\circ_{d'}\subset\hat{Y}_{d'}$. Then one has the everywhere unramified finite cover of schemes:
\begin{equation}
\text{Nm}: \hat{Y}_{d'}\longrightarrow\hat{Y}_{3,d'}
\end{equation}
with deck transformation group $\Gamma_{d'}$. Then $\text{Nm}_*\mathbb{Q}_l$ is the local system determined by the representation $\text{Ind}^{\Gamma_{d'}}_{S_{d'}}\textbf{1}$. Using the decomposition \eqref{repdecomp}, one gets:
\begin{equation}
\text{Nm}_*\mathbb{Q}_l\xrightarrow{\sim}\bigoplus_{d_1,d'_2\geq 0, d_1+d'_2=d'}V_{(d_1,d'_2)}
\end{equation}
Here right hand side means the local systems determined by the representations $V_{(d_1,d'_2)}$. Define $\mathcal{A}_d^\circ$ to be the open dense part of $\mathcal{A}_d$ such that $\div{\xi}-D_3'$ is multiplicity free, i.e.:
\begin{equation}
\xymatrix{
\mathcal{A}_d^\circ \ar[r] \ar[d] & \mathcal{A}_d \ar[d]\\
\hat{Y}_{3,d'}^\circ \ar[r] & \hat{Y}_{3,d'}
}
\end{equation}
By the proper base change theorem of \'etale cohomology, one has:
\begin{align}
\mathrm{R}f_{d,!}\mathbb{Q}_l|_{\mathcal{A}_d^\circ}&\xrightarrow{\sim}\text{Nm}_!\mathbb{Q}_l|_{\mathcal{A}_d^\circ}\\
\mathrm{R}g_{(d_1,d_2),!}L_{(d_1, d_2)}|_{\mathcal{A}_d^\circ}&\xrightarrow{\sim}V_{(d_1,d_2)}|_{\mathcal{A}_d^\circ}
\end{align}
Now by the finiteness of of $f_d$ and $g_d$, one knows that $\mathrm{R}f_{d,!}\mathbb{Q}_l[d'-g+1]$ is the middle extension of $\mathrm{R}f_{d,!}\mathbb{Q}_l[d'-g+1]|_{\mathcal{A}_d^\circ}$ and $\mathrm{R}g_{(d_1,d_2),!}L_{(d_1, d_2)}[d'-g+1]$ is the middle extension of $\mathrm{R}g_{(d_1,d_2),!}L_{(d_1, d_2)}[d'-g+1]|_{\mathcal{A}_d^\circ}$. Therefore the isomorphism above extends to the entire $\mathcal{A}_d$.

For the Hecke correspondence part, it is easy to check combinatorially that when restricted to $\mathcal{A}_d^\circ$, the action $[\mathcal{H}^d]$ on $V_{(d_1,d'_2)}|_{\mathcal{A}_d}$ coincide with $H$ in \eqref{repdecomp}, hence a scalar $d_1-d'_2=d_1-d_2+N$. Again by the theory of middle extension, this extends to all the $\mathrm{R}g_{(d_1,d_2),!}L_{(d_1, d_2)}$ over the entire $\mathcal{A}_d$.

\subsection{Comparison of the intersection number and the orbital integral}
First let's review some theorems about the Hecke algebra $\mathscr{H}^\Sigma$ in Yun-Zhang.

Define $\widetilde{\mathscr{H}}_l^{\Sigma}$ as the image of the following map:
\begin{equation*}
\mathscr{H}_G^\Sigma\otimes\mathbb{Q}_l\longrightarrow End_{\mathbb{Q}_l}(V)\times End_{\mathbb{Q}_l}(\mathcal{A}(K)\otimes\mathbb{Q}_l)\times\mathbb{Q}_l[Pic_X(k)]^{\iota_Pic}
\end{equation*}

It is a finitely generated algebra over $\mathbb{Q}_l$. Then they proved:
\begin{proposition}
Let $\mathscr{H}^{\Sigma,\dagger}\subset\mathscr{H}_G^\Sigma$. be the linear span of $h_D$ for effective divisors $D$ of degree greater then a fixed positive number, say $M$. Then the map:
\begin{equation*}
\mathscr{H}^{\Sigma,\dagger}\otimes\mathbb{Q}_l\longrightarrow\mathscr{H}^\Sigma_G\otimes\mathbb{Q}_l\longrightarrow\widetilde{\mathscr{H}}_l^\Sigma
\end{equation*}
is surjective. 
\end{proposition}

This lead to:
\begin{proposition}
For any $f\in\mathscr{H}^{\Sigma}_G$, one has:
\begin{equation}\label{fundamental equation}
\mathbb{I}^\mu(f)=(\text{log}q)^{-r}\mathbb{J}_r(f)
\end{equation}
\end{proposition}

First one check the equation for $f=h_D$ for $\text{deg}D\geq 2g_3-1+N$. The left hand side:
\begin{equation*}
\sum_{\xi\in\mathcal{A}_D(k)}\text{Tr}(f_{d,!}[\mathcal{H}_d^\mu]\circ Fr_\xi, (\mathrm{R}f_{d,!}\mathbb{Q}_l)_{\overline{\xi}})
\end{equation*}
From the last part, we have seen $\mathrm{R}f_{d,!}\mathbb{Q}_l\xrightarrow{\sim}\bigoplus_{d_1\geq 0, d_2\geq N, d_1+d_2=2d}\mathrm{R}g_{(d_1,d_2),!}L_{(d_1, d_2)}$. Therefore one can rewrite the above formula as:
\begin{equation}
\sum_{\xi\in\mathcal{A}_D(k)}\sum_{d_1\geq 0, d_2\geq N, d_1+d_2=2d}\text{Tr}(f_{d,!}[\mathcal{H}_d^\mu]\circ Fr_\xi, (\mathrm{R}g_{(d_1,d_2),!}L_{(d_1, d_2)})_{\overline{\xi}})
\end{equation}
Since $\mathcal{H}^\mu_d$ acts on $\mathrm{R}g_{(d_1,d_2),!}L_{(d_1, d_2)}$ by the scalar, one has:
\begin{equation*}
\sum_{\xi\in\mathcal{A}_D(k)}\sum_{d_1\geq 0, d_2\geq N, d_1+d_2=2d}(d_1-d_2+N)^r\text{Tr}(Fr_\xi, (\mathrm{R}g_{(d_1,d_2),!}L_{(d_1, d_2)})_{\overline{\xi}})
\end{equation*}
The right hand side:  
\begin{equation*}
\sum_{\xi\in\mathcal{A}_D(k)}\sum_{(d_1+d_2=2d, d_1\geq 0, d_2\geq N)}q^{(d_1-d_2+N)\cdot s}\text{Tr}(Fr_\xi, (\mathrm{R}g_!L_{(d_1,d_2)})_{\overline{\xi}})
\end{equation*}
Taking the $r$th derivative and evaluate it at $s=0$, one gets:
\begin{equation*}
\sum_{\xi\in\mathcal{A}_D(k)}\sum_{(d_1+d_2=2d, d_1\geq 0, d_2\geq N)}(d_1-d_2+N)^r\text{Tr}(Fr_\xi, (\mathrm{R}g_!L_{(d_1,d_2)})_{\overline{\xi}})
\end{equation*}
It is the same as the left hand side. 

Now for any $f\in\mathscr{H}^\Sigma_G$, Yun and Zhang proved that $\mathbb{I}^\mu(f)$ and $(\text{log}q)^{-r}\mathbb{J}_r(f)$ only depend on $f$'s image in $\widetilde{\mathscr{H}}_l^{\Sigma}$. From the last proposition, one can get the conclusion.

\subsection{Final conclusion} A larger part of Yun-Zhang's papers was devoted to proving a ``coarse" spectral expansion of the cohomology of the moduli of shtukas (integral model, not just the generic fiber). I'm going to restate it and use it.

Let $\mathcal{Y}=\text{Spec}\widetilde{\mathscr{H}}^\Sigma_l$ and $V(\xi)$ be $H_c^r(Sht^{\underline{\mu}}_G(\Sigma;\xi)$ for a $\xi:\mathfrak{S}_\infty\longrightarrow\widetilde{\mathfrak{S}}_\infty$. They proved the following:
\begin{theorem}
\begin{enumerate}
\item There is a decomposition of the scheme $\mathcal{Y}^{red}$:
\begin{equation*}
\mathcal{Y}^{red}=Z_{Eis,\mathbb{Q}_l}\coprod\mathcal{Y}_0
\end{equation*}
where $\mathcal{Y}_0$ is a finite set of closed points. Also there is a corresponding decomposition:
\begin{equation*}
\widetilde{\mathscr{H}}^\Sigma_l=\widetilde{\mathscr{H}}^\Sigma_{l,Eis}\times\widetilde{\mathscr{H}}^\Sigma_{l,0}
\end{equation*}
such that Spec$\widetilde{\mathscr{H}}^\Sigma_{l,Eis}=Z_{Eis,\mathbb{Q}_l}$ and Spec$\widetilde{\mathscr{H}}^\Sigma_{l,0}=\mathcal{Y}_0$.
\item There is a unique decomposition:
\begin{equation*}
V(\xi)=V_0(\xi)\oplus V_{Eis}(\xi)
\end{equation*}
such that Supp$V_{Eis}(\xi)\subset Z_{l,Eis}$ and Supp$V_0(\xi)\subset\mathcal{Y}_0$.  $V_0$ is of finite dimension over $\mathbb{Q}_l$.
\item Base change the scheme and the module over it to $\overline{\mathbb{Q}}_l$. Then the above decomposition of module becomes:
\begin{equation*}
V(\xi)\otimes\overline{\mathbb{Q}}_l=V_{Eis}(\xi)\otimes\overline{\mathbb{Q}}_l\oplus(\oplus_{\mathfrak{m}\in\mathcal{Y}_0(\overline{\mathbb{Q}}_l)}V_\mathfrak{m}(\xi))
\end{equation*}
where $V_{\mathfrak{m}}(\xi)$ is the eigenspace of the character $\mathfrak{m}$ of $\widetilde{\mathscr{H}}^\Sigma_l$.
\end{enumerate}
\end{theorem}

Let $Z_i^\mu$ be the image cohomology class in $V(\xi)$ of $\mathcal{Z}_i^\mu(\xi)$ for $i=1,2$. Further let $Z^\mu_{i,\mathfrak{m}}$ to be the projection of $Z^\mu_i$ to $V_{\mathfrak{m}}(\xi)$.
Identify $h\in\widetilde{\mathscr{H}}^\Sigma_{l,0}$ as $(0,h)$ in $\widetilde{\mathscr{H}}^\Sigma_l$. By the self-adjointness of the $\mathscr{H}_G^\Sigma$ action with respect to the cup product, one has:
\begin{equation}
\mathbb{I}^\mu(h)=\sum_{\mathfrak{m}\in\mathcal{Y}_0(\overline{\mathbb{Q}}_l)}(Z^\mu_{1,\mathfrak{m}}, h*Z^\mu_{2,\mathfrak{m}})
\end{equation}
On the analytic side, looking at the spectral expansion side:
\begin{equation}
\mathbb{J}_r(h)=\sum_{\pi\in\Pi_{\Sigma}(\overline{\mathbb{Q}}_l)}\lambda_\pi(h)(\displaystyle\frac{d}{ds})^r|_{s=0}(q^{Ns}\mathbb{J}_\pi(h^\Sigma_J,s))
\end{equation}
where $\Pi_\Sigma(\overline{\mathbb{Q}}_l)$ is the set of automorphic representations of level $G(\mathbb{O}^\Sigma)\times\prod_{x\in\Sigma}Iw_x$. This set can be viewed as a subset of $\mathcal{Y}_0(\overline{\mathbb{Q}}_l)$: For any $\pi$, take its $G(\mathbb{O}^\Sigma)\times\prod_{x\in\Sigma}Iw_x$ fixed line, then $\widetilde{\mathscr{H}}^\Sigma_{l,0}\otimes\overline{\mathbb{Q}}_l$ acts on this fixed line by the character $\lambda_\pi: \widetilde{\mathscr{H}}^\Sigma_{l,0}\otimes\overline{\mathbb{Q}}_l\longrightarrow\overline{\mathbb{Q}}_l$. 

Now take the test function to be the idempotent $e_\pi$ corresponding to the point $\pi\in\mathcal{Y}_0(\overline{\mathbb{Q}}_i)$, then one has:
\begin{align*}
\mathbb{I}^\mu(e_\pi)&=(Z^\mu_{1,\pi}, Z^\mu_{2,\pi})\\
\mathbb{J}_r(e_\pi)&=(\displaystyle\frac{d}{ds})^r|_{s=0}(q^{Ns}\mathbb{J}_\pi(f,s))
\end{align*}
Applying \eqref{fundamental equation}, one gets:
\begin{equation*}
(Z^\mu_{1,\pi}, Z^\mu_{2,\pi})=(logq)^{-r}(\displaystyle\frac{d}{ds})^r|_{s=0}(q^{Ns}\mathbb{J}_\pi(f,s))
\end{equation*}
In conclusion:
\begin{theorem}
Let $\mathcal{Z}^\mu_1$ and $\mathcal{Z}^\mu_2$ be the Heegner-Drinfeld cycles attached to $Y_1, Y_2$ and $\underline{\mu}$ and $Z^\mu_1, Z^\mu_2$ be their cohomology class in $H^r_c(Sht^r_G(\Sigma)')$. Then one has:
\begin{equation*}
(Z^\mu_{1,\pi}, Z^\mu_{2,\pi})=(logq)^{-r}(\displaystyle\frac{d}{ds})^r|_{s=0}(q^{Ns}\frac{\mathscr{P}_0(\phi,s)\mathscr{P}_{3,\eta}(\overline{\phi_3})}{\langle \phi, \phi \rangle})
\end{equation*}
\end{theorem}

\subsection{Some more conlusions}
In the Yun-Zhang case, the analytic side is just the base change $L$-function times some additional terms related to the genus of the base curve and the level $\Sigma$. Here the period integral is not an $L$-function of any standard form. It is just a period integral. However, using the classical result due to Waldspurger (the function field version was proved by Chuang and Wei, look at their paper \cite{CW}), we still have the following:
\begin{theorem}
Let $\pi$ be antomorphic representation who has spherical level away from $\Sigma$ and is a nonramified twist of Steinberg at each point in $\Sigma$. In the notations as above, one has:
\begin{equation*}
(Z^\mu_{1,\pi}, Z^\mu_{1,\pi})=0
\end{equation*}
is equivalent to:
\begin{equation*}
\mathscr{L}^{(r)}(\pi, 1/2)\cdot L(\pi\otimes\chi_1, 1/2)\cdot L(\pi\otimes\chi_2, 1/2)=0
\end{equation*}
in which $\mathscr{L}(\pi,s+1/2)=q^{(2g-2+N)s}\cdot L(\pi,s+1/2)$
\end{theorem}

To see this, first, according to Yun-Zhang's papers, one has:
\begin{equation*}
\mathscr{P}_0(\phi,s)=q^{(2g-2)s}\cdot L(\pi, s+1/2)
\end{equation*}
For $\phi$ a $G(\mathbb{O}^\Sigma)\times\prod_{x\in\Sigma}Iw_x$ fixed vector in $\pi$.

The other part, $\mathscr{P}_{3,\eta}(\overline{\phi}_3)$, one consider the $T_3(\mathbb{A})\times T_3(\mathbb{A})$-invariant linear functional:
\begin{align*}
P^\pi_3: \pi\times\pi&\longrightarrow\mathbb{C}\\
\phi_1\otimes\phi_2&\longmapsto\mathscr{P}_3(\phi_1,\eta)\mathscr{P}_3(\phi_2,\eta)
\end{align*}
One has the local $T_3(F_x)\times T_3(F_x)$-bilinear form:
\begin{equation*}
\alpha_x(\phi_1\otimes\phi_2)=\displaystyle\frac{L(\chi_{3,x},1)L(\pi_x,ad,1)}{\zeta_{F_x}(2)L(\pi_x, \eta_x, 1/2)}\int_{K_{3,x}^{\times}/F_x^{\times}}(\pi_x(t)\phi_1, \phi_2)_x\eta_x(t)dt
\end{equation*}
and the global-local relation:
\begin{equation*}
P^\pi_3=\displaystyle\frac{\zeta_F(2)L(\pi, \eta, 1/2)}{8L(\chi_3,1)^2L(\pi,ad,1)}\prod_{x\in|X|}\alpha_x
\end{equation*}

Now take $\phi$ as above, and its corresponding form $\phi_3$ on $G(\mathbb{A}_3)$ is a $G_3(\mathbb{O}^\Sigma)\times\prod_{x\Sigma}Iw_{w_x}$ fixed vector. One has $T_3(O_x)\subset G_3(O_x)$ for $x\notin\Sigma$ and $T_3(O_x)\subset Iw_{w_x}$. Therefore one can apply the proposition by Gross-Prasad\cite{GP} to conclude that:
\begin{equation*}
\prod_{x\in|X|}\alpha_x(\phi_x\otimes\overline{\phi}_x)\neq 0
\end{equation*}
Since the term $\displaystyle\frac{\zeta_F(2)}{8L(\chi_3,1)^2L(\pi,ad,1)}$ is a nonvansihing constant, whether $P^\pi_3(\phi\otimes\overline{\phi})$ vanishes or not depends on whether $L(\pi,\eta, 1/2)$ vanishes or not, but because one has:
\begin{equation*}
L(\pi,\eta,s)=L(\pi\otimes\chi_1,s)\cdot L(\pi\otimes\chi_2,s)
\end{equation*}
Therefore $\mathscr{P}_3(\overline{\phi}_3)=0$ is equivalent to $L(\pi\otimes\chi_1,s)\cdot L(\pi\otimes\chi_2,s)=0$. $\square$

Howard and Shnidman also computed the intersection number of the entire $Z^\mu_1$ and $Z^\mu_2$ by computing the orbital integral side. Here one can apply their method directly to get a similar but slightly different result:
\begin{theorem}
No matter how many points of modification you have i.e. for any $r\geq 0$, you'll always get:
\begin{equation*}
(Z^\mu_1,Z^\mu_2)=0
\end{equation*}
\end{theorem}

Take the test function $f$ to be the characteristic function of $G_0(\mathbb{A})$, one gets:
\begin{equation*}
(Z^\mu_1,Z^\mu_2)=(logq)^{-r}(\displaystyle\frac{d}{ds})^r|_{s=0}(q^{Ns}\mathbb{J}(f_J^\Sigma,s))
\end{equation*}
The point is that $\mathbb{J}(f_J^\Sigma,s))$ is always $0$. To get this, let me first cite a lemma proved by Howard and Shnidman:
\begin{lemma}
Fix $\gamma\in J(F)$, let $a\in K_3$ be its invariant. If there exist $t_0\in T_0(\mathbb{A})$ and $t_3\in T_3(\mathbb{A})$ such that $f^\Sigma_J(t^{-1}_0\gamma t_3)\neq 0$, then $a\in k$ and $2a=1$.
\end{lemma}

Notice that in their paper $f$ is the characteristic function of the image of $\widetilde{J}(\mathbb{O})$ in $J(\mathbb{A})$, but here one changes the local components at $\Sigma$ to $1_{Iw_{J,x}}$, but the hypothesis in the lemma still forces the existence of $\phi\in\widetilde{J}(\mathbb{A})$ lifting $t^{-1}_0\gamma t_3$ and send $\mathbb{O}_3$ to $\mathbb{O}\oplus\mathbb{O}$. This is because to be in $Iw_{J,x}$ is more restrictive than to be in $Iso(O_{K_3,x}, O_x\oplus O_x)$. Therefore all of their arguments still apply here.

Now the orbital integral simplies to just one term $\mathbb{J}(\gamma, f^\Sigma_J, s)$, and from their lemma, one can even write down the element $\gamma$: First take an element of $K_3: \epsilon$ such that $K_3=F(\epsilon)$ and $\epsilon$ satisfies an equation of the form $x^2-a=0$ for $a\in F$. Then as an $F$-vector space, $K_3=F+F\epsilon$, and there is an $F$-linear isomorphism: $K_3\longrightarrow F\oplus F: x+y\epsilon\longmapsto(x,y)$. Then this $\phi$ is a representative of $\gamma$ with invariant $1/2$: base change it $K_3$, one gets:
\begin{align*}
K_3\oplus K_3&\longrightarrow K_3\oplus K_3\\
(a,b)&\longmapsto (\displaystyle\frac{a+b}{2}, \displaystyle\frac{a-b}{2\epsilon})
\end{align*}
Then invairant, as we have seen, is given by: $\displaystyle\frac{-1}{4\epsilon}/\displaystyle\frac{-1}{2\epsilon}=\displaystyle\frac{1}{2}$.

Factorize the test function and the character $\eta$: $f^\Sigma_J=\prod_{x\in|X-\Sigma|}f_x\otimes\prod_{x\in\Sigma}1_{Iw_{w_x}}$ and $\eta=\prod_{x\in|X|}\eta_x$. The points away from $\Sigma$ is exactly the same as in Howard-Shnidman's paper. So one only has to care about $x\in\Sigma$. According to our set up, any $x\in\Sigma$ splits in $K_3$. As they pointed out, here one can choose $c\in F_x$ such that $(c\epsilon_x)^2=1$, and define the idempodents:
\begin{equation*}
e=\displaystyle\frac{1+c\epsilon_x}{2},\quad f=\displaystyle\frac{1-c\epsilon_x}{2}
\end{equation*}
These two idempotents actually split $K_{3,x}$ into $K_{3,w_x}\oplus K_{3,w'_x}$, though one doesn't know the order the two over points. The integral over $T_3(F_x)$ is just summing over the set of lattices:
\begin{equation*}
F_x^\times\backslash K_{3,x}^\times/O^\times_{K_3,x}=\{e+\varpi^lf: l\in\mathbb{Z}\}
\end{equation*}
Moding out $O^\times_{K_3,x}$ is because $O^\times_{K_3,x}$ is contained in $Iw_{w'_x}$ or $Iw_{w_x}$. Similarly, the integration over the split torus $T_0(F_x)$ is just summing over $F^\times_x\backslash F_x\times F_x/O^\times_{F,x}\times O^\times_{F,x}$. 

For $1_{Iw_x}((1,\varpi^{-k})\cdot\gamma\cdot (e+\varpi^lf))$ to be not $0$, one must first have $\phi(O^\times_{K_3,x}\cdot(e+\varpi^lf)=(O^\times_{F,x}\times O^\times_{F,x})\cdot(1,\varpi^{-k})$ up to rescaling by $F_x^\times$. From Howard-Shnidman, this is equal to saying: $|\varpi|^k=|c|$ and $l=0$. This is because of the following: by the definition of $\phi$, one has:
\begin{equation*}
e\longmapsto (\displaystyle\frac{1}{2}, \displaystyle\frac{c}{2})\quad \varpi^lf\longmapsto (\displaystyle\frac{1}{2}\varpi^l, -\displaystyle\frac{c}{2}\varpi^l)
\end{equation*}
they should span the same lattice as $(1,\varpi^k)$ up to some rescaling by $F^\times_x$. So there exists a matrix:
\begin{equation*}
\begin{pmatrix} 
a_{11} & a_{12} \\
a_{21} & a_{22}
\end{pmatrix}\in \text{GL}_2(O_x)
\end{equation*}
and a power of the uniformizer $\varpi^\alpha$, such that:
\begin{align*}
\displaystyle\frac{1}{2}\varpi^\alpha=a_{11} &\quad \displaystyle\frac{1}{2}\varpi^\alpha\varpi^l=a_{12}\\
\displaystyle\frac{c}{2}\varpi^\alpha=a_{21}\varpi^k &\quad -\displaystyle\frac{c}{2}\varpi^\alpha\varpi^l=a_{22}\varpi^k
\end{align*}
Since the matrix is in GL$_2(O_x)$, we must have the following: If $ord_x(a_{12}a_{21})>0$, then one must have $ord_x(a_{11})=ord_x(a_{22})=0$. Then from the first equation above, one gets $\alpha=0$. Then from the last equation, one gets $ord_x(c)=k-l$, plugging in the third equation, one gets $ord_x(a_{21})=-l\geq 0$. From the second equation, one gets $ord_x(a_{12})=l\geq 0$. These two inequalities forces $l=0$. This in turn tells us $ord_x(c)=k$. Their claim is correct; 

If $ord_x(a_{11}a_{22})>0$, then $ord_x(a_{12})=ord_x(a_{21})$, one plays a similar game: By the second equation, one gets $\alpha=-l$. Then plug it into the third equation, one sees $ord_x(c)=k+l$. Now by the last equation one gets $ord_x(a_{22})=l\geq 0$ and from the first equation one gets $ord_x(a_{11})=-l$. Again, both $l$ and $-l$ should be non-negative, therefore $l=0$ and $ord_x(c)=k$.

In case all four coefficients has valuation $0$, either of the two arguments above works. This is the proof of their claim. 

In their case this is enough, but now there is one more condition, i.e. $\phi$ should send the sublattice $O_{K_3,x}\cdot(e+\varpi^{l+1})$ or $O_{K_3,x}\cdot(e+\varpi^{l-1})$ to $O_{F,x}\cdot (1,\varpi^{k+1})$, up to rescaling by $\varpi^\alpha$. The same argument as above shows that one must have $l=1$ or $l=-1$ and $ord_x(c)=k+1$, this contradicts what we got above. This implies that for all $t_0$ and $t_3$, $1_{Iw_{w'_x}}(t^{-1}_0\gamma t_3)=0$. The local factors at the level are all $0$, therefore the global integrand is always $0$, the integral just vanish. $\square$



\bibliographystyle{plain}
\bibliography{bib-2}
\end{document}